\documentclass[10pt,a4paper]{amsart} 
\usepackage{amsfonts,amssymb,amsthm,amsmath,ifpdf,stmaryrd,amsrefs,longtable}
\usepackage{enumitem} 
\usepackage[hidelinks]{hyperref}
\usepackage{tikz}
\usetikzlibrary{matrix,arrows}
\usepackage{tikz-cd}

\usepackage[overload]{textcase}

\ifpdf
  \usepackage[final,expansion=alltext,protrusion=alltext]{microtype} 
\fi

\usepackage{xargs} 
\usepackage{mathtools} 
\usepackage{ifthen} 

\theoremstyle{plain}

\newtheorem*{Th*}{Theorem}
\newtheorem*{ThA}{Theorem A}
\newtheorem*{ThB}{Theorem B}
\newtheorem*{ThC}{Theorem C}

\newtheorem*{Cor*}{Corollary}

\theoremstyle{definition}

\theoremstyle{remark}

\makeatletter

\newif\ifTNS 
\TNStrue

\def\printtheoremname#1{\csname#1name\endcsname}
\def\printtheoremnames#1{\csname#1names\endcsname}

\def\thmref#1#2{\printtheoremname{#1}\ifTNS~\fi\ref{#1:#2}}

\def\uc#1#2{\MakeUppercase{#1}{#2}} 

\newcommand{\DefTheorem}[2]{\newenvironmentx{#1}[2][1=\empty,2=\empty]{%
    \ignorespaces%
    \ifx##2\empty%
      \begin{#2}%
    \else%
      \begin{#2}[{\uc##2}]%
    \fi%
    \ifx##1\empty%
      {}%
    \else%
      \label{#1:##1}%
    \fi%
    \ignorespaces}{\end{#2}\ignorespacesafterend}}

\makeatother

\DefTheorem{Th}{theorem}
\DefTheorem{Prop}{proposition}
\DefTheorem{Cor}{corollary}
\DefTheorem{Lem}{lemma}
\DefTheorem{Def}{definition}
\DefTheorem{Rem}{remark}
\DefTheorem{Par}{para}
\DefTheorem{Not}{notation}
\DefTheorem{Exer}{exercise}
\DefTheorem{Ex}{example}
\DefTheorem{Cons}{construction}
\DefTheorem{Scho}{scholium}

\newenvironment{Par*}{\ignorespaces\noindent\ignorespaces}{\ignorespacesafterend}

\numberwithin{equation}{section}

\newcommand\Define[2][\empty]{\ignorespaces%
  \emph{#2}}%

\tikzset{
  commutative diagrams/.cd,
  arrow style=tikz,
  diagrams={>=stealth},
  shift up/.style={
    to path={([yshift=#1]\tikztostart.east) -- ([yshift=#1]\tikztotarget.west) \tikztonodes}},
  shift up left/.style={
    to path={([yshift=#1]\tikztostart.west) -- ([yshift=#1]\tikztotarget.east) \tikztonodes}},
  mathdouble/.style={-,double equal sign distance}
}

\newenvironment{tikzmat}{\ignorespaces%

  \begin{minipage}{0.9\textwidth}\centering %
  \begin{tikzpicture}[>=stealth,description/.style={fill=white,inner sep=2pt},mathmat/.style={matrix of math nodes,row sep=2.5em,column sep=2.3em, text height=1.5ex, text depth=0.25ex},pathmat/.style={->,font=\scriptsize},descr/.style={fill=white,inner sep=2pt},mathdouble/.style={-,double equal sign distance}]%
  \ignorespaces}%
  {\ignorespaces\end{tikzpicture}%
  \end{minipage}\\%
  \ignorespacesafterend}

\def\ger{\mathfrak}

\newcommand\CategoryTypeface{\mathbf}
\def\cat{\CategoryTypeface}

\newcommand\SheafTypeface{\mathcal}
\def\sh{\SheafTypeface}

\def\DMO{\DeclareMathOperator}

\newcommand\ev{{\bar 0}}
\newcommand\odd{{\bar 1}}

\newcommand{\defi}{\coloneqq}     

\newcommand{\fa}{for all }

\newcommand{\fs}{for some }
\newcommand\mathfa[1][{}]{\quad\text{\fa{#1} }}

\newcommand{\scth}{such that }
\newcommand{\AND}{and}

\newcommand\mathtxt[1]{\quad\text{{#1}}\quad}

\newcommand{\nd}{\mathtxt{\AND}}

\newcommand\vq{\emph{v.}~}
\newcommand\eg{\emph{e.g.}~}
\newcommand\ie{\emph{i.e.}~}

\newcommand\via{\emph{via}~}
\newcommand\loccit{\emph{loc.~cit.}}

\newcommand\vs{\emph{vs.}~}

\newcommand\vphi{\varphi}
\newcommand\vrho{\varrho}

\newcommand\vkappa{\varkappa}
\newcommand\eps{\varepsilon}
\newcommand\nats{\mathbb{N}}
\newcommand\ints{\mathbb{Z}}

\newcommand\reals{\mathbb{R}}
\newcommand\cplxs{\mathbb{C}}
\newcommand\knums{\mathbb K}

\def\aff{{\mathbb A}}
\def\sdual{{\mathbb D}}

\newcommand\vvoid{\varnothing}
\newcommand\sle{\leqslant}
\newcommand\sge{\geqslant}

\DMO\dom{\mathrm{dom}}
\DMO\rk{\mathrm{rk}}
\DMO\Ad{\mathrm{Ad}}
\DMO\ad{\mathrm{ad}}
\DMO\GL{\mathrm{GL}}
\DMO\id{\mathrm{id}}
\DMO\pr{\mathrm{pr}}
\DMO\gr{\mathrm{gr}}
\DMO\sll{\ger{sl}}
\DMO\sdim{\mathrm{sdim}}
\DMO\sgn{\mathrm{sgn}}
\DMO\re{\mathrm{Re}}
\DMO\Gal{\mathrm{Gal}}
\DMO\Ann{\mathrm{Ann}}
\DMO\coker{\mathrm{coker}}
\DMO\im{\mathrm{im}}
\DMO\coim{\mathrm{coim}}
\DMO\Spec{\mathrm{Spec}}
\DMO\codim{\mathrm{codim}}
\DMO\chr{\mathrm{char}}
\DMO\supp{\mathrm{supp}}
\DMO\str{\mathrm{str}}
\DMO\tr{\mathrm{tr}}

\DMO\Top{\cat{Top}}
\DMO\Sets{\cat{Sets}}
\DMO\SMan{\cat{SMan}}
\DMO\SRSp{\cat{SRSp}}
\DMO\SSp{\cat{SSp}}
\DMO\SVSp{\cat{SVSp}}
\DMO\SVec{\cat{SVec}}
\DMO\SAlg{\cat{SAlg}}
\DMO\Mod{\cat{Mod}}
\DMO\Op{\cat{Op}}
\DMO\Cov{\cat{Cov}}
\DMO{\Shv}{\cat{Sh}} 

\DMO\Ob{\mathrm{Ob}}

\newcommand{\lBr}{[\kern-.65ex[}
\newcommand{\rBr}{]\kern-.65ex]}

\makeatletter
\newcommand\Size[7][1]{
                                 \ifx#20%
                                        \def\r@l{}\def\r@m{}\def\r@r{}%
                                 \else%
                                    \ifx#21%
                                           \def\r@l{\bigl}\def\r@r{\bigr}\def\r@m{\bigm}%
                                    \else%
                                           \ifx#22%
                                                 \def\r@l{\Bigl}\def\r@r{\Bigr}\def\r@m{\Bigm}%
                                            \else%
                                                 \ifx#23%
                                                        \def\r@l{\biggl}\def\r@r{\biggr}\def\r@m{\biggm}%
                                                  \else
                                                        \ifx#24%
                                                        \def\r@l{\Biggl}\def\r@r{\Biggr}\def\r@m{\Biggm}%
                                                        \fi%
                                                  \fi%
                                            \fi%
                                      \fi%
                                 \fi%
                                 \ifx#10%
                                       \def\r@m{}%
                                 \fi%
                                 \r@l#3{#4}\r@m#5{#6}\r@r#7%
}%

\makeatother

\makeatletter
\def\Set@Scallop[#1]#2#3{{#1}\Parens{#2}{#3}}
\newcommand\DeclareScalableOperator[2]{%
  \expandafter\def\csname#1\endcsname{\@ifnextchar[{{#2}\Set@Scallop}{{#2}\Set@Scallop[{}]}}
}
\makeatother

\def\DSO{\DeclareScalableOperator}

\DSO{Presh}{\cat{Presh}}
\DSO{Sh}{\cat{Sh}}

\DSO{ShHom}{\sh Hom} 
\DSO{ShGHom}{\underline{\sh Hom}} 
\DSO{ShEnd}{\sh End} 
\DSO{ShGEnd}{\underline{\sh End}} 
\DSO{ShDer}{\sh Der} 
\DSO{ShGDer}{\underline{\sh Der}} 
\DSO{ShCliff}{\sh C\text{\emph{liff}}} 

\DSO{Ct}{\mathcal{C}} 
\DSO{Cc}{\mathcal{C}_c} 
\DSO{Lp}{\mathbf{L}} 
\DSO{Hom}{\mathrm{Hom}} 
\DSO{End}{\mathrm{End}} 
\DSO{Aut}{\mathrm{Aut}} 
\DSO{Uenv}{\mathfrak{U}} 
\DSO{Cuenv}{\widehat{\mathfrak{U}}} 
\DSO{ABer}{\Abs0{\mathrm{Ber}}}
\DSO{Proj}{\mathrm{Proj}}
\DSO{Ber}{\mathrm{Ber}}
\DSO{Vol}{\mathrm{Vol}}
\DSO{Form}{\Omega}
\DSO{Ind}{\mathrm{Ind}}
\DSO{Coind}{\mathrm{Coind}}
\DSO{Der}{\mathrm{Der}}
\DSO{GDer}{\underline{\mathrm{Der}}}
\DSO{Alg}{\mathrm{Alg}}
\DSO{GHom}{\underline{\mathrm{Hom}}} 
\DSO{GEnd}{\underline{\mathrm{End}}} 
\DSO{Maps}{\mathrm{Maps}}
\DSO{Cliff}{\mathrm{Cliff}} 
\DSO{Jac}{\mathrm{Jac}} 

\newcommand\Set[3]{
                                 \Size{#1}{\{}{#2}{|}{#3}{\}}%
}%
\newcommand\Dual[3]{
                                 \Size[0]{#1}{\langle}{#2}{,}{#3}{\rangle}%
}%
\newcommand\Parens[2]{
  \Size[0]{#1}{(}{#2}{}{}{)}
}
\newcommand\Bracks[2]{
  \Size[0]{#1}{[}{#2}{}{}{]}
}

\newcommand\Abs[2]{
  \Size[0]{#1}{\lvert}{#2}{}{}{\rvert}
}
\newcommand\Span[2]{
  \Size[0]{#1}{\langle}{#2}{}{}{\rangle}
}
\makeatletter

\makeatletter

\newif\if@smallmat
\newif\if@none
\newif\if@paren
\newif\if@brack
\newif\if@brace
\newif\if@vline
                                 {\ifx#20%
                                        \@smallmattrue%
                                  \else%
                                         \@smallmatfalse
                                  \fi%
                                  \ifx#11%
                                         \@nonefalse\@parentrue\@brackfalse\@bracefalse\@vlinefalse%
                                  \else%
                                       \ifx#12%
                                            \@nonefalse\@parenfalse\@bracktrue\@bracefalse\@vlinefalse%
                                        \else%
                                            \ifx#13%
                                                 \@nonefalse\@parenfalse\@brackfalse\@bracetrue\@vlinefalse%
                                            \else%
                                                 \ifx#14%
                                                       \@nonefalse\@parenfalse\@brackfalse\@bracefalse\@vlinetrue
                                                 \else%
                                                       \ifx#15%
                                                             \@nonefalse\@parenfalse\@brackfalse\@bracefalse\@vlinefalse%
                                                       \else%
                                                             \@nonetrue\@parenfalse\@brackfalse\@bracefalse\@vlinefalse%
                                                       \fi%
                                                 \fi%
                                            \fi%
                                        \fi%
                                   \fi%
                                   \if@smallmat%
                                        \if@none%
                                             \begin{smallmatrix}%
                                        \else%
                                            \if@paren%
                                                  \left(\begin{smallmatrix}%
                                            \else%
                                                  \if@brack%
                                                          \left[\begin{smallmatrix}%
                                                  \else%
                                                          \if@brace%
                                                               \left\{\begin{smallmatrix}%
                                                          \else%
                                                               \if@vline%
                                                                    \left\lvert\begin{smallmatrix}%
                                                                \else%
                                                                    \left\lVert\begin{smallmatrix}%
                                                                \fi%
                                                          \fi%
                                                  \fi%
                                            \fi%
                                        \fi%
                                   \else%
                                        \if@none%
                                             \begin{matrix}%
                                        \else%
                                            \if@paren%
                                                  \begin{pmatrix}%
                                            \else%
                                                  \if@brack%
                                                          \begin{bmatrix}%
                                                  \else%
                                                          \if@brace%
                                                               \begin{Bmatrix}%
                                                          \else%
                                                               \if@vline%
                                                                    \begin{vmatrix}%
                                                                \else%
                                                                    \begin{Vmatrix}%
                                                                \fi%
                                                          \fi%
                                                  \fi%
                                            \fi%
                                        \fi%
                                   \fi}%
                                  {\if@smallmat%
                                        \if@none%
                                             \end{smallmatrix}%
                                        \else%
                                            \if@paren%
                                                  \end{smallmatrix}\right)%
                                            \else%
                                                  \if@brack%
                                                          \end{smallmatrix}\right]%
                                                  \else%
                                                          \if@brace%
                                                               \end{smallmatrix}\right\}%
                                                          \else%
                                                               \if@vline%
                                                                    \end{smallmatrix}\right\rvert%
                                                                \else%
                                                                    \end{smallmatrix}\right\rVert%
                                                                \fi%
                                                          \fi%
                                                  \fi%
                                            \fi%
                                         \fi%
                                   \else%
                                        \if@none%
                                             \end{matrix}%
                                        \else%
                                            \if@paren%
                                                  \end{pmatrix}%
                                            \else%
                                                  \if@brack%
                                                          \end{bmatrix}%
                                                  \else%
                                                          \if@brace%
                                                               \end{Bmatrix}%
                                                          \else%
                                                               \if@vline%
                                                                    \end{vmatrix}%
                                                                \else%
                                                                    \end{Vmatrix}%
                                                                \fi%
                                                          \fi%
                                                  \fi%
                                            \fi%
                                        \fi%
                                   \fi}%
                         
\makeatother

\newcommand\ssplfg[2][\empty]{%
  \smash{\SSp
    \ifx#1\empty%
      ^{\mathrm{lfg}}_{#2}
    \else
      ^{#1,\mathrm{lfg}}_{#2}
    \fi
  }
}

\def\clap#1{\hbox to 0pt{\hss#1\hss}} 

\begin{document}

\title{Singular superspaces}

\author[Alldridge]
{Alexander Alldridge}

\address{Universit\"at zu K\"oln\\
Mathematisches Institut\\
Weyertal 86-90\\
50931 K\"oln\\
Germany}
\email{alldridg@math.uni-koeln.de}

\author[Hilgert]
{Joachim Hilgert}

\address{Universit\"at Paderborn\\
Institut f\"ur Mathematik\\
33095 Paderborn\\
Germany}
\email{hilgert@math.upb.de}

\author[Wurzbacher]
{Tilmann Wurzbacher}

\address{%
Institut \'E.~Cartan\\
Universit\'e de Lorraine et
C.N.R.S.\\
57045 Metz, France}

\curraddr{%
  Fakult\"at f\"ur Mathematik\\
  Ruhr-Universit\"at Bochum\\
  44780 Bochum\\
  \mbox{Germany}
}
\email{tilmann.wurzbacher@univ-lorraine.fr}

\thanks{Research funded by Deutsche Forschungsgemeinschaft (DFG), grant nos.~SFB/TR 12 (all authors), ZI 513/2-1 (A.A.), and the Institutional Strategy of the University of Cologne within the German Excellence Initiative (A.A.). The final publication is available at Springer \emph{via} \texttt{http://dx.doi.org/10.1007/s00209-014-1323-5}}



\begin{abstract}
  We introduce a wide category of superspaces, called locally finitely generated, which properly includes supermanifolds, but enjoys much stronger permanence properties, as are prompted by applications. Namely, it is closed under taking finite fibre products (\ie is finitely complete) and thickenings by spectra of Weil superalgebras. Nevertheless, in this category, morphisms with values in a supermanifold are still given in terms of coordinates. This framework gives a natural notion of relative supermanifolds over a locally finitely generated base. Moreover, the existence of inner homs, whose source is the spectrum of a Weil superalgebra, is established; they are generalisations of the Weil functors defined for smooth manifolds. 
\end{abstract}

\maketitle

\section{Introduction}

Supermanifolds form a natural $\ints/2\ints$-graded generalisation of smooth manifolds. In contrast to the latter, certain basic constructions in supergeometry automatically lead to `singular' superspaces. One is thus compelled to work within a wider category. 

In this article, we consider superspaces $X$ over the field $\knums=\reals,\cplxs$, which locally admit embeddings into (smooth or analytic) affine superspace $\aff^{p|q}$. If these embeddings obey a general condition dubbed \emph{tidiness} (which in fact makes sense over any base field or ring), then we call $X$ \emph{locally finitely generated}. The following are our main results:

\begin{ThA}
  The category $\smash{\SSp_\knums^{\mathrm{lfg}}}$ of locally finitely generated superspaces admits finite fibre products. Supermanifolds form a subcategory closed under finite products. 
\end{ThA}

\begin{ThB}
  If $X$ is locally finitely generated and $Y$ is a superdomain, then there is a natural bijection between morphisms $X\to Y$ and their coordinate expressions. 
\end{ThB}

Theorem B is classical when $X$ is a supermanifold \cite{leites}, but demands substantially new methods in the general case. Within this framework, we define a notion of \emph{relative supermanifolds} over a locally finitely generated base, generalising the notion of \emph{families of supermanifolds}. A generalisation to more general base fields is conceivable, since the concept of tidiness, which is crucial here, does not make any assumptions thereon.

If $A$ is a \emph{Weil superalgebra}, \ie a finite dimensional local supercommutative superalgebra, then $\Spec A\defi(*,A)$ is locally finitely generated, and by Theorem A, $\Spec_SA\defi S\times\Spec A$, the \emph{Weil thickening} of $S\in\smash{\SSp^{\mathrm{lfg}}_\knums}$, exists in $\smash{\SSp_\knums^{\mathrm{lfg}}}$. Our methods readily allow to establish the following:

\begin{ThC}
  For any $\smash{S\in\SSp^{\mathrm{lfg}}_\knums}$, any Weil superalgebra $A$, and any relative supermanifold $X/S$, the relative inner hom functor $\GHom[_S]0{\Spec_SA,X}$ on locally finitely superspaces over $S$ is representable by a relative fibre bundle over $S$.
\end{ThC}

These functors are the super-analogues of the so-called \emph{Weil functors} \cite{weil-pointproches}. Examples include the (relative) even and odd tangent bundle, as well as all higher tangent and jet bundles. Weil functors and relative supermanifolds play a prominent role \eg in the Stolz--Teichner programme on Euclidean field theories, see Ref.~\cite{hkst-eft}.

\smallskip\noindent
Let us explain why supermanifolds do not suffice for applications. Firstly, fibre products of supermanifolds typically do not exist as such. Besides the problems arising from non-transversal intersections, well-known from differential topology, relevant new examples include the isotropies of supergroup actions at odd points.

Furthermore, for applications in \eg Lie theory or variational calculus, it is natural to consider a category stable under Weil thickenings. Finally, relative supermanifolds over a general base are required for a full-blown theory of Berezin integration along the fibres.

In algebraic geometry, the objective of constructing a category that meets these requirements is attained by the introduction of the concept of schemes, including \emph{ab initio} singular, \ie non-smooth, objects. Since their inception in the work of Berezin, Kostant, and Leites, supermanifolds have been defined within the larger category of superspaces, \ie locally super-ringed spaces, and thus, one may attempt to find a larger subcategory in which the above difficulties are no longer present. In this paper, we single out the full subcategory of locally finitely generated superspaces mentioned above. 

Morally, the local finite generation of a superspace is analogous to the condition that a scheme be \emph{locally of finite type}. However, the ring of germs of $\sh C^\infty$ functions on $\reals^n$ is far from being Noetherian, so the study of these spaces is not merely a translation of methods from algebraic geometry, but requires different techniques. (Similar considerations have recently lead Carchedi--Roytenberg \cite{cr} to study supergeometry from the point of view of $\sh C^\infty$-algebras.)

Indeed, for the embeddings in the definition of locally finite generation, we are obliged to require the property of \emph{tidiness}. Technically, this means that their vanishing ideals satisfy the Whitney condition \cite{moerdijk-reyes}, which in the smooth case is equivalent to closedness of the ideal in the $\sh C^\infty$ topology. We do not apply it in this form here, but prefer to use a formulation that makes no assumption on base fields. Without the condition of tidiness, the existence of equalisers would not be given. Moreover, in their construction, tidiness has to be enforced by a universal construction, which we call \emph{tidying}. 

Another respect, in which we have to depart from the paradigm of algebraic geometry, concerns the description of morphisms to affine superspace $\aff^{m|n}$. Morphisms from any locally ringed space to an affine scheme are in one-to-one correspondence with morphisms of the algebras of global sections, a fact which is quickly deduced by localisation. In particular, morphisms to \emph{algebraic} affine space $\aff^m$ are given by coordinates; that is, they are in bijection with tuples of functions $(x_1,\dotsc,x_m)$.

By contrast, to determine to which extent morphisms with range in \emph{smooth} (or analytic) affine superspace $\aff^{m|n}$ can be described by coordinates is a subtle matter. For the case of supermanifolds, the fact that they can is a fundamental result due to Leites \cite{leites}. We call a superspace for which this theorem holds for all open subspaces \emph{Leites regular} and investigate the stability of this property under passage to colimits, locally closed subspaces, and thickenings. An important outcome of this study is that tidiness is crucial to ensure the hereditarity of Leites regularity; in particular, locally finitely generated superspaces are Leites regular. 

An important aspect of Leites regularity is that in the category of Leites regular superspaces, the smooth (or analytic) affine superspace $\aff^m$ is the $m$-fold product of $\aff^1$. This is not true in the larger category of superspaces, for the algebra of smooth functions on $\reals^2$ is not the tensor square of the ring of univariate smooth functions. Thus, unless one is willing to work with sheaves of topological vector spaces, there is no simple-minded construction of the products $\aff^1\times\dotsm\times\aff^1$. Instead, we construct smooth (and analytic) affine superspace $\aff^m$ directly. The general case of $\aff^{m|n}$ appears as a special case of Weil thickening. 

By tidying, the existence of products for affine superspaces carries over to locally finitely generated superspaces. Together with the existence of equalisers (which exist for locally Hausdorff tidy superspaces), we arrive by the statement that fibre products of locally finitely generated superspaces exist (Theorem A).  

\medskip\noindent
We end this introduction by giving a brief synopsis of the article's contents; herein, the reader will find references to its quintessential concepts and results, as alluded to above. In Section \ref{sec:gen-ssp}, we recall standard notions from the theory of superspaces. In Section \ref{sec:loc-emb}, we embark upon a study of the fine structure of embeddings (called immersions in algebraic geometry). In particular, we introduce the notion of \emph{girth} (\thmref{Def}{girth-def}), which quantifies the relative size of an ambient space around an embedding. This is crucial below, to prove that Leites regularity is hereditary under certain embeddings. Important examples of embeddings are furnished by the \emph{Weil thickenings}, which are also introduced in this section. Finally, we discuss the concept of tidiness and show how it can be imposed on an embedding by a process called \emph{tidying} (\thmref{Prop}{tidying}).

In Section \ref{sec:leites}, we define several variants of the standard affine superspace $\aff^{m|n}$: real smooth, real analytic, complex analytic, and real smooth and analytic versions with complexified function sheaves. The latter two serve as local models for J.~Bernstein's category of \emph{cs} manifolds. We then define \emph{Leites regular} superspaces; this framework gives meaning to the concept of \emph{local coordinate systems}. 

In Section \ref{sec:lfg}, we show that Leites regularity is stable under colimits (\thmref{Prop}{indlim-regular}) and investigate its stability with respect to embeddings (\thmref{Prop}{emb-regular}) and thickenings (\thmref{Prop}{reg-thick}), the latter of which Weil thickenings are an example of. We then introduce locally finitely generated superspaces. We show that these are Leites regular (\thmref{Prop}{fg-countht}), admit finite limits (\thmref{Cor}{tidy-finlimits}), and Weil thickenings. In this framework, we introduce \emph{relative supermanifolds} over a base that is locally finitely generated. At this level, base change becomes meaningful, and for a fixed base, the relative category enjoys the same properties as the category of usual supermanifolds. Finally, we introduce (relative) \emph{Weil functors} as inner homs with respect to the spectra of Weil superalgebras, and give sharp representability results in $\SMan_S$ and the category of relative fibre bundles. 

\medskip\noindent
\emph{Acknowledgements.} We extend our thanks to the following institutions for providing their stimulating research environments during the preparation of this article: Mathe\-matisches Forschungsinstitut Oberwolfach, Max-Planck-Institut f\"ur Mathematik Bonn, and Ruhr-Universit\"at Bochum. We thank Torsten Wedhorn and Christoph Schabarum for useful comments on a preliminary version of the text, and the referee for remarks leading to a substantial improvement of the manuscript. 

\section{Generalities on superspaces}\label{sec:gen-ssp}
 
In this section, we collect basic facts and definitions related to superspaces in the sense of Ref.~\cite{manin}. Most of these are more or less straightforward generalisations of the ungraded case, \vq Ref.~\cite{gro-dieu-ega1new}. Accordingly, we will only give full proofs in such cases where we deviate from the standard lore. 

\subsection{Superspaces and their morphisms}

We start with the definition of the category of super-ringed spaces.

\begin{Def}[superringed-def][super-ringed spaces]
  A \Define{super-ringed space} is a pair $X=(X_0,\sh O_X)$, where $X_0$ is a topological space and $\sh O_X$ is a sheaf with values in the symmetric monoidal category of supercommutative unital superrings (\ie $\ints/2\ints$ graded rings) and even, unital ring morphisms. It is called the \Define{structure sheaf} of $X$.
  
  By definition, a \Define[super-ringed space!morphism]{morphism of super-ringed spaces} $\vphi:X\to Y$ is a pair $(\vphi_0,\vphi^\sharp)$ where $\vphi_0:X_0\to Y_0$ is a continuous map and $\vphi^\sharp:\vphi_0^{-1}\sh O_Y\to\sh O_X$ is a morphism of sheaves (of superrings).  These form a category denoted by $\SRSp$.
%
\end{Def}

\begin{Rem}
  In the literature, morphisms $\vphi:X\to Y$ of (super)ringed spaces are usually defined as pairs $(\vphi_0,\vphi^\sharp)$ where $\vphi^\sharp$ is a sheaf morphism $\sh O_Y\to\vphi_{0*}\sh O_X$. Explicitly, this means a collection of ring homomorphisms
  \[
    \vphi^\sharp_U:\sh O_Y(U)\to\sh O_X(\vphi_0^{-1}(U)),
  \]
  for any open subset $U\subseteq Y_0$, that commute with restrictions. 

  This point of view is equivalent to the one advocated here, see \eg Refs.~\cite{bredon}*{I.4, Equation (5)}, \cite{gelman}*{II.16, Proposition 17}, or \cite{iversen}*{II.4, Theorem 4.8}. We will need to use both perspectives, and will do so interchangeably. 
\end{Rem}

\begin{Def}[mor-factor][open subspaces]
  Whenever $X$ is a super-ringed space and $U\subseteq X_0$ is an open subset, let $X|_U=(U,\sh O_X|_U)$. Such a super-ringed space is called an \Define{open subspace} of $X$. It comes with a morphism $j_{X|_U}:X|_U\to X$, defined by $\smash{j_{X|_U}\defi(j_{X|_U,0},j_{X|_U}^\sharp)}$ where $j_{X|_U,0}$ is the inclusion of $U$ in $X_0$, and the sheaf morphism $\smash{j_{X|_U}^\sharp:j_{X|_U,0}^{-1}\sh O_X=\sh O_X|_U\to\sh O_X|_U}$ is the identity; $j_{X|_U}$ is called the \Define[inclusion morphism]{inclusion} of $X|_U$ in $X$. We will systematically identify open subspaces and the underlying open sets. Thus, unions and finite intersections of open subspaces make sense.
\end{Def}

If the sections of the structure sheaf are to represent `superfunctions' and non-vanishing `superfunctions' are to be invertible, one needs the stalks of the structure sheaf to be local rings. Ringed spaces with this property are called \emph{locally ringed spaces} in algebraic geometry. Recall that a ring $R$ is called \Define[ring!local]{local} if it possesses a unique maximal ideal. Similarly, one says that a {superring} $R$ is \Define[superring!local]{local}, if it possesses a unique maximal graded ideal. It turns out that in the supercommutative case, one does not have to distinguish between maximal ideals and maximal graded ideals.

\begin{Prop}[super-local][local superrings]
  Let $R$ be a supercommutative superring. The following are equivalent:
  \begin{enumerate}
    \item $R$ is local as a superring;
    \item $R$ is local as an ungraded ring;
    \item for homogeneous $a,b\in R$ with $a+b=1$, one of $a$ and $b$ is invertible; and
    \item for any $a,b\in R$ \scth $a+b=1$, one of $a$ and $b$ is invertible.
  \end{enumerate}
  
  In this case, the unique maximal ideal is graded and consists of those elements of $R$ which are not invertible.
\end{Prop}

In view of \thmref{Prop}{super-local}, the following definition of the category of superspaces is an extension of the definition of locally ringed spaces.

\begin{Def}[superspace-def][superspaces]
  A super-ringed space $X=(X_0,\sh O_X)$ is a \Define{superspace} if for each $x\in X_0$, the stalk $\sh O_{X,x}$ is a local superring. We denote the maximal ideal by $\ger m_{X,x}$. A morphism $\vphi:X\to Y$ of super-ringed spaces where $X$ and $Y$ are superspaces is called \Define[local morphism]{local} if $\vphi^\sharp(\ger m_{Y,\vphi_0(x)})\subseteq\ger m_{X,x}$. The category of superspaces and local morphisms is denoted by $\SSp$. The locally defined sections of $\sh O_X$, where $X$ is a superspace, are called \Define[superfunction]{superfunctions}.
  
  Given a field $k$, a superspace $X$ is called a \Define{$k$-superspace} if $\sh O_X$ is a sheaf of $k$-algebras. The category of $k$-superspaces and local morphisms $\vphi$ \scth $\vphi^\sharp$ is $k$-linear is denoted by $\SSp_k$.
\end{Def}

There is also a relative version of the category of superspaces.

\begin{Def}[rel-def][relative superspaces]
  Let $p:X\to S$ be a morphism of superspaces. We will say that $X$ is a \Define[superspace!over $S$]{superspace over $S$} (or a \Define[relative superspace|see{superspace}]{relative superspace}) and write $X/S$, the morphism $p=p_X$ then being understood.
  
  A \Define[morphism!of superspaces over $S$]{morphism over $S$}, written $f:X/S\to Y/S$, is a morphism $f:X\to Y$ \scth $p_Y\circ f=p_X$, \ie the following diagram commutes:

\begin{center}
  \begin{tikzcd}
    X\arrow[swap]{rd}{p_X}\arrow{rr}{f}&&Y\arrow{dl}{p_Y}\\
    &S 
  \end{tikzcd}
\end{center}

We denote the set of morphisms $X\to Y$ over $S$ by $\Hom[_S]0{X,Y}$ and the category of superspaces over a superspace $S$ by $\SSp_S$. 
\end{Def}

More generally, given any morphism $h:S\to T$, a morphism $f:X\to Y$, where $X/S$ and $Y/T$ are relative superspaces, is \Define[morphism!over $h$]{over $h$} if $p_Y\circ f=h\circ p_X$. This enables us to consider relative superspaces over varying bases.

\begin{Rem}[terminal][terminal objects]
  The category of $k$-superspaces has the terminal object $\Spec k\defi(*,k)$. Any $k$-superspace may therefore be considered as a superspace over $(*,k)$. In fact, the following holds.
\end{Rem}

\begin{Prop}[speck-spaces]
  Let $k$ be a (purely even) field. Given a superspace $X$, there is a bijection between the $k$-superspace structures on $X$ and morphisms $X\to\Spec k$. Moreover, if $X$ and $Y$ are $k$-superspaces, and $\vphi:X\to Y$ is a morphism of superspaces, then $\vphi$ is a morphism of $k$-superspaces if and only if it is over $\Spec k$.
\end{Prop}

We end this subsection with a simple but fundamental construction.

\begin{Cons}[basefieldch][change of base field]
  Let $\ell$ be a field extension of the field $k$. Any $\ell$-superspace $X$ can be naturally considered as a $k$-superspace by forgetting the $\ell$-structure on $\sh O_X$. In other words, by composing the given morphism $X\to\Spec\ell$ of superspaces with the natural morphism $\Spec\ell\to\Spec k$. We will denote the $k$-superspace associated with $X$ by the same letter. 

  Conversely, let $X$ be a $k$-superspace. We define $X_\ell$ to be the $\ell$-superspace 
  \[
  X_\ell\defi (X_0,\sh O_X\otimes_k\ell).
  \] 
  By considering the sheaf embedding $\sh O_X\to\sh O_{X_\ell}$ given by construction, we obtain a morphism $X_\ell\to X$ of $k$-superspaces. Moreover, there is a natural morphism $X_\ell\to\Spec\ell$ of $k$-superspaces. Let $Y$ be a $k$-superspace, and assume given morphisms $\vphi:Y\to X$ and $Y\to\Spec\ell$ of $k$-superspaces. The latter means that $Y$ has an $\ell$-superspace structure compatible with the embedding $k\subseteq\ell$. Thus, $\vphi^\sharp$ extends to a unique $\ell$-linear map $\sh O_X\otimes_k\ell\to\vphi_{0,*}\sh O_Y$. This shows that $X_\ell$ has the universal property of the fibre product $X\times_{\Spec k}\Spec\ell$ in $\SSp_k$. Thus, the functor $(-)_\ell$ is right adjoint to the forgetful functor $\SSp_\ell\to\SSp_k$. We denote the application of the functor $(-)_\ell$ to a morphism $\vphi$ by $\vphi_\ell$.
\end{Cons}

\subsection{Embeddings}

In this subsection, we discuss a particular class of morphisms, which we call embeddings, following the standard terminology of differential geometry. In algebraic geometry, they are called \emph{immersions}. However, the immersions of differential geometry in general are neither injective, nor is their image locally closed. 
To avoid confusions in the applications to supermanifolds, we therefore stick to the wording common in the study of smooth manifolds.

Recall that a continuous map is called an \Define{embedding} if it induces a homeomorphism onto its image (endowed with the relative topology). An embedding is called \Define[embedding!open]{open} resp.~\Define[embedding!closed]{closed} if its image is, in addition, open resp.~closed. An open (resp.~closed) embedding is indeed an open (resp.~closed) map.

\begin{Def}[embedding][embeddings]
  Let $\vphi:Y\to X$ be a morphism of superspaces. It is a called an \Define[embedding!open]{open embedding} if it factors as $\vphi=j_{X|_U}\circ\psi$ where $U\subseteq X_0$ is an open subset, and $\psi:Y\to X|_U$ is an isomorphism of superspaces. In this case, $U=\vphi_0(Y_0)$, and we denote by $\vphi(Y)$ the open subspace $X|_U\subseteq X$.

  The morphism $\vphi$ is called a \Define[embedding!closed]{closed embedding} if $\vphi_0$ is a closed embedding, and the sheaf morphism $\vphi^\sharp:\sh O_X\to\vphi_{0,*}\sh O_Y$ is {surjective}, \ie it induces an isomorphism
    \[
      \sh O_X/\sh I_Y\to\vphi_{0,*}\sh O_Y,
    \]
    where $\sh I_Y\defi\ker\vphi^\sharp$ is by definition the \Define{vanishing ideal} of $Y$ (or of $\vphi$).
  
  Moreover, $\vphi$ is called an \Define{embedding} if it factors as $\vphi=\vphi''\circ\vphi'$ where $\vphi'$ is a closed embedding and $\vphi''$ is an open embedding. Finally, an embedding $\vphi$ is called a \Define{thickening} if $X_0=Y_0$ as topological spaces and $\vphi_0$ is the identity map.
\end{Def}

Open and closed embeddings are embeddings, and embeddings are preserved under field extension. Compositions of open (resp.~closed) embeddings are open (resp.~closed) embeddings. General embeddings are stable under composition, as follows for instance from the following proposition.

\begin{Prop}[open-emb-char]
  Let $\vphi:Y\to X$ be a morphism of superspaces. The following assertions are equivalent:
  \begin{enumerate}[wide]
    \item The morphism $\vphi$ is an embedding.
    \item The map $\vphi_0$ is an embedding, and there exists an open subset $U\subseteq X_0$ \scth $\vphi_0(Y_0)$ is closed in $U$ and $\vphi^\sharp:\sh O_X|_U\to\vphi_{0,*}\sh O_Y$ is surjective.
    \item The map $\vphi_0$ is an embedding, $\vphi_0(Y_0)$ is locally closed, and the sheaf map $\vphi^\sharp:\sh O_X|_U\to\vphi_{0,*}\sh O_Y$ is surjective, where $U\defi\vphi_0(Y_0)\cup X_0\setminus\overline{\vphi_0(Y_0)}$.
  \end{enumerate}
\end{Prop}

The proposition is implied by the following basic topological fact.

\begin{Lem}[open-nbh-locclosed]
  Let $T$ be a topological space and $S$ a locally closed subset. The largest open subset $U\subseteq T$ such that $S$ is closed in $U$ is $U\defi S\cup T\setminus\overline S$.
\end{Lem}

The stability of embeddings under composition admits a partial converse.

\begin{Lem}[emb-factor]
  Let $\vphi:X\to Y$, $\psi:X\to X'$, $\psi':Y\to Y'$, and $\vphi':X'\to Y'$ be given \scth $\psi'\circ\vphi=\vphi'\circ\psi$. If $\psi'_0$ is injective and $\vphi,\psi$ are (open resp.~closed) embeddings, then so is $\vphi$.
\end{Lem}

The characterisation of embeddings in \thmref{Prop}{open-emb-char} also allows us to decide when a morphism factors through  an embedding.

\begin{Prop}[emb-factor][factorisation through embeddings]
  Let $\vphi:Y\to X$ be an embedding with ideal $\sh I_Y\subseteq\sh O_X|_U$, where $U\defi\vphi_0(Y_0)\cup X_0\setminus\overline{\vphi_0(Y_0)}$, and $\psi:Z\to X$ a morphism. Then $\psi$ factors through $\vphi$ if and only if $\psi_0(Z_0)\subseteq\vphi_0(Y_0)$ and $\sh I_Y\subseteq\ker\psi^\sharp|_U$. The factorisation is unique.
\end{Prop}

If the map underlying an embedding is an inclusion we talk about a subspace.

\begin{Def}[closed-subsp][subspaces and closed subspaces]
  Let $X$ be a superspace. If $j:Y\to X$ is a (closed) embedding of superspaces where $Y_0\subseteq X_0$ and $j_0$ is the inclusion of $Y_0$ in $X_0$, then we say that $Y$ is a (closed) \Define{subspace} of $X$.

Given a locally closed subset $Y_0$ of $X_0$, the structure of a subspace $Y$ over $Y_0$ (if it exists) is determined uniquely up to canonical isomorphism by a graded ideal $\sh I_Y$ of $\sh O_X|U$, where $U\defi Y_0\cup X_0\setminus\overline{Y_0}$, \scth $\supp\sh O_X|_U/\sh I_Y=Y_0$. This ideal is called the \Define{vanishing ideal} of $Y$. The isomorphism classes of such subspaces $Y$ are parametrised by the ideals $\sh I_Y$.
\end{Def}

For any subspace $Y$ of a superspace $X$, there is an increasing sequence of subspaces all with underlying space $Y_0$. The sections of the corresponding structure sheaves should be viewed as superfunctions on $X$ that do not necessarily \emph{vanish} on $Y_0$, but are \emph{nilpotent} of a certain order on $Y_0$. 

\begin{Def}[infin-nbh][infinitesimal normal neighbourhoods]
  Let $j:Y\to X$ be an embedding of superspaces with ideal $\sh I_Y\subseteq\sh O_X|_U$, where $\smash{U\defi j_0(Y_0)\cup X_0\setminus\overline{j_0(Y_0)}}$. For any integer $n\sge0$, let $\sh I_{\smash{Y^{(n)}}}\defi\sh I_Y^{n+1}\subseteq\sh O_X|_U$ and $\sh O_{\smash{Y^{(n)}}}\defi j_0^{-1}\Parens0{\sh O_X|_U/\sh I_{\smash{Y^{(n)}}}}$.
  
  The superspace $Y^{(n)}\defi(Y_0,\sh O_{\smash{Y^{(n)}}})$ is called the \Define[embedding!infinitesimal normal neighbourhood (order $n$)]{$n$th infinitesimal normal neighbourhood} of $Y$ in $X$. It comes with a natural embedding $\smash{j^{(n)}:Y^{(n)}\to X}$. In particular, $Y=Y^{(0)}$ and $j=j^{(0)}$. Whenever $n\sle m$, the embedding $\smash{j^{(n)}:Y^{(n)}\to X}$ factors uniquely through $\smash{j^{(m)}:Y^{(m)}\to X}$ \via a closed embedding $\smash{j^{(nm)}:Y^{(n)}\to Y^{(m)}}$, by virtue of \thmref{Prop}{emb-factor}.
\end{Def}

The $(Y^{(n)},j^{(nm)})$ form an inductive system of superspaces, so naturally the question arises whether we can form inductive limits in $\SRSp_S$ and $\SSp_S$. It turns out that not only inductive limits but \emph{all} small colimits exist in these categories.

\begin{Prop}[local-cocomplete]
  Let $S$ be a super-ringed space. The category $\SRSp_S$ is cocomplete, \ie all {small colimits} exist, and $\SSp_S$ is closed under all of them.
\end{Prop}

\begin{proof}
  By \thmref{Prop}{super-local}, superrings are local if and only if they are local as ungraded rings, so the statement follows from the classical case \cite{demazure-gabriel}*{I.1.6}.
 \end{proof}

We can identify the structure sheaves of small colimits of superspaces. 

\begin{Rem}[indlim-functions]
  Let $\smash{X=\varinjlim_kX_k}$ be an inductive limit in $\SRSp_S$ (or in $\SSp_S$). Then $\smash{\sh O_X=\varprojlim_k(f_k)_{0,*}\sh O_{X_k}}$ in $\Sh0{X_0}$, where $f_k:X_k\to X$ are the natural morphisms. Because this projective limit is the projective limit of presheaves and the Yoneda embedding commutes with projective limits, we have that $\smash{\Gamma(\sh O_X)=\varprojlim_k\Gamma(\sh O_{X_k})}$.
  
  The statement also holds with `small colimit' in place of `inductive limit' if one replaces  `projective limit' by `small limit'.
\end{Rem}

\begin{Ex}
  We return to the inductive system of infinitesimal normal neighbourhoods. Let $\smash{Y^{(\infty)}}$ with the morphisms $\smash{j^{(n\infty)}:Y^{(n)}\to Y^{(\infty)}}$ denote the inductive limit of the inductive system $\smash{(Y^{(n)},j^{(nm)})}$ of superspaces. The superspace $\smash{Y^{(\infty)}}$ is called the \Define[embedding!infinitesimal normal neighbourhood (order $\infty$)]{infinitesimal normal neighbourhood of order $\infty$} of $Y$ in $X$. The morphisms $\smash{j^{(nm)}}$, where $m\in\nats\cup\infty$, are thickenings.
\end{Ex}

\begin{Ex}[weil-indlim]
  Let $A_n$ be the local $k$-algebra $k[T]/(T^{n+1})$, where $T$ is an even indeterminate. Since $(T^{n+1})\subseteq(T^n)$, these algebras form a projective system 
  \[
    A_{n+1}\to A_n:f+(T^{n+1})\mapsto f+(T^n).
  \]
  We have a corresponding inductive system of $k$-superspaces 
  \[
    X_n\defi(*,A_n)\to X_{n+1}=(*,A_{n+1}).
  \]
  Here, by $*$, we denote the the singleton topological space. The inductive limit of these is the superspace $X=(*,A)$, where $A\defi k\llbracket T\rrbracket$ is the ring of formal power series in $T$.
\end{Ex}

\subsection{Gluing}

The gluing construction is standard in geometry. It depends on the possibility to have fibre products with open embeddings. 

\begin{Prop}[open-fibreprod][fibre products with open embeddings]
  Let $\vphi:X\to Z$ be an open embedding of superspaces, and $\psi:Y\to Z$ a morphism of superspaces. Then {$X\times_ZY$} exists in $\SSp$, and $p_2:X\times_ZY\to Y$ is an open embedding. The statement carries over to $k$-superspaces.
\end{Prop}

Now we can show how to glue superspaces and their morphisms; this technique will be fundamental. To have it available in its most general form, we introduce the concept of `gluing data'.

\begin{Def}[gluing-data][coverings and gluing data]
  Let $S$ be a superspace. An \Define[superspace!open cover]{open cover of $S$} is a family $(U_i)$ of open subspaces \scth $\bigcup_iU_{i,0}=S_0$.
  
  More generally, a morphism $\vphi:T\to S$ is called a \Define[superspace!covering]{covering of $S$} if there exists an open cover $(U_i)$ of $T$ \scth $T=\coprod_iU_i$, $\vphi|_{U_i}:U_i\to S$ are open embeddings, and the open subspaces $\vphi(U_i)\subseteq S$ form an open cover of $S$. In other words, $\smash{\vphi=\coprod_i\vphi_i}$ where $\vphi_i:U_i\to S$ are open embeddings. We identify $\vphi$ with the collection $(\vphi_i)$.
  
  Given two coverings $\vphi=(\vphi_i:U_i\to S)$ and $\psi=(\psi_j:V_j\to S)$ of $S$, the fibre product $\vphi\times_S\psi=(\vphi_i\times_S\psi_j)$ exists, in view of \thmref{Prop}{open-fibreprod}. In case $\psi=\vphi$, we write $\vphi_{ij}:U_{ij}\defi U_i\times_SU_j\to S$ for this fibre product. The triple intersections $U_i\times_SU_j\times_SU_k\to S$ will be denoted by $\vphi_{ijk}:U_{ijk}\to S$.
  
  Let $\vphi=(\vphi_i:U_i\to S)$ be a covering of $S$, and assume given superspaces $X_i/U_i$ and isomorphisms $\psi_{ij}:X_j\times_{U_j}U_{ij}\to X_i\times_{U_i}U_{ij}$. The collection of data $(\vphi_i:U_i\to S,X_i/U_i,\psi_{ij})$ will be called \Define[superspace!gluing data]{gluing data for a superspace over $S$} if
  \begin{equation}\label{eq:gluing}
    p_{12}^*(\psi_{ij})\circ p_{23}^*(\psi_{jk})\circ p_{13}^*(\psi_{ki})={\id_{X_i\times_{U_i}U_{ijk}}}.
  \end{equation}
  Here, $p_{12}^*{\psi_{ij}}$ is the morphism $X_j\times_{U_j}U_{ijk}\to X_i\times_{U_i}U_{ijk}$ induced by $\psi_{ij}$ by pulling back along the left-hand vertical face of the pullback square
  \begin{tikzmat}
    \matrix (m) [mathmat] {U_{ijk} & U_k\\ U_{ij} & S\\ };
    \path [pathmat] (m-1-1) edge node[auto] {$p_3$} (m-1-2)
              edge node [left] {$p_{12}$} (m-2-1)
          (m-1-2) edge node [auto] {$\vphi_k$} (m-2-2)
          (m-2-1) edge node [below] {$\vphi_{ij}$} (m-2-2);
  \end{tikzmat}
  Similarly for the other quantites in Equation~\eqref{eq:gluing}.
\end{Def}

When the $U_i$ are simply open subspaces of $S$ and the $\vphi_i$ are just the canonical embeddings $j_{U_i}:U_i\to S$, we will forgo mentioning $\vphi$ explicitly and simply write $(U_i,X_i/U_i,\psi_{ij})$ for the gluing data. In this case, somewhat abusing notation, we will write Equation~\eqref{eq:gluing} in the simpler form
  \[
    \psi_{ij}\circ\psi_{jk}\circ\psi_{ki}={\id}\mathtxt{on}X_i\times_{U_i}U_{ijk}.
  \]

Once one has a set of gluing data, it follows from the corresponding statements for sheaves that they lead to glued superspaces.

\begin{Prop}[relative-glue]
  Let $S$ be a superspace and $(\vphi_i:U_i\to S,X_i/U_i,\psi_{ij})$ gluing data for a superspace over $S$. Then there exist a superspace $X/S$ and isomorphisms $\psi_i:X_i/U_i\to (X\times_SU_i)/U_i$ with $\psi_j=\psi_i\circ\psi_{ij}$ on $X_j\times_{U_j}U_{ij}$.
  
  Moreover, these data are uniquely characterised up to unique isomorphism by the following universal property: For any superspace $Y/S$ and any choice of morphisms $\vrho_i:X_i/U_i\to (Y\times_SU_i)/U_i$ \scth $\vrho_j=\vrho_i\circ\psi_{ij}$ on $X_j\times_{U_j}U_{ij}$, there exists a unique morphism $\vrho:X/S\to Y/S$ \scth $\vrho\circ\psi_i=\vrho_i$.
\end{Prop}

The statement of \thmref{Prop}{relative-glue} immediately carries over to $k$-superspaces.

\section{The local structure of embeddings}\label{sec:loc-emb}

In this section, we introduce some concepts that will help us in the quantitative study of embeddings. We also define for any superspace certain distinguished subspaces (namely, the reduction and the body), as well as a useful class of thickenings (namely, the Weil thickenings). We analyse the fine structure of embeddings by the use of the new concept of tidiness. Finally, we show how to improve embeddings by a regularisation process called tidying.

\subsection{Girth and retractions of embeddings} 

There are thickenings of various sizes, as is illustrated by the example of infinitesimal normal neighbourhoods. We introduce the concept of girth of an embedding to quantify this phenomenon. Conversely, one can consider the question whether a thickening can be reversed. This leads to the notion of a retraction.

\begin{Def}[girth-def][finite and countable girth]
  Let $j_Y:Y\to X$ be an embedding of $S$-superspaces. We say that $X$ (or $j_Y$) has \Define[finite girth!around $Y$]{finite girth around $Y$} if the vanishing ideal $\sh I_Y$ is nilpotent, \ie if there exists some non-negative integer $q$ \scth $\smash{\sh I_Y^{q+1}=0}$. In this case, the minimal such $q$ will be called the \Define{girth} of $X$ around $Y$. We shall say that $X$ (or $j_Y$) is \Define[finite girth!local, around $Y$]{locally of finite girth around $Y$} if there exists an open cover of $U$ (the largest open subset of $X_0$ in which $j_{Y,0}(Y_0)$ is closed) \scth the restriction of $\sh I_Y$ to every patch of the cover is nilpotent.
  
  Similarly, we say that $X$ (or $j_Y$) has \Define[countable girth!around $Y$]{countable girth around $Y$} if $\smash{\bigcap_{q=1}^\infty\sh I_Y^q=0}$. If $X$ has locally finite girth around $Y$, then it has countable girth around $Y$.
\end{Def}

Given an embedding $j:Y\to X$, we obtain for the infinitesimal normal neighbourhoods $Y^{(k)}$: If $m\sge n\sge0$ is finite, then the thickening $j^{(nm)}:Y^{(n)}\to Y^{(m)}$ has finite girth at most $\lceil m/n\rceil$.

\begin{Ex}[grass-girth]
  Let $k$ be a field of characteristic $\neq2$ and $A$ the Grassmann algebra on $n$ generators $\tau_1,\dotsc,\tau_n$, \ie
  \[
    A\defi k\Span1{\tau_1,\dotsc,\tau_n}/\Parens1{\tau_i\tau_j+\tau_i\tau_j\bigm|1\sle i\sle j\sle n}.
  \]
  Let $X$ be the superspace $(*,A)$ (where $*$ the singleton space), later to be denoted by $\aff^{0|n}$. Then the $k$-point $X_0\defi(*,k)$ is embedded in $X$ \via 
  \[
    j:X_0\to X,\quad j_0\defi\id_*,\quad j^\sharp(\tau_j)\defi0,\quad j=1,\dotsc,n.
  \]
  It is easy to check that $j$ is in fact the unique morphism $X_0\to X$ of $k$-superspaces.

  The vanishing ideal $J\subseteq A$ of $j$ is generated by $\tau_j$, $j=1,\dotsc,n$. Hence, $J^{n+1}=0$ and $J^\ell\neq0$ for $\ell\sle n$. The girth of $X$ around $X_0$ is $n$.
\end{Ex}

\begin{Ex}
  Let $X_n$ be the $k$-superspace constructed in \thmref{Ex}{grass-girth}. There is an embedding $j^{(n)}:X_n\to X_{n+1}$, given by 
  \[
    j^{(n)\sharp}(\tau_j)\defi
    \begin{cases}
      \tau_j & j\sle n,\\
      0 & j=n+1.
    \end{cases}
  \]
  Let $X$ be the inductive limit of this family of superspaces, which exists because of \thmref{Prop}{local-cocomplete}. Consider the morphism $j:X_0=(*,k)\to X$, where $*$ is the singleton space, induced by the inductive family $j_n\defi j^{(n-1)}\circ\dotsm\circ j^{(0)}:X_0\to X_n$ of morphisms.

  Then $X$ does not have finite girth around $X_0$, since all the finite products $\tau_1\dotsm\tau_n$ are non-zero in $\sh O_X$. However, it is easy to see that it has countable girth. 
\end{Ex}

\begin{Ex}[formal-power]
  Fix $\knums\in\{\reals,\cplxs\}$, and let $X$ be $(\reals,\sh O_X)$, where $\sh O_X$ is the sheaf of $\knums$-valued smooth functions, or the sheaf of $\knums$-valued real analytic functions.

  Let $A_n$ and $A$ be the local algebras defined in \thmref{Ex}{weil-indlim} (for $k=\knums$) and $Y_n\defi(*,A_n)$, $Y\defi(*,A)$ the corresponding $\knums$-superspaces, where $*$ is the singleton space. We define for every $n$ embeddings $j_n:Y_n\to X$ by $j_{n,0}(*)\defi 0$ and 
  \[
    j_n^\sharp(f)\defi\sum_{\ell=0}^n\frac1{\ell!}f^{(\ell)}(0)T^\ell.
  \]
  Then the morphisms $j_n$ are compatible with the inductive system $Y_n\to Y_{n+1}$, and thus define an embedding $j:Y\to X$. Explicitly, we have
  \[
    j^\sharp(f)=\sum_{\ell=0}^\infty\frac1{\ell!}f^{(\ell)}(0)T^\ell, 
  \]
  where the infinite sum is a formal power series.

  The stalk at $0$ of the ideal $\sh J$ of $j$ is the set of all $f\in\sh O_{X,0}$ \emph{flat} at $0$, \ie $f^{(\ell)}(0)=0$ for all $\ell$. In the case of analytic functions, the ideal $\sh J_0$ is zero, so that $j$ has girth zero. By contrast, each of the $j_n$ is of countable, but not of finite girth.

  Next, consider the case of smooth functions. One has $\sh J_0=\bigcap_{k=0}^\infty\ger m^k$ where $\ger m\defi\ger m_{X,0}$. For $R\defi\sh O_{X,0}$, we have that $R/\sh J_0$ is Noetherian, so we find that $\sh J_0\ger m^k=\bigcap_N\ger m^{k+N}=\sh J_0$ and $\sh J_0^2=\bigcap_k\sh J_0\ger m^k=\sh J_0$ by Krull's Theorem (\thmref{Prop}{krull}). Thus, $\bigcap_k\sh J_0^k=\sh J_0$. But $\sh J_0\neq0$. For instance, consider 
  \[
    f:\reals\to\reals,\quad f(x)\defi
    \begin{cases}
      e^{-1/x^2}&x>0,\\
      0&x\sle0.
    \end{cases}
  \]
  Then $f^{(n)}(0)=0$ for all $n\in\nats$, so that the germ $f_0$ of $f$ at $0$ is contained in $\sh J_0$, by \thmref{Cor}{hadamard-higherorder} below. On the other hand, since for every real number $\eps>0$, we have $f(\eps)>0$, we conclude that $f_0$ is non-vanishing. We conclude that $j$ does not have countable girth. In particular, none of the $j_n$ has countable girth.
\end{Ex}

We have the following useful fact, of which we omit the simple proof. 

\begin{Prop}[finht-trans][transitivity of girth]
  Let $j_Z:Z\to Y$ and $j_Y:Y\to X$ be embeddings of $S$-superspaces. If $X$ has (locally) finite girth around $Y$, and $Y$ has (locally) finite girth around $Z$, then $X$ has (locally) finite girth around $Z$ with respect to the embedding $j_Y\circ j_Z:Z\to X$.
\end{Prop}

We now introduce retractions, \ie left inverses for embeddings.

\begin{Def}[retraction-def][retractions]
  Let $X$ be an $S$-superspace and $j_Y:Y\to X$ a morphism over $S$ which is an embedding. Let $U\subseteq X_0$ be the largest open subset in which $j_{Y,0}(Y_0)$ is closed. A \Define{retraction around $Y$} is a left inverse $r:X|_U\to Y$ of $j_Y$, \ie a morphism over $S$ \scth the following diagram commutes:
  \begin{tikzmat}
    \matrix (m) [mathmat] {& X|_U & \\ Y & & Y\\ };
    \path [pathmat] (m-2-1)   edge node [auto] {$j_Y$} (m-1-2)
              edge [mathdouble] (m-2-3)
          (m-1-2) edge node [auto] {$r$} (m-2-3);
  \end{tikzmat}
  If $j_Y$ admits a retraction, it will be called \Define{retractable}. More generally, if $j_Y$ admits an open cover by retractable embeddings $j_{Y_i}$, \ie there are covers by open subspaces $(Y_i)$ of $Y$ and $(X_i)$ of $X$ \scth the embeddings $j_{Y_i}:Y_i\to X_i$ induced by $j_Y$ are retractable, it will be called \Define[retractable!locally]{locally retractable}. We shall then also say that $X$ is (locally) retractable around $Y$.
\end{Def}

  Let $j_Y:Y\to X$ be an embedding and $U$ the largest open subset of $X_0$ \scth $j_{Y,0}(Y_0)$ is closed in $U$. Then $j_Y$ is retractable if and only if the short exact sequence
  \begin{tikzmat}
    \matrix (m) [mathmat] {0 & \sh I_Y & \sh O_X|_U & j_{Y_0,*}\sh O_Y & 0\\ };
    \path [pathmat] (m-1-1)   edge (m-1-2)
          (m-1-2) edge (m-1-3)
          (m-1-3) edge node [auto] {$j_Y^\sharp$} (m-1-4)
          (m-1-4)   edge (m-1-5);
  \end{tikzmat}
  splits as a sequence of superring sheaves. (The splitting is required to be unital on the support $j_{Y,0}(Y_0)$ of $j_{Y_0,*}\sh O_Y$.) If $j_Y$ is a morphism over $S$, then it is retractable as such if and only if the sequence splits as a sequence of $\smash{p_{X,0}^{-1}\sh O_S}$-algebras.

\medskip
As in the case of girth, we have the following useful fact. 

\begin{Prop}[retr-trans][transitivity of retractability]
  Let $j_Z:Z\to Y$ and $j_Y:Y\to X$ be embeddings of $S$-superspaces. If $X$ is (locally) retractable around $Y$, and $Y$ is (locally) retractable around $Z$, then $X$ is (locally) retractable around $Z$ with respect to the embedding $j_Y\circ j_Z:Z\to X$.
\end{Prop}

\begin{Ex}
  Let $k=\knums\in\{\reals,\cplxs\}$ and $A_n$ be the Grassmann superalgebra, constructed in \thmref{Ex}{grass-girth}. Define $X_n\defi(\reals,\sh O_X)$, where $\sh O_X\defi\sh C^\infty_\reals\otimes_\knums A_n$ and $\sh C^\infty_\reals$ is the sheaf of $\knums$-valued smooth functions on $\reals$. There is a natural embedding $j:X_0\to X_n$, given by 
  \[
    j_0\defi\id_\reals,\quad
    j^\sharp(f\tau_{i_1}\dotsm\tau_{i_k})\defi\delta_{k,0}f,
  \]
  for any locally defined smooth function $f$.

  For simplicity, consider the case $n=2$. We may define for any $\lambda\in\knums$ a retraction 
  \[
    r_\lambda:X\to X_0,\quad 
    r_\lambda^\sharp(f)\defi f+\lambda f'\tau_1\tau_2.
  \]
  Thus, retractions of embeddings are far from unique, even in the simplest examples. This fact causes complications in the theory of integration on supermanifolds. 
\end{Ex}

\subsection{Reduction and body}

Reduction of locally ringed spaces removes nilpotency. Both for schemes and for supermanifolds, this results in a space with a sheaf of functions whose sections are determined completely by their values. 

For superspaces, one has yet another reduction, by which everything generated by odd superfunctions is removed. The outcome is called the body of the superspace. There is a further way of associating an even superspace to an arbitrary superspace, the so-called even part, which is constructed by simply taking the even part of the structure sheaf. This construction is of lesser importance, but technically useful. 

\begin{Cons}[even-embedding][reduction of superspaces]
  Let $X$ be a superspace and $x\in X_0$ a point. We define $\vkappa(x)\defi\vkappa_X(x)\defi\sh O_{X,x}/\ger m_{X,x}$; this is a field called the \Define{residue field} at $x$. (The superspaces we will be mostly considering will be defined over a field that at any point coincides with the residue field, but such a restriction is not \emph{a priori} necessary.) If $U\subseteq X_0$ is open, $x\in U$, and $f\in\sh O_X(U)$, we define $f(x)\in\vkappa(x)$ to be the image of the germ $f_x$ in $\vkappa(x)$. This quantity is called the \Define{value} of $f$ at $x$.
  
  For any open subset $U\subseteq X_0$, define an ideal $\sh N_X\subseteq\sh O_X$ by
  \begin{align*}
    \sh N_X(U)&\defi\Set1{f\in\sh O_X(U)}{f(x)=0\text{ \fa} x\in U}\\
    &=\Set1{f\in\sh O_X(U)}{f_x\in\ger m_{X,x}\text{ \fa}x\in U}.
  \end{align*}
  
  Define $\sh O_{X_0}\defi\sh O_X/\sh N_X$. Endowed with this sheaf, $X_0$ is a superspace, which, by abuse of notation, we will denote by $X_0$ and call the \Define{reduced superspace} or \Define[superspace!reduction]{reduction} of $X$.  If $X$ is a $k$-superspace, then so is $X_0$. In the latter case, $\vkappa(x)$ is an extension of $k$ for any $x\in X_0$. There is a canonical morphism of ($k$-)superspaces $j_{X_0}:X_0\to X$, 
  \[
    j_{X_0}\defi\Parens1{j_{X,0}={\id_{X_0}},\sh O_X\to\sh O_X/\sh N_X=\sh O_{X_0}},
  \]
  which is by its mere definition a  thickening with vanishing ideal $\sh N_X$. If $j_{X_0}$ is an isomorphism, then $X$ is called \Define[superspace!reduced]{reduced}.
\end{Cons}

The reduction satisfies a universal property and thus defines a functor.

\begin{Lem}[red-universal]
  Let $X$ be a superspace. Then any morphism $Y\to X$, where $Y$ is reduced, factors uniquely through $j_{X_0}$.
\end{Lem}

\begin{Prop}[red-nat][functoriality of reduction]
  There is a \Define{reduction functor}, denoted by $(-)_0:\SSp\to\SSp$, defined on objects by $X\mapsto X_0$, and on morphisms by letting $\vphi_0:X_0\to Y_0$ for $\vphi:X\to Y$ be the unique factorisation of $\vphi\circ j_{X_0}$ through $j_{Y_0}$. It is right adjoint to the inclusion of reduced superspaces in superspaces.
  
  The canonical morphism $j_{X_0}:X_0\to X$ defines a natural transformation of functors $(-)_0\to{\id}$. The statement carries over to $k$-superspaces.
\end{Prop}

  We will call a retraction of the embedding $j_{X_0}:X_0\to X$ simply a \Define{retraction} of $X$. A retraction is the same as the datum of a sheaf morphism $\sh O_{X_0}\to\sh O_X$ of local superrings \scth for the induced $\sh O_{X_0}$-algebra structure on $\sh O_X$, the canonical morphism $\sh O_X\to\sh O_{X_0}$ is $\sh O_{X_0}$-linear.

  In more conventional terminology, if $r$ is a retraction, one says that the image under $r^\sharp$ of $\Gamma(\sh O_{X_0})$ in $\Gamma(\sh O_X)$ is a \emph{function factor}. Retractions play an important role in the theory of integration on supermanifolds. 

  Similarly, we shall say that $X$ has \Define{finite girth}, \Define{girth} $h$, \Define{locally finite girth}, or \Define{countable girth} if it does so around $X_0$. 

  We denote by $\sh N_X^\infty$ the ideal $\bigcap_{N=1}^\infty\sh N_X^N$, and call it the \Define{Whitney ideal} of the superspace $X$. Thus, $X$ has countable girth if and only if $\sh N_X^\infty=0$.

\begin{Def}[body][body of a superspace]
  Let $X$ be a superspace. Define its \Define[superspace!body]{body} $X_\ev$ as the superspace $X_\ev:=(X_0,\sh O_{X_\ev})$, where $\sh O_{X_\ev}\defi\sh O_X/\sh I_{X_\ev}$ and $\sh I_{X_\ev}\defi\sh O_{X,\odd}+\sh O_{\smash{X,\odd}}^2$ is the ideal generated by $\sh O_{\smash{X,\odd}}$. When $X$ is a $k$-superspace, then so is $X_\ev$.
  
  The body is equipped with a thickening $j_{X_\ev}:X_\ev\to X$, which is a morphism of $k$-superspaces if $X$ is a $k$-superspace. When it is an isomorphism, $X$ is called \Define[superspace!even]{even}.
\end{Def}

The body of a superspace enjoys analogous properties to the reduction. 

\begin{Lem}[body-univ]
  Let $X$ be a superspace. The body $X_\ev$ enjoys the following universal property: Any morphism $Y\to X$, where $Y$ is even, factors uniquely through $j_{X_\ev}$.
\end{Lem}

\begin{Prop}[even-nat][functoriality of the body]
  There is a \Define{body functor}, denoted by $(-)_\ev:\SSp\to\SSp$, which is defined on objects by $X\mapsto X_\ev$ and on morphisms by letting $\vphi_\ev:X_\ev\to Y_\ev$ for $\vphi:X\to Y$ be the unique factorisation of $\vphi\circ j_{X_\ev}$ through $j_{Y_\ev}$. It is right adjoint to the inclusion of even superspaces in superspaces.
  
  The canonical morphism $j_{X_\ev}:X_\ev\to X$ defines a natural transformation of functors $(-)_\ev\to{\id}$. The statement carries over to $k$-superspaces.
\end{Prop}

The inclusion of even superspaces in superspaces also possesses a left adjoint, the so-called \Define{even part}.

\begin{Def}[ev-part][even part of a superspace]
  Let $X$ be a superspace. Define its \Define[superspace!even part]{even part} by $X^\ev\defi(X_0,\sh O_{X,\ev})$. When $X$ is a $k$-superspace, then so is $\smash{X^\ev}$. A canonical morphism $\psi_X:X\to\smash{X^\ev}$ is defined by setting $\psi_{X,0}\defi{\id}_{X_0}$ and taking $\smash{\psi_X^\sharp:\sh O_{X^\ev}=\sh O_{X,\ev}\to\sh O_X}$ to be the canonical inclusion. 
\end{Def}

\begin{Lem}[ev-part-univ]
  Let $X$ be a superspace. Then $(X^\ev,\psi_X)$ has the following universal property: Any morphism $X\to Y$, where $Y$ is even, factors uniquely through $\psi_X$.
\end{Lem}

\begin{Prop}[ev-part-func][functoriality of the even part]
  There is an \Define[even part!functor]{even part functor} with image in the even superspaces, denoted by $(-)^\ev:\SSp\to\SSp$, defined on objects by $X\mapsto X^\ev$ and on morphisms by letting $\vphi^\ev:X^\ev\to Y^\ev$ for $\vphi:X\to Y$ be the unique factorisation of $\psi_Y\circ\vphi$ through $\psi_X$. It is left adjoint to the inclusion of even superspaces in superspaces. 

  The canonical morphism $\psi_X:X\to X^\ev$ defines a natural transformation of functors ${\id}\to(-)^\ev$. The statement carries over to $k$-superspaces.
\end{Prop}

\begin{Ex}
  Consider the $k$-superalgebra
  \[
    A\defi k[T|\tau]/(T^{N+1})
  \]
  and $X\defi(*,A)$, where it is understood that $T$ and $\tau$ are even and odd indeterminates, respectively, and $*$ denotes the singleton space. Then $X_\ev=X^\ev=(*,A_\ev)$ where 
  \[
    A_\ev= k[T]/(T^{n+1}).
  \]
  On the other hand, we have $X_0=(*,k)$.
\end{Ex}

\begin{Ex}
  Let $X=(*,A)$ where $A=A_2$ is the Grassmann algebra on two generators, see \thmref{Ex}{grass-girth}, and $*$ denotes the singleton space. Then $X^\ev=(*,B)$ where
  \[
    B= k[T]/(T^2).
  \]
  On the other hand, we have $X_0=X_\ev=(*,k)$. 
\end{Ex}

The reduction, body, and even part functors preserve embeddings.

\begin{Prop}[redbodev-emb][reduction, body, and even part of embeddings]
  If $j:X\to Y$ is an (open resp.~closed) embedding, then so are $j_0:X_0\to Y_0$, $j_\ev:X_\ev\to Y_\ev$, and $j^\ev:X^\ev\to Y^\ev$.
\end{Prop}

\begin{proof}
  The statement is entirely obvious for open embeddings. Thus, let $j:X\to Y$ be closed embeddding. We have $\vphi\circ j_{X_0}=j_{Y_0}\circ\vphi_0$ where $j_{X_0}$ and $j_{Y_0}$ are closed embeddings. Thus, $\vphi_0$ is a closed embedding, by \thmref{Lem}{emb-factor}. For the body, the same argument goes through. 

  For the even part, observe that the vanishing ideal of $j$ is graded, so that $(j^\ev)^\sharp$ is manifestly a surjective sheaf map.
\end{proof}

The residue fields in the construction of the reduction allow us to speak about values of `superfunctions', \ie sections of the structure sheaf in single points of $X_0$. If $X$ is a $k$-superspace $k$ is contained in  $\vkappa(x)$ for all $x\in X_0$. Thus one can define sections with values in $A$ for any subset $A$ of $k$.

\begin{Def}[aval-sfn][superfunctions with specified values]
  Let $(X,\sh O_X)$ be a $k$-superspace. For any $A\subseteq k$, we define a sheaf $\sh O_{X,A}$ by
  \[
    \sh O_{X,A}(U)\defi\Set1{f\in\sh O_X(U)}{f(x)=j_{X_0}^\sharp(f)(x)\in A\text{ \fa}x\in U}.
  \]
  The local sections are called \Define[superfunction!with values in A]{superfunctions with values in $A$}. 
\end{Def}

Note that if $A$ is a subfield of $k$, then $\sh O_{X,A}$ is an $A$-superalgebra sheaf, because odd superfunctions have value zero. We have the following corollary of \thmref{Prop}{red-nat}.

\begin{Cor}
  Let $A\subseteq k$ and $\vphi:X\to Y$ a morphism of $k$-superspaces. Then $\vphi^\sharp$ induces a sheaf morphism $\vphi^{-1}_0\sh O_{Y,A}\to\sh O_{X,A}$.
\end{Cor}

Let $\vphi:X\to Y$ be an embedding of $k$-superspaces and $A\subseteq k$ such that $0\in A$. We will say that $\vphi$ is \Define{$A$-valued} if $\vphi$ induces a surjective sheaf map $\vphi_0^{-1}\sh O_{Y,A}\to\sh O_{X,A}$. Any thickening with a locally nilpotent ideal is $A$-valued, and so are the canonical embeddings $j_{X_0}:X_0\to X$ and $j_{X_\ev}:X_\ev\to X$. If $\vphi$ is an $A$-valued embedding, then so are $\vphi_\ev$, $\smash{\vphi^\ev}$, and $\vphi_0$. 

  Let $k$ be a field and $X$ be a $k$-superspace. Define
  \[
    X_0(k)\defi\Set1{x\in X_0}{\vkappa(x)=k}.
  \]
  Recall that the elements of $X_0(k)$ are called \Define{$k$-rational points} of $X_0$. If $X_0=X_0(k)$, then the local sections of $\sh O_{X_0}$ may be considered as $k$-valued functions. In this case, we call $X$ \Define[superspace!rational]{$k$-rational}. The following lemma is standard in algebraic geometry.

\begin{Lem}[rat-points][{$k$}-rational points]
  Let $*=(*,k)$ be the terminal object in $\SSp_k$. For any $k$-superspace $X$, there is a natural bijection
  \[
    X(*)\defi\Hom0{*,X}\to X_0(k).
  \]
  If $x\in X_0$, then $x\in X_0(k)$ if and only if $\sh O_{X,x}=k\oplus\ger m_{X,x}$. 
\end{Lem}

A morphism between super-ringed spaces which happen to be superspaces is not necessary local, \ie a morphism of superspaces. If the superspaces are $k$-rational and the ring morphisms $k$-linear, this is no longer an option.

\begin{Prop}[mor-local][locality of morphisms]
  Let $X$ and $Y$ be $k$-rational $k$-superspaces and $\vphi:X\to Y$ a morphism of super-ringed space such that $\vphi^\sharp$ is $k$-linear. Then $\vphi^\sharp(f)(x)=f(\vphi_0(x))$, for any $x\in X_0$ and $f\in\sh O_{Y,\vphi_0(x)}$. In particular, $\vphi$ is local, and therefore, a morphism of $k$-superspaces.
\end{Prop}

\begin{proof}
  Let $y=\vphi_0(x)$. In view of \thmref{Lem}{rat-points} and \thmref{Prop}{super-local}, for $f_y\in\sh O_{Y,y}$, $f(y)\in\vkappa_Y(y)=k$ is the unique $\lambda\in k$ \scth $f_y-\lambda$ is not invertible. By assumption, $\vphi^\sharp_y(f_y)-\lambda=\vphi_y^\sharp(f_y-\lambda)$, and this is invertible whenever $\lambda\neq f(y)$. Let $g=\vphi^\sharp(f)$. Since $x$ is $k$-rational, there is some $\lambda\in k=\vkappa_X(x)$ for which $g_x-\lambda$ is not invertible, namely, $\lambda=g(x)$. Hence, $\vphi^\sharp(f)(x)=g(x)=f(y)$. In particular, if $f_y\in\ger m_{Y,y}$, \ie $f(y)=0$, then $\vphi^\sharp(f)(x)=0$, \ie $\vphi^\sharp(f)_x\in\ger m_{X,x}$.
\end{proof}

\begin{Rem}
  In the literature, morphisms of supermanifolds (to be defined below) are often not assumed to be local. Since supermanifolds are rational superspaces, the preceding proposition shows that the locality of morphisms of supermanifolds is automatic. However, in considering morphisms from more general superspaces to supermanifolds, it is more natural to assume locality, as we do.
\end{Rem}

\subsection{Weil thickenings}

In this subsection, we introduce a construction called Weil thickening, which in its original form goes back to A.~Weil \cite{weil-pointproches}, compare also Ref.~\cite{kolar-michor-slovak}. It does not change the underlying topological space, but only the structure sheaf. Basic examples of superspaces arise by Weil thickening; moreover, it is useful for the infinitesimal study of morphisms. It has two ingredients, a superspace and a Weil superalgebra, which is a special type of local superalgebra.

\begin{Def}[weil-def][Weil superalgebras]
  Let $k$ be a field. A \Define[Weil superalgebra]{Weil $k$-super\-al\-gebra} is a finite-dimensional $k$-super\-al\-gebra $A$ with a graded nilpotent ideal $\ger m$ such that $A=k\oplus\ger m$. Any Weil super\-al\-gebra is a local superring, and $\ger m$ is the maximal (graded) ideal; it consists exactly of the nilpotent elements of $A$. The minimal $h$ \scth $\ger m^{h+1}=0$ is called the \Define{girth} of $A$. Equivalently, one may define Weil superalgebras as finite-dimensional local superalgebras (by Krull's Theorem). 
\end{Def}

Weil superalgebras are obtained by truncation of {polynomial superalgebras}.

\begin{Ex}
  Let $A$ be a $k$-superalgebra where $\chr k=0$. Then $A$ is a Weil superalgebra if and only if $A\cong k[x_1,\dotsc,x_k|\xi_1,\dotsc,\xi_\ell]/I$ for some graded ideal $I$ \scth $I\supseteq(x_1,\dotsc,x_k,\xi_1,\dotsc,\xi_\ell)^N$, for some $N>0$. In particular, $\bigwedge(k^q)^*$ and $k[\eps]/(\eps^{k+1})$ are Weil $k$-superalgebras, of girth $q$ and $k$, respectively.
\end{Ex}

\begin{Ex}[super-dual][super-dual numbers and (multi)jet superalgebras]
  Let $T_i$ and $\tau_j$ be even and odd indeterminates, respectively. We define
  \[
    \sdual_m^{p|q}\defi k[T_1,\dotsc,T_p|\tau_1,\dotsc,\tau_q]/I
  \]
  where the ideal $I$ is generated by all
  \[
    T_1^{\alpha_1}\dotsm T_p^{\alpha_p}\tau_1^{\beta_1}\dotsm\tau_q^{\beta_j},\quad\alpha_i\in\nats,\beta_j\in\{0,1\},\textstyle\sum_i\alpha_i+\sum_j\beta_j=m+1.
  \]
  Then $\sdual_m^{p|q}$ is a Weil $k$-superalgebra of girth $m$, called the \Define[multijet superalgebra]{$m$-multijet superalgebra} of $k$. For $p|q=1|1$ and $m=1$, it is denoted by $\sdual$ and called the superalgebra of \Define{super-dual numbers}. Then $\sdual_\ev=\sdual^{1|0}_1$ is $k[T]/(T^2)$, the algebra of dual numbers familiar in algebraic geometry. Remarkably, it coincides with the even part of the Grassmann algebra $k[\tau_1,\tau_2]$ on two generators.
\end{Ex}

The Weil thickening $X^A$ of a superspace $X$ by a Weil algebra $A$ is obtained simply by tensoring its ring of sections by $A$. It turns out that there is a canonical thickening $X\to X^A$ which is a closed embedding and admits a canonical retraction.

\begin{Cons}[weil-extended-superspace][Weil thickened superspaces]
  Let $X$ be a $k$-superspace. Given a Weil $k$-superalgebra $A=k\oplus\ger m_A$, form $X^A\defi(X_0,\sh O_X\otimes A)$, so
  \begin{equation}\label{eq:weil-ext-max}
    \ger m_{X^A,x}=\ger m_{X,x}\otimes A+\sh O_{X,x}\otimes\ger m_A.
  \end{equation}
  There is a canonical morphism $j_X^A:X\to X^A$, given by $j^A_{X,0}\defi{\id_{X_0}}$, and
  \[
    j_X^{A\sharp}\defi{\id_{\sh O_X}}\otimes\eps:\sh O_{X^A}=\sh O_X\otimes A\to\sh O_X\otimes k=\sh O_X,
  \]
  where $\eps:A\to k$ is the unique algebra morphism. It is clear that $j_X^A$ is a closed embedding, indeed, a thickening. There is a canonical retraction $r_{X^A}:X^A\to X$, given by $r^A_{X,0}\defi\id_{X_0}$ and
  \[
    r_X^{A\sharp}\defi{\id_{\sh O_X}}\otimes\eta:\sh O_X=\sh O_X\otimes k\to\sh O_{X^A}=\sh O_X\otimes A,
  \]
  where $\eta:k\to A$ is the unique algebra morphism.
\end{Cons}

By \thmref{Prop}{retr-trans}, $r_X^A:X^A \to X$ can be combined with retractions of embeddings of other superspaces into $X$.

\begin{Prop}[weil-extended-retractable][retractability of Weil thickenings]
  Let $j_Y:Y\to X$ be a retractable embedding of $k$-superspaces and $A$ a Weil $k$-superalgebra. Then $X^A$ is retractable around $Y$. 
\end{Prop}

Weil thickenings can be reinterpreted as products with superspaces whose underlying topological space is a singleton. 

\begin{Prop}[weil-extended]
  Let $A$ be a Weil $k$-superalgebra, and let $\Spec A=(*,A)$, where $*$ is the singleton space. For any $k$-superspace $X$, $X\times\Spec A$ exists in $\SSp_k$, and is given by $X^A$. If $B$ is another Weil $k$-superalgebra, then $\Spec(A\otimes B)=\Spec A\times\Spec B$.  
\end{Prop}

\begin{proof}
  Let $Y$ be a $k$-superspace. Then $\Hom0{Y,\Spec A}$ is the same as the set of all morphisms of $k$-superalgebras $\psi^\sharp:A\to\Gamma(\sh O_Y)$. Consider the map
  \[
    \Hom0{Y,X}\times\Hom0{Y,\Spec A}\to\Hom0{Y,X^A}:(\vphi,\psi^\sharp)\mapsto(\vphi_0,\vphi^\sharp\otimes\psi^\sharp),
  \]
  given by $(\vphi^\sharp\otimes\psi^\sharp)(f\otimes a)=\vphi^\sharp(f)\cdot\psi^\sharp(a)$. It is clearly a natural bijection. This proves that $X^A=X\times\Spec A$, and the second statement follows.
\end{proof}

\thmref{Prop}{weil-extended} shows that $j_X^A:X\to X^A=X\times\Spec A$ is just $j_X^A=({\id_X},*)$, where $*\to\Spec A$ is the canonical embedding and $*$ is the terminal $k$-superspace. By the same token, the morphism $r_X^A:X^A=X\times\Spec A\to X$ is just $r_X^A=p_1$, the first projection.

\begin{Rem}
  In algebraic geometry, one defines $\Spec R$ for any commutative ring $R$ to be the collection of prime ideals with the Zariski topology, endowed with the structure sheaf obtained by localisation of $R$. 
  
  Note that for a Weil $k$-superalgebra $A$ any prime ideal is maximal, justifying our above definition. Indeed, let $\ger p$ be a prime ideal in $A$. We have $A=k[x_a]/I$ for some set $(x_a)$ of homogeneous indeterminates, and some ideal $I$ containing $(x_a)^N$ for some $N$. Then $\ger p=\ger q/I$ for some prime ideal $\ger q$ of $k[x_a]$ containing $I$. This implies $x_a^N\in\ger q$ for all $a$, but since $\ger q$ is a radical ideal, it follows that $\ger q=(x_a)$, so that $\ger p=\ger m_A$. 
\end{Rem}

  Using \thmref{Prop}{finht-trans}, we obtain that for any finite girth embedding $j_Y:Y\to X$ the embedding $j_X^A\circ j_Y:Y\to X^A$ into a Weil thickening $X^A$ of $X$ has finite girth. We can, however, do better and determine the girth of $j_X^A$ precisely.

\begin{Prop}[weil-extended-girth][girth of Weil thickenings]
  Let $j_Y:Y\to X$ be an embedding of $k$-superspaces of girth $q$ and $A$ a Weil $k$-superalgebra of girth $h$.
  
  Then the embedding $j\defi j_X^A\circ j_Y:Y\to X^A$ has girth $q+h$. In particular, $X\times \Spec A=X^A$ has girth $h$ around $j_X^A:X\to X^A$. Similarly, if $j_Y$ has countable girth, then so has $j$.
\end{Prop}

\begin{proof}
  We may assume that $j_Y:Y\to X$ is a closed embedding. The vanishing ideal of $j_X^A:X\to X^A$ is simply $\sh O_X\otimes\ger m$, where $\ger m=\ger m_A$ is the maximal ideal of $A$. Moreover, by definition of $j_X^A$, we have
  \[
    j^\sharp(f\otimes a)=j_Y^\sharp(f)\cdot\eps(a),\quad
    f\otimes a\in\sh O_{X^A}(U)=\sh O_X(U)\otimes A.
  \]
  Since $\eps(a)\in k$ and $k$ has no zero divisors, we have $j^\sharp(f\otimes a)=0$ if and only if $j_Y^\sharp(f)=0$ or $\eps(a)=0$. It follows that the vanishing ideal $\sh I$ of $j$ is given by
  \[
    \sh I=\sh I_Y\otimes A+\sh O_X\otimes\ger m
  \]
  where $\sh I_Y$ is the vanishing ideal of $j_Y:Y\to X$. Next, we compute inductively that
  \[
    \textstyle\sh I^j/\sh I^{j+1}=\bigoplus_{i=0}^j\sh I_Y^i/\sh I_Y^{i+1}\otimes\ger m^{j-i}/\ger m^{j-i+1},
  \]
  where we agree to write $\sh I^0=\sh O_{X^A}$, $\sh I_Y^0=\sh O_X$, and $\ger m^0=A$. Therefore, the left-hand side is zero if and only $j>q+h$. 

  If $j$ has countable girth, we argue similarly. Indeed, for $N>h$, we have
  \[
    \textstyle\sh I^N\subseteq\sum_{k=0}^h\sh I^{N-k}_Y\otimes\ger m^k,
  \]
  so that $\sh I^N\subseteq\sh I_Y^{N-h}\otimes A$. Let $f\in\bigcap_{N=1}^\infty\sh I^N(U)$. Choosing a $k$-basis $(a_\ell)$ of $A$, we have a unique representation $\smash{f=\sum_\ell f_\ell\otimes a_\ell}$ where $f_\ell\in\sh O_X(U)$. By the above, $\smash{f_\ell\in\bigcap_{N=1}^\infty\sh I_Y^N(U)=0}$, which proves our claim.
\end{proof}

\subsection{Formally Noetherian and tidy superspaces}

In this subsection we introduce two regularity conditions, one for superspaces and one for embeddings of superspaces. Given a superspace $X$, we consider the `formal' ring $\Hat{\sh O}_{X,x}\defi\sh O_{X,x}/\ger m_{X,x}^\infty$ for any $x\in X_0$. If these rings are all Noetherian, the sections of the structure sheaf are not too far from being determined by their restrictions to infinitesimal normal neighbourhoods of points.
 
An embedding $j:Y\to X$ of superspaces is called \emph{tidy} if the restriction of $j^\sharp$ to the stalk at some point $x\in j_0(Y_0)$ is determined by the cosets of germs in $\Hat{\sh O}_{X,y}$ for $y$ in a neighbourhood of $x$ in $j_0(Y_0)$.  A superspace with Noetherian formal rings, for which the identity is tidy, will be called a tidy superspace. The concept of tidiness will be instrumental for the coordinate description of morphisms.

In this subsection, we describe a number of general properties of tidy embeddings and superspaces.

\medskip
  Let $R$ be a local supercommutative superring with maximal ideal $\ger m$. We let $\ger m^\infty\defi\bigcap_{N=1}^\infty\ger m^N$ be the Whitney ideal of the superspace $(*,R)$, $*$ denoting the singleton space.

\begin{Prop}[krull][Krull's Theorem]
  Let $M'\subseteq M$ be an inclusion of $R$-supermod\-ules. If $R/\ger m^\infty$ is Noetherian and $M/\ger m^\infty M$ is finitely generated, then $\ger m^\infty M+M'=\bigcap_{N=0}^\infty(\ger m^NM+M')$. This applies for $M=R$ and $M'$ any graded ideal.
\end{Prop}

\begin{proof}
  The ring $R/\ger m^\infty$ is local and Noetherian, with maximal ideal $\ger n\defi\ger m/\ger m^\infty$. Consider $P\defi M/(\ger m^\infty M+M')$. Since $M/\ger m^\infty M$ is a finitely generated graded $R/\ger m^\infty$-module, so is $P$. Let $Q\defi\bigcap_{n=1}^\infty(\ger n^nP)$. Then $\ger nQ=Q$ by the Artin--Rees Lemma, so by the Nakayama Lemma, $Q=0$. The canonical map $M\to P$ sends $\ger m^nM+M'$ to $\ger n^nP$. We find that $\bigcap_{n=1}^\infty(\ger m^nM+M')$ is mapped to $0$, as desired. 
\end{proof}

\begin{Def}[noeth-jets][Formally Noetherian superspaces]
  Let $X$ be a superspace. If $\sh O_{X,x}/\ger m_{X,x}^\infty$ is Noetherian for all $x\in X_0$, then we say that $X$ is \Define{formally Noetherian}.
\end{Def}

Formal Noetherianity is stable under Weil thickenings and passage to subspaces.

\begin{Lem}[weilext-noetherianjets]
  Let $A$ be a Weil $k$-superalgebra and $X$ a formally Noetherian $k$-superspace. Then $X^A$ is formally Noetherian.
\end{Lem}

\begin{proof}
  Clearly, $\ger m_{X^A,x}$ contains $\ger m_{X,x}\otimes k$, so $\ger m_{X^A,x}^\infty$ contains $\ger m_{X,x}^\infty\otimes k$, so there is a surjective homomorphism $\smash{\sh O_{X,x}/\ger m_{X,x}^\infty\otimes A\to\sh O_{X^A,x}/\ger m_{X^A,x}^\infty}$. The left-hand side is finitely generated over the Noetherian ring $\sh O_{X,x}/\ger m_{X,x}^\infty$, so is Noetherian itself.
\end{proof}

\begin{Lem}[noetherian-emb]
  Let $j:Y\to X$ be an embedding of superspaces. If $X$ is formally Noetherian, then so is $Y$.
\end{Lem}

\begin{proof}
  Observe simply that for $x=j_0(y)$, the surjection $\sh O_{X,x}\to\sh O_{Y,y}/\smash{\ger m_{Y,y}^\infty}$ induced by $j^\sharp$ factors through the Noetherian superring $\sh O_{X,x}/\ger m_{X,x}^\infty$.
\end{proof}

\begin{Def}[tidy][tidiness]
  Let $j:Y\to X$ be an embedding of superspaces with vanishing ideal $\sh I$. We say that $j$ is \Define[embedding!tidy]{tidy}\index{tidy!embedding} if for any $x\in j_0(Y_0)$, any open neighbourhood $U\subseteq X_0$ of $x$, and any $f\in\sh O_X(U)$, the following holds:
  \[
    \Parens1{\forall y\in U\cap j_0(Y_0),N\in\nats\,:\,f_y\in\sh I_y+\ger m_{X,y}^N}\ \Rightarrow\ f_x\in\sh I_x.
  \]

  We say that $X$ is \Define[superspace!tidy]{tidy}\index{tidy!superspace} if ${\id_X}$ is a tidy embedding. The latter condition means that the zero ideal is tidy; this implies that the Whitney ideal $\sh N^\infty_X=0$, which is equivalent to the requirement that $X$ have countable girth. Thus, if $X$ is tidy, then it has countable girth. 

  Although no counterexample appears to be known, it seems somewhat doubtful that the converse holds, even under the assumption of formal Noetherianity and finite girth. Trivially, however, any $X$ that is reduced (\ie of girth zero) is tidy. 
\end{Def}

\begin{Rem}[tidy]
Now assume that $X$ is formally Noetherian and let $j:Y\to X$ be an embedding of vanishing ideal $\sh I$. Then $j$ is tidy if and only if for any $x\in j_0(Y_0)$, any open neighbourhood $U\subseteq X_0$ of $x$, and any $f\in\sh O_X(U)$, the following holds:
  \[
    \Parens1{\forall y\in U\cap j_0(Y_0),N\in\nats:f_y\in\sh I_y+\ger m^\infty_{X,y}}\ \Rightarrow\ f_x\in\sh I_x.
  \]
  This follows immediately from \thmref{Prop}{krull}. Notice that this is also the case when $X=Y$ and $j=\id_X$, even if $X$ is not formally Noetherian. 

  In particular, if $\sh O_{X,x}$ is Noetherian \fa $x\in j_0(Y_0)$, then this is always trivially verified, since in this case $\ger m_{X,x}^\infty=0$, by the Nakayama Lemma. 
\end{Rem}

Tidiness can be rephrased on the level of the ideals, as follows. 

\begin{Def}[tidy-ideal][tidying of a vanishing ideal]
  Let $j:Y\to X$ be an embedding of superspaces with ideal $\sh I$. We define the \Define[tidying!of an ideal]{tidying} $\overline{\sh I}$ of $\sh I$ by $\overline{\sh I}\defi\bigcap_{N=1}^\infty\overline{\sh I}_N$, where $\overline{\sh I}_N$ is the ideal sheaf 
  \[
    \overline{\sh I}_N(U)\defi\Set1{f\in\sh O_X(U)}{\forall y\in U\cap j_0(Y_0):f_y\in\ger m_{X,y}^N+\sh I_y}
  \]
  \fa open subsets $U\subseteq X_0$. Then $j$ is tidy if and only if $\sh I=\overline{\sh I}$. Observe that $\sh I$ depends on $j_0$. If $X$ is formally Noetherian, then
  \[
    \overline{\sh I}(U)=\Set1{f\in\sh O_X(U)}{\forall y\in U\cap j_0(Y_0):f_y\in\ger m_{X,y}^\infty+\sh I_y}
  \]
  \fa open $U\subseteq X_0$. (Compare \thmref{Rem}{tidy}.)
\end{Def}




\begin{Prop}[weil-tidy][tidiness of Weil thickenings]
  Let $X$ be a tidy $k$-superspace and $A$ a Weil $k$-superalgebra. Then $X^A$ is tidy.
\end{Prop}

\begin{proof}
  To see that $\id_{X^A}$ is tidy, we argue as in the proof of \thmref{Prop}{weil-extended-girth} in the case of countable girth. Namely, $m_{X^A,x}=\ger m_{X,x}\otimes A+\sh O_{X,x}\otimes\ger m_A$, so for $N>h$, where $h$ is the height of $A$, we have $\ger m_{X^A,x}^N\subseteq\sum_{k=0}^h\ger m_{X,x}^{N-k}\otimes\ger m_A^k$. Given $f\in\sh O_{X^A}(U)$ \scth $f_x\in\ger m_{X^A,x}^\infty$ \fa $x\in U$, we may expand $f_x$ in terms of a basis of $A$. By the assumption that $X$ be tidy, the coefficients must vanish. 
\end{proof}

Tidy embeddings give rise to tidy subspaces.

\begin{Prop}[noetherianjets-countht][tidiness of subspaces]
  For any tidy embedding $j:Y\to X$ of superspaces, $Y$ is tidy. 
\end{Prop}

\begin{proof}[\protect{Proof of \thmref{Prop}{noetherianjets-countht}}]
  Let $V=j_0^{-1}(U)$, $U\subseteq X_0$ being open. Take $f\in\sh O_Y(V)$ \scth $f_y\in\ger m_{Y,y}^N$ \fa $y\in V$ and $N\in\nats$. Passing to an open cover, we may assume that $f=j^\sharp(g)$ \fs $g\in\sh O_X(U)$.
  
  Let $y\in V$ and $x\defi j_0(y)$. Since $j^\sharp$ is a surjective sheaf map, $\ger m_{Y,y}=j^\sharp(\ger m_{X,x})$. For any $N$, we have $f_y\in\ger m_{Y,y}^N$. Hence, there are finitely many $g_{kn}\in\ger m_{X,x}$, $k=1,\dotsc,N$, \scth $\smash{f_y=j^\sharp\Parens1{\sum_ng_{1n}\dotsm g_{Nn}}}$. So for $h\defi g_x-\sum_ng_{1n}\dotsm g_{Nn}$, we have $j^\sharp(h)=0$, \ie $h\in\sh I_x$. Thus, $\smash{g_x\in\bigcap_{N=1}^\infty(\ger m_{X,x}^N+\sh I_x)}$. Since $x$ was arbitrary, it follows that $g_x\in\sh I_x$ for $x\in U\cap j_0(Y_0)$. Hence, $f=j^\sharp(g)=0$, and $Y$ is tidy. 
\end{proof}

\begin{Cor}[open-tidy]
  Let $j:Y\to X$ be an open embedding. Then $j$ is tidy if and only if $j(Y)$ is tidy.
\end{Cor}

Conversely, tidy subspaces give rise to tidy embeddings. 

\begin{Prop}[tidy-compos][preservation of tidiness]
  Let $j:Z\to Y$ and $i:Y\to X$ be embeddings. If $j$ is tidy, then so is $i\circ j$.
\end{Prop}

\begin{proof}
  Denote the vanishing ideals of $i$, $j$ and $k\defi i\circ j$ by $\sh I$, $\sh J$ and $\sh K$, respectively. Let $U\subseteq X_0$ be open, $V\defi i_0^{-1}(U)$, and $W\defi j_0^{-1}(V)$. Consider $f\in\sh O_X(U)$ and $z\in W$. We put $y\defi j_0(z)$ and $x\defi i_0(y)$.
  
  We assume that for all $u\defi i_0(v)$, $v\defi j_0(w)$, $w\in W_0$, and all $N\in\nats$, we have $f_u\in\ger m_{X,u}^N+\sh K_u$. There is $g\in\ger m_{X,u}^N$ \scth $j^\sharp(i^\sharp(f_u-g))=k^\sharp(f_u-g)=0$, so 
  \[
    i^\sharp(f)_v\in i^\sharp(g)+\ker j^\sharp_v\subseteq\ger m_{Y,v}^N+\sh J_v.
  \]
  Because $v$, $N$ were arbitrary, we find $\smash{i^\sharp(f)_y\in\sh J_y}$, thanks to the tidiness of $i$. This gives $k^\sharp(f_x)=j^\sharp(i^\sharp(f)_y)=0$, \ie $f_x\in\sh K_x$, so that $k$ is tidy.
\end{proof}

\begin{Cor}[tidy-emb]
  Any embedding of a tidy superspace is itself tidy.
\end{Cor}

\begin{proof}
  Apply \thmref{Prop}{tidy-compos} with $Z=Y$ and $j={\id_Y}$.
\end{proof}

We give an example of an embedding that is not tidy.

\begin{Ex}[fatpt][the fat point]
  Let $X\defi(\reals,\sh O_X)$, where $\sh O_X=\sh C^\infty_\reals$ is the sheaf of $\knums$-valued smooth functions, where $\knums\in\{\reals,\cplxs\}$. Let $Y$ be the subspace given by $Y_0=\{0\}$ and $\sh O_Y=\sh O_{X,0}$, the algebra of germs of smooth functions at $0$. 

  The canonical embedding $j:Y\to X$ is given by letting $j_0$ be the embedding of the point $0$ in $\reals$ and setting $j^\sharp:j_0^{-1}\sh O_X=\sh C^\infty_{\reals,0}\to\sh O_Y=\sh C^\infty_{\reals,0}$ equal to the identity of the algebra of germs. Then $j$ is manifestly not a tidy embedding, since $\ger m_{X,0}^\infty$ is not reduced to zero, by \thmref{Ex}{formal-power}.
\end{Ex}

\begin{Prop}[tidying][existence of tidying]
  Let $X$ be a superspace. There is a closed embedding $t:X^\circ\to X$ \scth $X^\circ$ is tidy, with the following universal property: Whenever $\vphi:Y\to X$ is a morphism, where $Y$ is tidy, there is a unique morphism $\vphi^\circ:Y\to X^\circ$ \scth $\vphi=t\circ\vphi^\circ$. 
\end{Prop}

\begin{proof}
  Define the superspace $X^\circ\defi(X_0,\sh O_{X^\circ})$ by $\sh O_{X^\circ}\defi \sh O_X/{\overline 0}$, where $0$ is the zero ideal and $\overline0$ is its tidying. The definition of the embedding $t:X^\circ\to X$ is obvious. Clearly, $t$ is tidy, so that $\id_{X^\circ}$ is tidy, by \thmref{Prop}{noetherianjets-countht}.
  
  Let $\vphi:Y\to X$ be a morphism, where $Y$ is tidy.  Let $V\subseteq Y_0$ be open, and $f$ a local section of $\sh O_X$ defined on a neighbourhood $U$ of $\vphi_0(V)$. Assume that $f\in{\overline0}(U)$, \ie for any $y\in V$, $x\defi\vphi_0(y)$, we have $f_x\in\ger m_{X,x}^\infty$. Then $\vphi^\sharp(f)_x=\vphi^\sharp(f_x)\in\ger m_{Y,y}^\infty$. Because $y$ was arbitrary and $Y$ is tidy, $\smash{\vphi^\sharp(f)=0}$, and $\vphi$ factors uniquely through $t$ to a morphism $\vphi^\circ:Y\to X^\circ$.
\end{proof}

\begin{Def}[tidying][tidying]
  The morphism $t:X^\circ\to X$ and the space $X^\circ$ constructed in \thmref{Prop}{tidying} are called the \Define{tidying} of $X$. They are unique up to canonical isomorphism.
\end{Def}

The following is immediate from the construction in \thmref{Prop}{tidying}.

\begin{Cor}[tidying][ideals of tidyings]
  Let $j:Y\to X$ be an embedding of superspaces with ideal $\sh I$. The ideal of the tidy embedding $j\circ t:Y^\circ\to X$ is $\overline{\sh I}$, the tidying of $\sh I$.
\end{Cor}

We give an easy example of a tidying. 

\begin{Ex}[fat-tidying][tidying of the fat point]
  Recall the embedding $j:Y\to X$ from \thmref{Ex}{fatpt}. The ideal of this embedding is $0$. Its tidying is the ideal sheaf $\overline0$ whose local sections are 
  \[
    \overline{0}(U)=\Set1{f\in\sh O_X(U)}{0\in U\ \Rightarrow\ f_0\in\ger m_{X,0}^\infty}
  \]
  \ie the functions that are flat at $0$. Hence, the tidying of $Y$ is given by 
  \[
    Y^\circ_0=\{0\},\quad\sh O_{Y^\circ}=\sh C^\infty_{\reals,0}/\ger m_{\reals,0}^\infty=\knums\llbracket T\rrbracket,
  \]
  by the theorem of \'E.~Borel, and $j^\circ\defi j\circ t^\circ:Y^\circ\to X$ is given by 
  \[
    j^\circ_0\defi j_0,\quad (j^\circ)^\sharp:(j^\circ_0)^{-1}\sh O_X=\sh C^\infty_{\reals,0}\to\knums\llbracket T\rrbracket,
  \]
  where $(j^\circ)^\sharp$ is the map associating to any function germ its Taylor series at the point $0$. In other words, $j^\circ$ is the embedding constructed in \thmref{Ex}{formal-power}. 
\end{Ex}

\begin{Lem}[tidying-properties]
  Let $X$ be a superspace and $U\subseteq X$ be an open subspace. Then $U^\circ$ is an open subspace of $X^\circ$.
\end{Lem}

\begin{proof}
  This is obvious by construction. Indeed, $U_0=(U^\circ)_0$ is open in $X_0=(X^\circ)_0$ and $\sh O_{U^\circ}=\sh O_X|_{U_0}/{\overline 0}|_{U_0}=(\sh O_X/\overline0)|_{U_0}=\sh O_{X^\circ}|_{U_0}$.
\end{proof}

Thickenings into normally Noetherian tidy superspaces can be approximated by tidy thickenings of finite girth.

\begin{Prop}[noether-tidy-approx][tidy approximation]
  Let $j:Y\to X$ be a thickening with $X$ formally Noetherian and tidy. The tidying $X_k$ of the $k$th infinitesimal normal neighbourhood of $Y$ has girth $\sle k$ around $Y$ and $X=\varinjlim_kX_k$.
\end{Prop}

\begin{proof}
  The $k$th order infinitesimal normal neighbourhood $X^{(k)}$ gives a natural thickening $\smash{X^{(k)}}\to X$ with vanishing ideal $\sh I^{k+1}$ where $\sh I\defi\ker j^\sharp$. Let $j_k:X_k\to X$ be the canonical thickening induced by tidying. Since $X^{(k)}$ has finite girth $\sle k$ and there is a thickening $X_k\to X^{(k)}$, so has $X_k$. 

  The $X_k$ form an inductive system, by \thmref{Prop}{tidying}. Let $\tilde\jmath:\tilde X\to X$ denote the inductive limit of the thickenings $j_k$. It is a thickening, since the associated projective system $\sh O_{X_{k+1}}\to\sh O_{X_k}$ of sheaves is surjective and thus satisfies the Mittag-Leffler condition. 

  Thus, to show that $j$ is an isomorphism, it suffices to prove that $\tilde\jmath^\sharp$ is injective. By assumption, $X$ is formally Noetherian, so on applying Krull's theorem (\thmref{Prop}{krull}) twice, we find 
  \[
    \ker\tilde\jmath^\sharp_x\subseteq\bigcap_{k,N}(\ger m^N+\sh I_x^k)=\bigcap_k(\ger m_{X,x}^\infty+\sh I_x^k)\subseteq\bigcap_k(\ger m_{X,x}^\infty+\ger m_{X,x}^k)=\ger m_{X,x}^\infty
  \]
  for any $x\in X_0$, since manifestly $\sh I_x\subseteq\ger m_{X,x}$. Since $X$ is tidy by assumption, we have $\smash{\ker\tilde\jmath^\sharp}=0$, so that $X\cong\tilde X$. 
\end{proof}

\begin{Prop}[tidy-equaliser][equalisers of tidy superspaces]
  Equalisers exist in the category of locally Hausdorff tidy superspaces, and they are tidy embeddings.
\end{Prop}

\begin{proof}
  Let $\phi,\psi:X\to Y$ be morphisms where $Y_0$ is locally Hausdorff. Form $Z_0\defi\Set1{x\in X_0}{\phi_0(x)=\psi_0(x)}$. This is a locally closed subset of $X_0$, because the diagonal of $Y_0$ is locally closed. We let $j_0:Z_0\to X_0$ be the embedding and set $\sh O_Z\defi j_0^{-1}\sh O_X|_U/\sh I$ where $U\subseteq X_0$ is largest open subset in which $Z_0$ is closed and $\sh I$ is the ideal generated by $\im(\phi^\sharp-\psi^\sharp)$. Then $Z\defi(Z_0,\sh O_Z)$ is a superspace and there is a canonical embedding $j:Z\to X$ \scth $\phi\circ j=\psi\circ j$. 
  
  Let $t:Z^\circ\to Z$ be the tidying of $Z$, so that $Z^\circ$ is tidy, and $\jmath^\circ\defi j\circ t$ is a tidy embedding, by \thmref{Cor}{tidy-emb}. Clearly, $\phi\circ\jmath^\circ=\psi\circ\jmath^\circ$.  Moreover, if $X_0$ is locally Hausdorff, then so is $Z_0$, since $j_0$ is an embedding.
  
  Now, let $W$ be tidy and $\vphi:W\to X$ be a morphism \scth $\phi\circ\vphi=\psi\circ\vphi$. Then $\vphi^\sharp(\sh I)=0$, so that $\vphi$ factors uniquely through $j$. Since $W$ is tidy, $\vphi$ even factors uniquely through $\jmath^\circ$, by virtue of \thmref{Prop}{tidying}, to a morphism $\vphi^\circ:W\to Z^\circ$. Let $\vrho:W\to Z^\circ$ be another morphism with $\jmath^\circ\circ\vrho=\vphi=\jmath^\circ\circ\vphi^\circ$. Seeing that $\jmath^\circ$ is an embedding and so a monomorphism of superspaces, we find that $\vrho=\vphi^\circ$.
\end{proof}

\begin{Prop}[open-eq][equalisers of open subspaces]
  Let $\phi,\psi:X\to Y$ be morphisms of locally Hausdorff tidy superspaces with equaliser $\vphi:Z\to X$. Let $U\subseteq X$ and $V\subseteq Y$ be open subspaces with $U\subseteq\phi^{-1}(V)\cap\psi^{-1}(V)$. Then the equaliser of $\phi|_U,\psi|_U:U\to V$ is an open subspace of $Z$.
\end{Prop}

\begin{proof}
  Let $W\defi\vphi^{-1}(U)$ and $\vphi|_W:W\to U$ the morphism induced by $\vphi$. In view of \thmref{Cor}{open-tidy}, $W$ is tidy, and it is certainly locally Hausdorff. We claim that $\vphi|_W:W\to U$ is the equaliser of $\phi|_U$ and $\psi|_U$.
  
  Let $\vrho:S\to U$ be a morphism where $W$ is locally Hausdorff and tidy and $\phi|_U\circ\vrho=\psi|_U\circ\vrho$. Suppose that $\sigma\defi j_U\circ\vrho$ factors uniquely through $\vphi:Z\to X$ to a morphism $\sigma':S\to Z$. Then $\vphi_0(\sigma'_0(S_0))=\sigma(S_0)\subseteq U_0$, so $\sigma'(S_0)\subseteq W_0$, and $\sigma'$ factors uniquely through $j_W:W\to Z$ to a morphism $\vrho':S\to W$. Thus, $j_U\circ\vphi|_W\circ\vrho'=\vphi\circ j_W\circ\vrho'=\sigma=j_U\circ\vrho$, so that $\vphi|_W\circ\vrho'=\vrho$. If $\vrho'':S\to W$ is another such morphism, then $\vrho'=\vrho''$, because $\vphi|_W$ is an embedding.
\end{proof}

\section{Morphisms and local coordinates}\label{sec:leites}

This section is devoted to the study of the question to which extent morphisms with range in a supermanifold may be given a coordinate description. We address this matter by introducing the concept of Leites regularity. It depends on parameters, which indicate whether the supermanifolds in question are smooth, real or complex analytic, and whether their functions take real or complex values. 

In the previous sections, we have considered superspaces without any restrictions on base fields. From now on, we will be working exclusively over the real or complex field; this restriction may not be necessary, and it is conceivable that the material we will be discussing here could be treated in greater generality, for instance, following the approach of Ref.~\cite{al-inf}. However, our focus is different here, so we will content ourselves with this setup.

\subsection{Premanifolds}

In this section, we introduce (pre)manifolds in the traditional manner, using atlases. Then we show how to view them as superspaces. Finally, we state the classical Hadamard lemma, which is instrumental in introducing coordinates into the study of superspaces. 

Applications show that it is necessary to consider, apart from cases of the real and complex supermanifolds, superspaces that have the mixed structure of a complex sheaf of superfunctions and an underlying real manifold as their body. We will from the very beginning work with a pair $(\knums,\Bbbk)$ of fields, where $\Bbbk\subseteq\knums$ are either $\reals$ or $\cplxs$. This will allow us to handle the case of \emph{cs} manifolds together with the more conventional cases of real and complex supermanifolds. 

\begin{Def}[premfd-def][{$\Bbbk$}-premanifolds]
  Let $\varpi=\infty,\omega$ and $X_0$ a topological space. A \Define[atlas!over $\Bbbk$ of class $\sh C^\varpi$]{$\Bbbk$-atlas of class $\sh C^\varpi$} is an open cover $(U_i)$ of $X_0$, together with topological embeddings $\phi_i:U_i\to\Bbbk^{p_i}$ \scth for any $i,j$,
  \[
    \phi_{ij}\defi\phi_i\circ\phi_j^{-1}:U_{ji}\defi\phi_j(U_i\cap U_j)\to U_{ij}
  \]
  is smooth (for $\varpi=\infty$) or analytic (for $\varpi=\omega$).  Two atlases are called \emph{equivalent} if their union is an atlas.
  
  Endowed with an equivalence class of $\Bbbk$-atlases of class $\sh C^\varpi$, $X_0$ is called a \Define[premanifold!of class $\sh C^\varpi$]{$\Bbbk$-premanifold of class $\sh C^\varpi$}. If in addition, $X_0$ is Hausdorff, then it is called a \Define[manifold!of class $\sh C^\varpi$]{manifold}.
  
  A continuous map $\vphi_0:X_0\to Y_0$ of $\Bbbk$-premanifolds of class $\sh C^\varpi$ is called \emph{of class $\sh C^\varpi$} (resp.~smooth for $\varpi=\infty$, resp.~analytic for $\varpi=\omega$), if there exist atlases $(U_i,\phi_i)$, $(V_i,\psi_i)$ of $X_0$ resp.~$Y_0$ \scth $\vphi_0(U_i)\subseteq V_i$ and $\psi_i\circ\vphi_0\circ\phi_i^{-1}$ is smooth (for $\varpi=\infty$) resp.~analytic (for $\varpi=\omega$).
  
  A map $\vphi_0:X_0\to Y_0$ of class $\sh C^\varpi$ is called an \Define[embedding!open]{open embedding} if $\vphi_0$ induces a homeomorphism $X_0\to\vphi_0(X_0)$, endowed with the relative topology from $Y_0$, and in addition, $\vphi_0(X_0)$ is open in $Y_0$. In this case, $\vphi_0$ is an open map; open embeddings in this sense will turn out to coincide with those previously defined, \vq \thmref{Cor}{man-emb}.
\end{Def}

We view premanifolds as superspaces for which we impose the regularity conditions and the options for the underlying fields mentioned in the introduction. Accordingly, we consider five different classes of premanifolds.

\begin{Def}[mfd-def][premanifolds over $(\knums,\Bbbk)$]
  Let $\varpi=\infty,\omega$ and $X_0$ a $\Bbbk$-premanifold of class $\sh C^\varpi$. We let $\sh O_X$ be the sheaf of commutative $\knums$-algebras of $\knums$-valued functions of class $\sh C^\varpi$ on $X_0$. The five possibilities are summarised in the following table:\vspace{0.1cm}
  \begin{center}\def\arraystretch{1.2}
  \begin{longtable}{|c||c|c|c|c|}
  \hline
  type of manifold&$\knums$&$\Bbbk$&$\varpi$&$\sh O_X$\\
  \hline\hline
  real smooth or $\mathcal C^\infty$ & $\reals$ & $\reals$ & $\infty$ &$\reals$-valued smooth functions\\
  \hline
  real analytic or $\mathcal C^\omega$ & $\reals$ & $\reals$ & $\omega$&$\reals$-valued real analytic functions\\
  \hline
  \emph{c} & $\cplxs$ & $\reals$ & $\infty$&$\cplxs$-valued smooth functions\\
  \hline
  \emph{c} analytic & $\cplxs$ & $\reals$ & $\omega$ &$\cplxs$-valued real analytic functions\\
  \hline
  complex analytic & $\cplxs$ & $\cplxs$ & $\omega$ & complex analytic functions\\
  \hline
  \end{longtable}
  \end{center}

  Here, the `\emph{c}' stands for `complex', following the terminology, introduced by J.~Bernstein \cite{deligne-morgan}, of `\emph{cs}' (`complex super') manifolds for $\cplxs$-superspaces whose structure sheaf is modeled on the sheaf of smooth functions of real variables with values in a complex Grassmann algebra. 

  The pair $X\defi(X_0,\sh O_X)$ will be called a \Define{(pre)manifold over $(\knums,\Bbbk)$} of class $\sh C^\varpi$. In the case $\varpi=\infty$ resp.~$\varpi=\omega$, we say that $X$ is smooth resp.~analytic. We also say that $X$ is a (pre)manifold of the type specified in the table. 
For instance, if $(\knums,\Bbbk)=(\cplxs,\reals)$ and $\varpi=\omega$, we say that $X$ is \Define[premanifold!\emph{c} analytic]{\emph{c} analytic}. Accordingly, the \emph{c} (pre)manifolds are ordinary smooth (pre)manifolds, endowed with the sheaf of complex-valued smooth functions.
\end{Def}

The following proposition shows that the premanifolds from \thmref{Def}{mfd-def} are superspaces.

\begin{Prop}[mfd-max]
  Let $X$ be a premanifold over $(\knums,\Bbbk)$ and $x\in X_0$. Then $\sh O_{X,x}$ is local and the maximal ideal $\ger m_{X,x}$ consists of the germs of functions vanishing at $x$. Hence, $X$ is a formally Noetherian reduced $\knums$-rational $\knums$-superspace, and in particular tidy. The morphisms of super-ringed spaces associated with maps of class $\sh C^\varpi$ are morphisms of $\knums$-superspaces.
\end{Prop}

\begin{proof}
  All of the statements are straightforward, save the formal Noetherianity. Since $X$ is reduced, this follows from the following \thmref{Lem}{aff-noetherianjets}.
\end{proof}

\begin{Lem}[aff-noetherianjets]
  Let $p\in\nats$. The $\knums$-superspace $\aff^p$ is formally Noetherian.
\end{Lem}

\begin{proof}
  If $\varpi=\omega$, then $\sh O_{\aff^p,x}$ is the convergent power series ring in $p$ indeterminates, so $\ger m_{\aff^p,x}^\infty=0$ and $\sh O_{\aff^p,x}$ is itself Noetherian. If $\varpi=\infty$, then 
  \[
    \sh O_{\aff^p,x}/\ger m_{\aff^p,x}^\infty=\knums\llbracket T_1,\dotsc,T_p\rrbracket
  \]
  by \'E.~Borel's theorem, and this is again a Noetherian ring.
\end{proof}

The $\knums$-superspace $X=(X_0,\sh O_X)$ is called the \Define[superspace!associated with a (pre)manifold]{$\knums$-superspace associated with $X_0$}. Any such $\knums$-superspace will be called a \Define{premanifold over $(\knums,\Bbbk)$} (of class $\sh C^\varpi$), and a  \Define{manifold over $(\knums,\Bbbk)$} (of class $\sh C^\varpi$) if $X_0$ is a manifold (\ie Hausdorff).

With any map $\vphi_0:X_0\to Y_0$ of $\sh C^\varpi$ of $\Bbbk$-premanifolds of class $\sh C^\varpi$, we may associate a morphism $\vphi=(\vphi_0,\vphi^\sharp):X\to Y$ of the associated $\knums$-superspaces, by $\vphi^\sharp(f)\defi f\circ\vphi_0$ \fa $f\in\sh O_Y(U)$, $U\subseteq Y_0$ open. This obviously defines a functor.

\medskip
Finally, we turn to the Hadamard lemma. It is valid in each of our five categories of premanifolds. To have an efficient formulation, we simply  write $\aff^p$\label{affp-def} for the space $\Bbbk^p$ together with the sheaf $\sh O_{\aff^p}$ of $\knums$-valued functions of class $\sh C^\varpi$, and view it as a manifold over $(\knums,\Bbbk)$ of class $\sh C^\varpi$.

Recall from \thmref{Def}{aval-sfn} that for any $\knums$-superspace $X$ the sheaf of $\Bbbk$-valued superfunctions on $X$ is denoted by $\sh O_{X,\Bbbk}$.

\begin{Prop}[hadamard-basic][Hadamard lemma]
  Let $U\subseteq\aff^p$ be open. Consider functions $f_1,\dotsc,f_n\in\sh O_{\aff^p,\Bbbk}(U)$, and set $Y\defi\bigcap_{j=1}^nf_j^{-1}(0)$. Assume that the derivatives $df_1(x)\dotsc,df_n(x)$ are linearly independent at any $x\in Y$. Then any point of $Y$ possesses an open neighbourhood $V\subseteq U$ \scth every $g\in\sh O_{\aff^p}(V)$, $g(V\cap Y)=0$, admits a representation $g=\sum_{j=1}^nf_jg_j$ \fs $g_j\in\sh O_{\aff^p}(V)$. Moreover, we may choose $V$ \scth there are $\sh C^\varpi$ functions $f_{n+1},\dotsc,f_p$ on $V$, for which the map $f=(f_1,\dotsc,f_p)$ is a $\sh C^\varpi$ diffeomorphism onto its open image.
\end{Prop}

The following three corollaries of the Hadamard lemma describe the algebraic nature of the sheaf $\sh O_{\aff^p}$. 

\begin{Cor}[hadamard-higherorder]
  Assume the conditions of \thmref{Prop}{hadamard-basic}. Take $x\in Y$, and fix $\smash{f_{n+1},\dotsc, f_p}$ as in its statement. There is an open neighbourhood $V\subseteq U$ of $x$ \scth for every $g\in\sh O_{\aff^p}(V)$ and any $N\sge 0$, there is a unique polynomial $P=\sum_\alpha P_\alpha(f_{n+1},\dotsc,f_p)\cdot T^\alpha$ of degree $\sle N$ in $T_1,\dotsc,T_n$, with coefficients of class $\sh C^\varpi$, \scth $g-P(f_1,\dotsc,f_n)\in I_Y(V)^{N+1}$, where
  \[
    I_Y(V)\defi\Set1{h\in\sh O_{\aff^p}(V)}{h(V\cap Y)=0}.
  \]
\end{Cor}

The following two special cases will be of particular importance.

\begin{Cor}[hadamard-pt][pointwise Hadamard lemma]
  Let $x\in\aff^p_0$ and $f\in\sh O_{\aff^p,x}$. For any $N\sge0$, there is a unique $P\in\knums[T_1,\dotsc,T_p]$, $\deg P\sle N$, \scth
  \[
    f-P(t-x)=f-P(t_1-x_1,\dotsc,t_p-x_p)\in\ger m_{\aff^p,x}^{N+1},
  \]
  where $t=(t_1,\dotsc,t_p)$ is the identity of $\aff^p_0$.
\end{Cor}

\begin{Cor}[hadamard-diag][diagonal Hadamard lemma]
  Let $\Delta:\aff^p\to\aff^{2p}$ be the morphism associated with the diagonal map. Then $\Delta$ is a closed embedding and its ideal $\sh I\subseteq\sh O_{\aff^{2p}}$ is generated by the sections $x_j-y_j$, $j=1,\dotsc,p$.
\end{Cor}

\subsection{Leites regularity}\label{sec:Leites regularity}

In this subsection, we introduce affine superspace, which is the local model for supermanifolds. Affine superspace carries a canonical system of `coordinate' functions. A natural question is to which extent morphisms to affine superspace are determined by their action on these coordinate functions. This leads to the concept of Leites regularity of superspaces. As we shall see in the following, this gives rise to a remarkably robust and versatile category.

\begin{Def}[aff-def][standard affine superspace]
  In what follows, we fix the fields $\knums$, $\Bbbk$, and the regularity class $\varpi=\infty,\omega$. We will write $\aff^{p|q}$ for the $\knums$-superspace $\aff^p\times\Spec\bigwedge(\knums^q)^*$, where $\aff^p$ is the $\knums$-superspace associated with the $\Bbbk$-manifold $\Bbbk^p$ of class $\sh C^\varpi$, and the Weil superalgebra $\bigwedge(\knums^q)^*=\knums[\theta_1,\dotsc,\theta_q]$ is the supercommutative $\knums$-superalgebra freely generated on $q$ odd generators. We call this the \Define{affine superspace} of dimension $p|q$ (of class $\sh C^\varpi$ over $(\knums,\Bbbk)$). 
\end{Def}

Note that $\aff^{p|0}$ coincides with $\aff^p$, as introduced in the previous subsection. 

\begin{Not}[std-coord][standard coordinate functions]
  We denote the standard coordinate functions on $\aff^p$ by $t_1,\dotsc,t_p$. Furthermore, we let $\theta_1,\dotsc,\theta_q$ be the standard generators of $\bigwedge(\knums^q)^*=\knums[\theta_1,\dotsc,\theta_q]$, considered as sections of $\sh O_{\aff^{0|q}}$. We use the same letters if $t_i$, $\theta_j$ are considered as sections of $\sh O_{\aff^{p|q}}$.
  
  We abbreviate the tuple $(t_1,\dotsc,t_p,\theta_1,\dotsc,\theta_q)$ by $(t,\theta)$. Sometimes, we will not wish to distinguish explicitly between even and odd members of $(t,\theta)$. In this case, we will write the standard coordinates $t=(t_a)$ where $a=1,\dotsc,p+q$, and we will not impose any particular order on the even and odd members of $t$. We will say that $t=(t_a)$ is in \Define[standard coordinates!in standard order]{standard order} if $\Abs0{x_a}=\ev$ \fa $a\sle p$ and $\Abs0{x_a}=\odd$ \fa $a>p$.
  
  On any open $U\subseteq\aff^{p|q}_0=\aff^p_0$, every $g\in\sh O_{\aff^{p|q}}(U)$ has a unique expression
  \begin{equation}\label{eq:std-coord-repn}
    g=\sum\nolimits_Ig_I\theta^I\mathtxt{where}\theta^I\defi\theta_{i_1}\dotsm\theta_{i_k}
  \end{equation}
  for $I=(1\sle i_1<\dotsm<i_k\sle q)$ and $g_I\in\sh O_{\aff^p}(U)$. In particular, the $\theta_j$ allow us to consider $\sh O_{\aff^p}\subseteq\sh O_{\aff^{p|q}}$ in a distinguished fashion.
\end{Not}

\begin{Rem}[inf]
  Conceivably, one might define affine superspaces modeled on more general (topological or bornological) super-vector spaces of possibly infinite dimension. This is a potentially interesting generalisation of the theory presented here. 
\end{Rem}

The idea behind the Leites regularity of a $\knums$-superspace is to guarantee that a morphism to a standard affine superspace is completely determined by the pullbacks of the standard coordinates under the morphism. It turns out that it is sufficient to demand this for the even standard coordinates.  

In what follows, for given $\knums$-superspaces $X$ and $S$, we write
\[
  X(S)\defi\Hom0{S,X}
\]
for the set of morphisms $S\to X$ in the category $\SSp_\knums$.

\begin{Def}[def-reg][Leites regularity and subregularity]
  We say that a $\knums$-superspace $X$ is \Define{(Leites) regular} (of class $\sh C^\varpi$ over $(\knums,\Bbbk)$) if for any open subspace $U\subseteq X$, and any $p$, the map
  \begin{equation}\label{eq:regular-def}
    \aff^p(U)=\Hom0{U,\aff^p}\to\Gamma(\sh O_{U,\Bbbk,\ev}^p):\vphi\mapsto(\vphi^\sharp(t_1),\dotsc,\vphi^\sharp(t_p)),
  \end{equation}
  where $t_j$ are the standard coordinate functions on $\aff^p$, is bijective. If these maps are merely all injective, then we say that $X$ is \Define{(Leites) subregular}.
  
  We denote the full subcategory of $\SSp_\knums$ whose objects are the regular superspaces of class $\sh C^\varpi$ over $(\knums,\Bbbk)$ by $\SSp_{\knums,\Bbbk}^\varpi$. Notice that $\aff^p$ is not the $p$-fold product of $\aff^1$ in $\SSp_\knums$, so that we are obliged to work with $\aff^p$. It will turn out that $\aff^p$ \emph{is} the $p$-fold product of $\aff^1$ in the subcategory of Leites regular superspaces. 
\end{Def}

\begin{Rem}
  Along the lines of \thmref{Rem}{inf}, it may be interesting to consider infinite-dimensional versions of the concept of Leites regularity, where one replaces finite tuples of superfunctions with a suitably defined set of vector-valued superfunctions, for instance using appropriate notions of tensor product or vector-valued differentiability. 
\end{Rem}

\begin{Ex}[flat-notsubregular][The fat point is not Leites subregular]
  Let $X=(*,\sh C^\infty_{\reals,0})$ be the fat point, constructed in \thmref{Ex}{fatpt}. We assume that $\Bbbk=\reals$ and $\varpi=\infty$. We show that $X$ is not Leites subregular.

  Indeed, let $\vphi:X\to\aff^1$ be the embedding constructed in \loccit, that is, $\vphi_0(*)=0$ and $\vphi^\sharp(f)= f_0$, the germ of $f$ at $0$. Then $\vphi^\sharp(t)=t_0$. 

  Let $T:\sh C^\infty_{\reals,0}\to\knums\llbracket t\rrbracket$ be the Taylor series map. By a result of Reichard \cite{reichard}*{Satz 2}, there is an algebra homomorphism $\phi:\knums\llbracket t\rrbracket\to\sh C^\infty_{\reals,0}$ \scth $T\circ\phi=\id$ and $\phi$ maps convergent power series to germs of real-analytic functions. 

  We define $\psi:X\to\aff^1$ by $\psi_0(*)\defi0$ and $\psi^\sharp(f)\defi\phi\circ T(f)$. Then $\psi^\sharp(t)=t_0+h$, where $h$ is flat at $0$. But since $h$ is real-analytic, it follows that $h=0$. On the other hand, $\psi^\sharp(f)=0$ for any smooth function $f$ that is flat at $0$. 

  However, there are smooth functions that are flat at $0$, but whose germ at $0$ is non-trivial, see \thmref{Ex}{fatpt}. Hence, $\vphi\neq\psi$ and $X$ is indeed not Leites subregular. 
\end{Ex}

The  following lemma shows that for Leites regular superspaces, morphisms to affine superspace are indeed completely described by their tuples of `component' superfunctions. In fact, it shows that the category of Leites regular superspaces is the largest subcategory of $\SSp_\knums$ in which Leites's morphism theorem holds: Morphisms to $\aff^{p|q}$ are in bijection with tuples of `component' superfunctions.

\begin{Lem}[affpq-mor]
  Let $X$ be a Leites regular superspace. The following natural map is a bijection, for any non-negative integers $p$ and $q$,
  \begin{gather*}
    \aff^{p|q}(X)\to\Gamma(\sh O_{X,\Bbbk,\ev})^p\times\Gamma(\sh O_{X,\odd})^q\\
    \vphi\mapsto(\vphi^\sharp(t_1),\dotsc,\vphi^\sharp(t_p),\vphi^\sharp(\theta_1),\dotsc,\vphi^\sharp(\theta_q)).
  \end{gather*}
\end{Lem}

\begin{proof}
  Since $\aff^{p|q}=\aff^{p|0}\times\aff^{0|q}$ and the statement holds for $q=0$ by definition, it is sufficient to prove it for $p=0$. Observe that $\Gamma(\sh O_{\aff^{0|q}})=\knums[\theta_1,\dotsc,\theta_q]$ is the free supercommutative algebra generated by the odd indeterminates $\theta_i$. Hence, the underlying map of any morphism $\vphi:X\to\aff^{0|q}$ necessarily is the constant map $X_0\to *=(\aff^{0|q})_0$, and $\vphi^\sharp(P)=P(\vphi^\sharp(\theta_1),\dotsc,\vphi^\sharp(\theta_q))$ for any $P\in\knums[\theta_1,\dotsc,\theta_q]$.
\end{proof}

As we will see, the affine superspace $\aff^{p|q}$ is Leites regular. This is not hard to prove, but we derive it from more general statements below (\thmref{Prop}{manifold-mor} and \thmref{Cor}{weil-ext-reg}). Once it has been established, \thmref{Lem}{affpq-mor} shows that in the category of Leites regular $\knums$-superspaces, we have $\aff^{p+r|q+s}=\aff^{p|q}\times\aff^{r|s}$.

\begin{Def}[coords][charts and coordinate systems]
  Let $X$ be a superspace. A \Define[chart!local]{local chart} is an open embedding of $\knums$-superspaces $\vphi:U\to\aff^{p|q}$ where $U\subseteq X$ is an open subspace. We say that $\vphi$ is \Define[chart!defined on $U$]{defined on $U$}. If $\vphi$ is defined on $X$, then we call $\vphi$ a \Define[chart!global]{global chart}. The tuple $p|q$ is called the \Define[chart!dimension (graded)]{(graded) dimension} of $\vphi$.

  Let $\vphi$ be a local chart. Then the tuple
  \[
   (x,\xi)\defi(x_1,\dotsc,x_p,\xi_1,\dotsc,\xi_q)\defi(\vphi^\sharp(t_1),\dotsc,\vphi^\sharp(t_p),\vphi^\sharp(\theta_1),\dotsc,\vphi^\sharp(\theta_q))
  \]
  is called a system of \Define[coordinate system!local]{local coordinates} (defined on $U$). If $\vphi$ is a global chart, we say that $(x,\xi)$ is a system of \Define[coordinate system!global]{global coordinates}. 
\end{Def}

Note that $\knums$-superspaces admitting enough charts to cover the underlying topological space are automatically Leites regular. Thus, in view of \thmref{Lem}{affpq-mor}, charts and coordinate systems are in bijection. We may and will also call $p|q$ the \Define[coordinate system!dimension (graded)]{(graded) dimension} of $(x,\xi)$.
  
  Occasionally, we will not wish to distinguish in our notation between even and odd members of a system of local coordinates. In this case, we will instead write $x=(x_a)$ where $x_a\defi\vphi^\sharp(t_a)$, $a=1,\dotsc,p+q$. We will say that $x=(x_a)$ is in \Define[local coordinates!in standard order]{standard order} if in its definition, $(t_a)$ is in standard order.
  
In view of Equation \eqref{eq:std-coord-repn}, given a system of local coordinates $(x,\xi)$ defined on $X|_U$, any $f\in\sh O_X(U)$ has a unique representation
\begin{equation}\label{eq:coord-repn}
  f=\sum\nolimits_If_I\xi^I\mathtxt{where}\xi^I\defi\xi_{i_1}\dotsm\xi_{i_k}
\end{equation}
and $f_I=x^\sharp(g_I)$ \fs $g_I\in\sh O_{\aff^p}(\vphi_0(U))$, where $x$ also denotes the morphism $X|_U\to\aff^p$ defined by $x^\sharp(t_i)=x_i$. Observe that the $f_I$ are \emph{not} local sections of $\sh O_{X_0}$; rather, they are even local sections of $\sh O_X$.

Although it is justifiable to write $f_I=g_I(x_1,\dotsc,x_p)$, we will avoid this notation since it is prone to provoke the misunderstanding that the $f_I$ are `ordinary functions'.

\begin{Ex}
  The standard coordinate system $(t,\theta)$ on $\aff^{p|q}$ is a global system of coordinates, and so is its restriction to any open subspace.
\end{Ex}

\begin{Ex}
  On $\aff^{1|2}$, a non-standard coordinate system is given by $(t+\theta_1\theta_2,\theta_1,\theta_2)$, where $(t,\theta_1,\theta_2)$ is the standard coordinate system. 
\end{Ex}

At this point it would be possible to introduce supermanifolds of class $\sh C^\varpi$ over $(\knums,\Bbbk)$ as superspaces that are covered by open subsets admitting global charts. The reason we postpone this is that at the moment, we do not have the tools available to also define relative supermanifolds. This will be possible once we have a good understanding of the category of Leites regular  $\knums$-superspaces and its full subcategory of locally finitely generated $\knums$-superspaces (see Section~\ref{sec:relative supermanifolds}).

As the following proposition shows, in the framework of Leites regular superspaces, morphisms admit a simple description in terms of coordinates.

\begin{Prop}[mor-coord][morphisms \vs coordinates]
  Let $X$ and $Y$ be Leites regular $\knums$-superspaces of class $\sh C^\varpi$ over $(\knums,\Bbbk)$, $y=(y_a)$ a global coordinate system in standard order on $Y$, and $p|q$ the graded dimension of $Y$. Then the map   
  \begin{gather*}
    Y(X)\to\Set1{x\in\Gamma(\sh O_X)_\ev^p\times\Gamma(\sh O_X)_\odd^q}{x(X_0)\subseteq y(Y_0)}\colon\vphi\mapsto \vphi^\sharp(y)
  \end{gather*}
  is a bijection. Here, $x(X_0)$---and similarly, $y(Y_0)$---stands for the set
  \[
    x(X_0)\defi\Set1{(x_1(o),\dotsc,x_p(o))}{o\in X_0}\subseteq\aff^p_0.
  \]
\end{Prop}

\begin{proof}
  The coordinate system is given by $y=\psi^\sharp(t)$ for a unique global chart $\psi:Y\to\aff^{p|q}$, and $\psi$ is an open embedding. Let $x\in\Gamma(\sh O_X)_\ev^p\times\Gamma(\sh O_X)_\odd^q$ \scth $x(X_0)\subseteq y(Y_0)$. Since $Y$ is Leites regular, there exists by \thmref{Lem}{affpq-mor} a unique morphism $\vrho:X\to\aff^{p|q}$ \scth $x=\vrho^\sharp(t)$.
  
  The condition $x(X_0)\subseteq y(Y_0)$ guarantees that $\vrho$ factors uniquely through $\psi$ to a morphism $\sigma:X\to Y$, by \thmref{Prop}{emb-factor}. In particular, $\sigma^\sharp(y)=\vrho^\sharp(t)=x$. The proof of the converse statement is similar.
\end{proof}

 Observe that the sets $x(X_0)$ and $y(Y_0)$ in the statement of \thmref{Prop}{mor-coord} are contained in $\aff^p_0=\Bbbk^p$ (and not only in $\knums^p$).

\begin{Def}[mapcond][mapping condition]
  Let $X$ and $Y$ be regular $\knums$-superspa\-ces. Given $x=(x_a)\in\Gamma(\sh O_X)^{p+q}$ and a global coordinate system $y=(y_a)$ of $Y$ \scth $\Abs0{x_a}=\Abs0{y_a}$ \fa $a$, the condition
  \begin{equation}\label{eq:mapcond}
    x(X_0)\subseteq y(Y_0)
  \end{equation}
  from \thmref{Prop}{mor-coord} is called the \Define{mapping condition} on $x$ and $y$.
\end{Def}

  \thmref{Prop}{mor-coord} admits an obvious modification for coordinate systems not in standard order. We will use it in this more general form without further notice.

\subsection{Morphisms with non-standard affine target}

Sometimes, it is useful to study the affine superspaces $\aff^{p|q}$ independent of coordinates. In this subsection, we digress briefly to show how this can be accomplished in full generality. The basis for this is given by the following definition.

\begin{Def}[kalpha-def][super-vector spaces over $(\knums,\Bbbk)$]
  Let $V=V_\ev\oplus V_\odd$ be a finite-dimensional super-vector space over $\Bbbk$, together with a fixed $\knums$-structure on $V_\odd$. We call $V$ a \Define[super-vector space!over $(\knums,\Bbbk)$]{$(\knums,\Bbbk)$-super-vector space}. For super-vector spaces $V,W$ over $(\knums,\Bbbk)$, we set
  \begin{align*}
    \GHom[_{\knums,\Bbbk}]0{V,W}\defi\,&(\Hom[_\Bbbk]0{V_\ev,W_\ev}\oplus\Hom[_\knums]0{V_\odd,W_\odd})\\
    &\oplus(\Hom[_\Bbbk]0{V_\ev,W_\odd}\oplus\Hom[_\Bbbk]0{V_\odd,W_\ev}),
  \end{align*}
  with the obvious $\Bbbk$-structure on the even part and the $\knums$-structure on the odd part obtained by acting on $W_\odd$ in the range and on $V_\odd$ in the argument. 

  In particular, we obtain the \Define[super-vector space!dual, over $(\knums,\Bbbk)$]{dual super-vector space over $(\knums,\Bbbk)$}
  \[
    V^*\defi\GHom[_{\knums,\Bbbk}]0{V,\Bbbk}=\Hom[_\Bbbk]0{V_\ev,\Bbbk}\oplus\Hom[_\knums]0{V_\odd,\knums}.
  \]
  Here, we use the $\knums$-structures on $V_\odd$ and $\Hom[_\Bbbk]0{V_\odd,\Bbbk}$ to identify the latter canonically with $\Hom[_\knums]0{V_\odd,\knums}$. The even part of $\GHom[_{\knums,\Bbbk}]0{V,W}$ is denoted by $\Hom[_{\knums,\Bbbk}]0{V,W}$, and its elements are called \Define[morphism! of super-vector spaces over $(\knums,\Bbbk)$]{morphisms of $(\knums,\Bbbk)$-super-vector spaces}. We also set 
  \[
    V\otimes W\defi(V_\ev\otimes_\Bbbk W_\ev\oplus V_\odd\otimes_\knums W_\odd)
    \oplus(V_\ev\otimes_\Bbbk W_\odd\oplus V_\odd\otimes_\Bbbk W_\ev),
  \]
  with the obvious $\Bbbk$- resp.~$\knums$-structure on the even resp.~odd part.
\end{Def}

\begin{Lem}
  With the above definitions, the $\otimes$ and $\underline{\mathrm{Hom}}$ functors on finite-dimensional super-vector spaces over $(\knums,\Bbbk)$ form a pair of adjoint functors and $\Bbbk$ acts as the unit of $\otimes$. 
\end{Lem}

\begin{proof}
  Let $U,V,W$ be finite-dimensional super-vector spaces over $(\knums,\Bbbk)$. We have 
  \begin{align*}
    \Hom[_\Bbbk]0{U_\ev\otimes_\Bbbk V_\ev,W_\ev}&=\Hom[_\Bbbk]0{U_\ev,\Hom[_\Bbbk]0{V_\ev,W_\ev}},\\    
    \Hom[_\Bbbk]0{U_\odd\otimes_\knums V_\odd,W_\ev}&=\Hom[_\knums]0{U_\odd,\Hom[_\Bbbk]0{V_\odd,W_\ev}},\\
    \Hom[_\knums]0{U_\ev\otimes_\Bbbk V_\odd,W_\odd}&=\Hom[_\Bbbk]0{U_\ev,\Hom[_\knums]0{V_\odd,W_\odd}},\\
    \Hom[_\knums]0{U_\odd\otimes_\Bbbk V_\ev,W_\odd}&=\Hom[_\knums]0{U_\odd,\Hom[_\Bbbk]0{V_\ev,W_\odd}}.
  \end{align*}
  It follows that there is a natural bijection 
  \[
    \Hom[_{\knums,\Bbbk}]0{U\otimes V,W}=\Hom[_{\knums,\Bbbk}]0{U,\GHom[_{\knums,\Bbbk}]0{V,W}}.
  \]
  Since $\Hom[_{\knums,\Bbbk}]0{\Bbbk,U}=U_\ev$, it follows that 
  \[
    \Hom[_{\knums,\Bbbk}]0{\Bbbk\otimes U,V}=\Hom[_{\knums,\Bbbk}]0{\Bbbk,\GHom[_{\knums,\Bbbk}]0{U,V}}=\Hom[_{\knums,\Bbbk}]0{U,V}.
  \]
  Hence, $\Bbbk$ is a tensor unit for $\otimes$, and the assertion follows. 
\end{proof}
  
\begin{Cons}[linear-ssp][affine superspace of a super-vector space]
  Let $V$ be a $(\knums,\Bbbk)$-super vector space. The vector space $V_\ev$ may be considered as a $\Bbbk$-premanifold. We will denote the associated premanifold over $(\knums,\Bbbk)$ by the same letter. Let $V$ be the $\knums$-superspace $V_\ev\times\Spec\bigwedge(V_\odd)^*$. We call this the \Define[affine superspace!of a super-vector space over $(\knums,\Bbbk)$]{affine superspace} of $V$. So the reduced space $V_0$ of $V$ is the premanifold over $(\knums,\Bbbk)$ given by $V_\ev$. The dual $V^*$ is contained in $\Gamma(\sh O_{V,\Bbbk})$.
\end{Cons}

  To state promised coordinate-free descriptions of morphisms with an affine target cleanly, the following concepts prove useful.

\begin{Def}[oxk-mod][{$\sh O_{X,\Bbbk}$}-modules]
  Let $X$ be a $\knums$-superspace. Recall the sheaf $\sh O_{X,\Bbbk}$ of $\Bbbk$-valued superfunctions, where $\sh O_{X,\Bbbk}(U)\subseteq\sh O_X(U)$ is defined by the condition $f(x)\in\Bbbk$ \fa $x\in U$, for $U\subseteq X_0$. It is a sheaf of algebras in the category of super-vector spaces over $(\knums,\Bbbk)$. 

  An \Define[oxkmodule@$\sh O_{X,\Bbbk}$-module]{$\sh O_{X,\Bbbk}$-module} is a sheaf $\sh M$ with values in the category of super-vector spaces over $(\knums,\Bbbk)$, which is a graded $\sh O_{X,\Bbbk}$-module. One has an obvious notion of \Define[oxkmodule@$\sh O_{X,\Bbbk}$-module!morphism]{morphism of $\sh O_{X,\Bbbk}$-modules} $\sh M\to\sh N$. The set of all these will be denoted by $\Hom[_{\sh O_{X,\Bbbk}}]0{\sh M,\sh N}$. If $\vphi:Y\to X$ is a morphism of $\knums$-superspaces and $\sh M$ is an $(\sh O_{X,\Bbbk})$-module, then 
  \[
    \vphi^*\sh M\defi\sh O_{Y,\Bbbk}\otimes_{\vphi_0^{-1}\sh O_{X,\Bbbk}}\vphi_0^{-1}\sh M
  \]
  naturally carries the structure of an $\sh O_{Y,\Bbbk}$-module.
\end{Def}

\begin{Cons}[oxk-mod-cons][free $\sh O_{X,\Bbbk}$-modules]
  Let $X$ be a $\knums$-superspace and $V$ a finite-dimensional $(\knums,\Bbbk)$-super-vector space. Let $\sh O_{X,\Bbbk}\otimes V$ be the tensor product in the category of super-vector spaces over $(\knums,\Bbbk)$. This is an $\sh O_{X,\Bbbk}$-module. Similarly, define $\GHom[_{\knums,\Bbbk}]0{V,\sh O_{X,\Bbbk}}$ in the category of super-vector spaces over $(\knums,\Bbbk)$. Its even part is denoted by $\Hom[_{\knums,\Bbbk}]0{V,\sh O_{X,\Bbbk}}$.

  An $\sh O_{X,\Bbbk}$-module $\sh P$ is called \Define[oxkmodule@$\sh O_{X,\Bbbk}$-module!free]{free} if it is isomorphic to $\sh O_X\otimes V$, for some finite-dimensional $(\knums,\Bbbk)$-super-vector space $V$. It is called \Define[oxkmodule@$\sh O_{X,\Bbbk}$-module!locally free]{locally free} if $X_0$ admits an open cover, the restrictions of $\sh P$ to the constituents of which are free.
  
  We have the following canonical isomorphism of $\sh O_{X,\Bbbk}$-modules,
  \[
    \GHom[_{\knums,\Bbbk}]0{V,\sh O_{X,\Bbbk}}\cong\sh O_{X,\Bbbk}\otimes V^*.
  \]
  In particular, $\GHom[_{\knums,\Bbbk}]0{V,\sh O_{X,\Bbbk}}$ is a free $\sh O_{X,\Bbbk}$-module.
\end{Cons}

With the help of this terminology, we can now give a coordinate-free description of morphisms with affine target.

\begin{Cor}[mor-coordfree]
  Let $X$ be a regular $\knums$-superspace and $V$ a $(\knums,\Bbbk)$-super-vector space. Let $V$ also denote the associated affine superspace. Then there is a bijection
  \begin{tikzmat}
    \matrix (m) [mathmat] {V(X) & \Gamma\Parens1{\Hom[_{\knums,\Bbbk}]0{V^*,\mathcal O_{X,\Bbbk}}}=\Gamma\Parens1{(\sh O_{X,\Bbbk}\otimes V)_\ev}\\ };
    \path [pathmat] (m-1-1) edge (m-1-2);
  \end{tikzmat}
  which is natural in $X$ and $V$. Under this map,
  \begin{align*}
    V_\ev(X)&\cong\Gamma\Parens1{\Hom[_{\Bbbk}]0{V_\ev^*,\sh O_{X,\Bbbk,\ev}}}=\Gamma\Parens1{\sh O_{X,\Bbbk,\ev}\otimes_\Bbbk V_\ev},\\
    V_\odd(X)&\cong\Gamma\Parens1{\Hom[_\knums]0{V_\odd^*,\mathcal O_{X,\odd}}}=\Gamma\Parens1{\sh O_{X,\odd}\otimes_\knums V_\odd}.
  \end{align*}
  
  More generally, for each open subspace $U\subseteq V$, there is a natural bijection
  \begin{align*}
    U(X)&\cong\!\Set1{f\in\Gamma(\Hom[_{\knums,\Bbbk}]0{V^*,\sh O_{X,\Bbbk}})}{\forall x\in X_0\exists u\in U_0:f(\cdot)(x)=\Dual0\cdot u}\\
    &=\!\Set1{f=\textstyle\sum_i\!f_i\otimes v_i\in\Gamma((\sh O_{X,\Bbbk}\otimes V)_\ev)}{\forall x\in X_0\,\exists u\in U_0:u=\textstyle\sum\!_if_i(x)v_i}
  \end{align*}
  Here, $\Dual0\cdot\cdot$ denotes the canonical pairing of $V^*$ and $V$.
\end{Cor}

\section{Singular superspaces and relative supermanifolds}\label{sec:lfg}

This section contains our main results. Firstly, we show that Leites regularity of superspaces is preserved under colimits and also under the passage to subspaces and thickenings, at least if they are subregular. We then construct the subcategory of locally finitely generated superspaces. These admit finite limits, and are stable both under passage to tidy subspaces and Weil thickenings.  

Thus, within the category of locally finitely generated superspaces, we may construct the category of relative supermanifolds over a possibly singular base. This also provides a natural framework for the study of Weil functors, which give a uniform description of such natural geometric objects as the even and odd tangent bundle and higher order versions of these. 

\subsection{Constructions of Leites regular superspaces}

So far, we have established that Leites regularity is a desirable property for $\knums$-superspaces. However, we have at this point not yet given any examples of superspaces with this property. In this subsection, we will derive a number of general stability results that will furnish a generous supply. In particular, these result imply that locally finitely generated superspaces, defined in Subsection \ref{sec:locally finitely generated superspaces}, are Leites regular (\thmref{Prop}{fg-countht}), vastly generalising Leites's Chart Theorem \cite{leites}.

We start by showing that premanifolds are Leites regular. Next, we show that Leites regularity is stable under colimits. Moreover, we investigate to which extent it is preserved under taking subspaces and thickenings, for example, Weil thickenings.

For the sake of brevity, from now on, we will drop the prefix `Leites' and just speak of `regular' resp.~`subregular' superspaces.

\begin{Prop}[manifold-mor][regularity of premanifolds]
  Let $X$ be a premanifold over $(\knums,\Bbbk)$. Then $X$ is regular. In particular, $\aff^p$ is regular. 
\end{Prop}

For the \emph{proof}, we make the following observation, which is of independent interest.

\begin{Lem}[finht-canmap-inj]
  Any tidy $\knums$-superspace is subregular.
\end{Lem}

\begin{proof}
  Let $X$ be tidy and let morphisms $\vphi,\psi:X\to\aff^p$ be given \scth $\vphi^\sharp(t_j)=f_j=\psi^\sharp(t_j)$. Then
  \[
    t_j(\vphi_0(x))=j_{X_0}^\sharp\vphi^\sharp(t_j)(x)=j_{X_0}^\sharp\psi^\sharp(t_j)(x)=t_j(\psi_0(x))
  \]
  for any $x\in X_0$ and $j=1,\dotsc,p$, so that $\vphi_0=\psi_0$.
  
  Let $f$ be a local section of $\sh O_{\aff^p}$. Let $N\in\nats$ and $x\in X_0$ be arbitrary, and consider the point $y\defi\vphi_0(x)=\psi_0(x)$. By \thmref{Cor}{hadamard-pt}, there is a polynomial $P\in\knums[t_1-y_1,\dotsc,t_p-y_p]$ \scth $(f-P)_y\in\ger m_{\aff^p,y}^{N+1}$. Then
  \[
    \vphi^\sharp(f)_x-\psi^\sharp(f)_x=\vphi^\sharp(f-P)_x-\psi^\sharp(f-P)_x\in\ger m_{X,x}^{N+1}.
  \]
  Since $x$ and $N$ were arbitrary and $X$ is tidy, we conclude that $\vphi^\sharp(f)-\psi^\sharp(f)$ vanishes. It follows that $\vphi^\sharp(f)=\psi^\sharp(f)$.
\end{proof}

\begin{proof}[\protect{Proof of \thmref{Prop}{manifold-mor}}]
  Given functions $f_1,\dotsc,f_p\in\Gamma(\sh O_{X,\Bbbk})$, there is a unique premanifold map $\vphi_0:X_0\to\Bbbk^p$ of class $\sh C^\varpi$ with components $t_j\circ\vphi_0=f_j$. The associated morphism $\vphi:X\to\aff^p$ (compare the remarks following \thmref{Prop}{mfd-max}) satisfies $\vphi^\sharp(t_j)=f_j$, proving surjectivity. All open subspaces of $X$ are premanifolds and hence tidy, so by \thmref{Prop}{mfd-max}, so \thmref{Lem}{finht-canmap-inj} implies the claim.
\end{proof}

\begin{Cor}[man-emb]
  Consider the functor that sends a $\Bbbk$-premanifold of class $\sh C^\varpi$ to its associated $\knums$-superspace. It maps open embeddings of premanifolds to open embeddings of $\knums$-superspaces, is fully faithful, and its image consists of the premanifolds over $(\knums,\Bbbk)$ of class $\sh C^\varpi$.
\end{Cor}

\begin{proof}
  Let $\vphi_0:X_0\to Y_0$ be a manifold map of class $\sh C^\varpi$ that is an open embedding of topological spaces. The associated morphism $\vphi:X\to Y$ is an open embedding of $\knums$-superspaces. The statement follows from \thmref{Prop}{relative-glue} and \thmref{Prop}{manifold-mor}.
\end{proof}

\begin{Cor}[aff-prod][products of affine superspaces]
  Let $V$ and $V'$ be finite-dimensional super-vector spaces over $(\knums,\Bbbk)$. Then $\aff(V\oplus V')$ is the direct product of $\aff(V)$ and $\aff(V')$ in the category of regular $\knums$-superspaces.
\end{Cor}

\begin{proof}
  For any regular $\knums$-superspace $X$, there is a natural bijection
  \begin{align*}
    \aff(V\oplus V')(X)&=\Gamma\Parens1{\sh O_{X,\Bbbk}\otimes(V\oplus V')}_\ev\\
    &=\Gamma(\sh O_{X,\Bbbk}\otimes V)_\ev\times\Gamma(\sh O_{X,\Bbbk}\otimes V')_\ev=\aff(V)(X)\times\aff(V')(X).  
  \end{align*}
  This proves the assertion.
\end{proof}

The category of regular superspaces is quite robust.

\begin{Prop}[indlim-regular][colimits of regular superspaces]
  The full subcategory of $\SSp_\knums$ of regular (resp.~subregular) $\knums$-superspaces is {cocomplete}.
\end{Prop}

\begin{proof}
  Let $I$ be a small category, $F:I\to\SSp_\knums$ a functor, and $\smash{X=\varinjlim_IF}$ in $\SSp_\knums$ (\vq \thmref{Prop}{local-cocomplete}). Denote the canonical morphisms $F(i)\to X$ by $j_i$. Assume that $F$ takes values in subregular superspaces. Let $\vphi,\psi:X\to\aff^p$ be morphisms \scth $\smash{f_k\defi\vphi^\sharp(t_k)=\psi^\sharp(t_k)}$, $k=1,\dotsc,p$. For any $i$, we form $f^i_k\defi j_i^\sharp(f_k)$. Then $\smash{(\vphi\circ j_i)^\sharp(t_k)=f^i_k=(\psi\circ j_i)^\sharp(t_k)}$, $k=1,\dotsc,p$, so by assumption, we have $\vphi\circ j_i=\psi\circ j_i$. Since $i$ was arbitrary, this shows $\vphi=\psi$, so $X$ is subregular.

  Assume now that $F$ even takes values in regular superspaces, fix even superfunctions $f_1,\dotsc,f_p\in\Gamma(\sh O_{X,\Bbbk,\ev})$, and form $f_k^i$ as above. By assumption, there are unique $\vphi_i:F(i)\to\aff^p$ \scth $\smash{\vphi_i^\sharp(t_k)=f^i_k}$, $k=1,\dotsc,p$. For any $\chi:i\to i'$ in $I$, we have
  \[
    F(\chi)^\sharp(f^{i'}_k)=(j_{i'}\circ F(\chi))^\sharp(t_k)=j_i^\sharp(t_k)=f^i_k
  \]
  By the uniqueness, this implies $\vphi_{i'}\circ F(\chi)=\vphi_i$, so there exists a unique morphism $\vphi:X\to\aff^p$ \scth $\vphi\circ j_i=\vphi_i$ \fa $i$. Then $\smash{j_i^\sharp(\vphi^\sharp(t_k))=\vphi_i^\sharp(t_k)=f^i_k}$ for any $i$. This implies $\smash{\vphi^\sharp(t_k)=f_k}$, $k=1,\dotsc,p$, in view of \thmref{Rem}{indlim-functions}.
\end{proof}

In particular, one may consider infinite versions of the affine superspaces $\aff^p$.

\begin{Not}[aff-inf][infinite affine superspace $\aff^\infty$]
  For any integers $p'>p\sge0$, let $\smash{j_{p'p}:\aff^p\to\aff^{p'}}$ be the morphism induced by the map $x\mapsto(x,0_{p'-p})$. Then the family $(\aff^p,j_{p'p})$ is an inductive system of $\knums$-superspaces. Let $\aff^\infty$ denote its colimit.
  
  By \thmref{Prop}{indlim-regular}, $\aff^\infty$ is a regular superspace, and by \thmref{Rem}{indlim-functions}, there is an infinite sequence $t_1,t_2,\dotsc\in\Gamma(\sh O_{\aff^\infty,\Bbbk})$ of coordinate functions. It is determined by $j_p^\sharp(t_j)=t_j$ \fa $j\sle p$, where $j_p:\aff^p\to\aff^\infty$ are the natural morphisms.
\end{Not}

The proof of the following lemma is similar to that of \thmref{Prop}{indlim-regular}.

\begin{Lem}
  Let $X$ be a $\knums$-superspace. Then $X$ is regular (resp.~subregular) if and only if for each open subspace $U\subseteq X$, the following natural map is bijective (resp.~injective):
  \begin{equation}\label{eq:regular-ainfty}
    \aff^\infty(U)\to\prod_{k=1}^\infty\Gamma(\sh O_{U,\Bbbk,\ev}):\vphi\mapsto(\vphi^\sharp(t_1),\vphi^\sharp(t_2),\dotsc).
  \end{equation}
\end{Lem}

We now derive an explicit description of the infinitesimal normal neighbourhoods of the diagonal of $\aff^p$.

\begin{Not}
  Let $p\in\nats$ and $\Delta:\aff^p\to\aff^{2p}$ be the diagonal morphism. Denote by $\smash{(\aff^p)^{(n)}}$, $n\in\nats\cup\infty$, the order $n$ infinitesimal normal neighbourhood of $\aff^p$ w.r.t.~$\Delta$.
\end{Not}

Recall the Weil superalgebra $\sdual_n^{p|q}$ from \thmref{Ex}{super-dual}, and let $\sdual_\infty^{p|q}\defi\varprojlim_n\sdual_n^{p|q}$.

\begin{Prop}[aff-nbh][infinitesimal neighbourhoods of the diagonal]
  Let $p,n\in\nats$. We have $\smash{(\aff^p)^{(n)}}=\aff^p\times\Spec\sdual_n^{p|0}$, so that $\smash{\sh O_{\smash{(\aff^p)^{(n)}}}=\sh O_{\aff^p}\otimes\sdual_n^{p|0}}$. The latter statement also holds for $n=\infty$, so that $\smash{\sh O_{\smash{(\aff^p)^{(\infty)}}}=\sh O_{\aff^p}\llbracket T_1,\dotsc,T_p\rrbracket}$. 
\end{Prop}

\begin{proof}
  Let us first assume that $n$ is finite. We let $f_j\defi x_j-y_j$ for $j=1,\dotsc,p$ and $f_j\defi x_{j-p}$ for $j=p+1,\dotsc,2p$. Then the assertion follows immediately from \thmref{Cor}{hadamard-higherorder}. The superspaces $\smash{(\aff^p)^{(n)}}$ form an inductive system for fixed $p$. Since $\sh O_{\aff^p}\otimes\sdual_\infty^{p|0}$ is the projective limit of $\sh O_{\aff^p}\otimes\sdual_n^{p|0}$, for any $p$, the assertion follows.
\end{proof}

We now investigate when the regularity of $\knums$-superspaces is inherited by subspaces. To that end, consider the following definition.

\begin{Def}[local-factor][local factorisation]
  Let $f:X\to Y$, $h:X\to Z$ be morphisms of $\knums$-superspaces. We say that \Define[factorisation (of morphisms)!local]{$h$ locally factors through $f$} resp.~\Define[factorisation (of morphisms)!local unique]{locally factors uniquely} if there are open covers $X_i$, $Y_i$ of $X$, $Y$ \scth $X_i\subseteq f^{-1}(Y_i)$ and $h:X_i\to Z$ factors resp.~factors uniquely through $f:X_i\to Y_i$, for any $i$.
\end{Def}

\begin{Prop}[emb-regular][regularity of subspaces]
  Let $j:Y\to X$ be a $\Bbbk$-valued embedding of $\knums$-superspaces where $X$ is regular. The following are equivalent:
  \begin{enumerate}[wide]
    \item The $\knums$-superspace $Y$ is regular;
    \item morphisms from open subspaces of $Y$ to $\aff^p$ locally factor through $j$; and
    \item the $\knums$-superspace $Y$ is subregular.
  \end{enumerate}

  \noindent
  The equivalent conditions hold if $j:Y\to X$ is locally retractable.
\end{Prop}

The \emph{proof} uses the following lemma, which is based on ideas from the theory of ideals of $\sh C^\infty$-algebras, \vq \cite{moerdijk-reyes}*{Chapter I, Proposition 1.2}.

\begin{Lem}[emb-mor-aff-welldef]
  Let $X$ be regular, $\vphi:Y\to X$ and $\phi,\psi:X\to\aff^p$ be morphisms of $\knums$-superspaces. If $(\phi\circ\vphi)^\sharp(t_i)=(\psi\circ\vphi)^\sharp(t_i)$, for $i=1,\dotsc,p$, then $\phi\circ\vphi=\psi\circ\vphi$.
\end{Lem}

\begin{proof}
  Since $X$ is regular, we may define $\vrho:X\to\aff^{2p}$ by the requirement that $\vrho^\sharp(x_i)=\phi^\sharp(t_i)$ and $\vrho^\sharp(y_i)=\psi^\sharp(t_i)$. Similarly, since $\aff^{2p}$ is regular, there are unique morphisms $p_1,p_2:\aff^{2p}\to\aff^p$ \scth $\smash{p_1^\sharp(t_i)=x_i}$ and $\smash{p_2^\sharp(t_i)=y_i}$. Obviously, $p_1\circ\vrho=\phi$ and $p_2\circ\vrho=\psi$, where again, we have used the regularity. 

  Now, let $k\in\sh O_{\aff^p}(W)$ where $W\subseteq\aff^p_0$ is open. According to \thmref{Cor}{hadamard-diag}, there exist functions $k_i$ defined on some open neighbourhood of $\Delta_0(W)\subseteq\aff_0^p\times\aff_0^p$ with 
  \[
    p_1^\sharp(k)-p_2^\sharp(k)=\sum\nolimits_{i=1}^p(x_i-y_i)k_i.
  \]  
  We thus compute
  \[
     \phi^\sharp(k)-\psi^\sharp(k)=\vrho^\sharp(p_1^\sharp(k)-p_2^\sharp(k))=\sum\nolimits_{i=1}^p(\phi^\sharp(t_i)-\psi^\sharp(t_i))\vrho^\sharp(k_i),
  \]
  so $\vphi^\sharp\phi^\sharp(k)=\vphi^\sharp\psi^\sharp(k)$. Since $k$ was arbitrary, we have $\phi\circ\vphi=\psi\circ\vphi$.
\end{proof}

\begin{proof}[\protect{Proof of \thmref{Prop}{emb-regular}}]
  We may assume that $j$ is a closed embedding, and by passing to open subspaces, it is sufficient to consider the maps $\aff^p(Y)\to\Gamma(\sh O_{Y,\Bbbk,\ev})^p$. Let $f_1,\dotsc,f_p\in\Gamma(\sh O_{Y,\Bbbk,\ev})$. There exist an open cover $(U_i)$ of $X_0$ and even superfunctions $g_{ij}\in\sh O_{X,\ev,\Bbbk}(U_i)$, $j=1,\dotsc,p$, \scth $\smash{j^\sharp}(g_{ij})=f_j|_{V_i}$, where $V_i\defi j_0^{-1}(U_i)$.
  
  By the assumption on $X$, there are morphisms $\vphi_i:X|_{U_i}\to\aff^p$ \scth $\smash{\vphi_i^\sharp(t_j)=g_{ij}}$. In view of \thmref{Lem}{emb-mor-aff-welldef}, the morphisms $\vphi_k\circ j$ and $\vphi_\ell\circ j$ coincide on $Y|_{V_k\cap V_\ell}$. Hence, by \thmref{Prop}{relative-glue}, there is a morphism $\vphi:Y\to\aff^p$ \scth $\vphi=\vphi_i\circ j$ on $Y|_{V_i}$. In particular, we have $\vphi^\sharp(t_j)=f_j$. Hence, the canonical map $\aff^p(Y)\to\Gamma(\sh O_{Y,\Bbbk,\ev})^p$ is surjective for any $p$, and (i) and (iii) are equivalent.
  
  On the other hand, let $\psi:Y\to\aff^p$ be another morphism \scth $\psi^\sharp(t_j)=f_j$, and assume that $\psi$ locally factors through $j$. Possibly after refining the open cover $(U_i)$, $\psi:Y|_{V_i}\to\aff^p$ factors through $j:Y|_{V_i}\to X|_{U_i}$ to morphisms $\psi_i:X|_{U_i}\to\aff^p$. Then $h_{ij}\defi\psi_i^\sharp(t_j)$ satisfies $j^\sharp(h_{ij})=f_j|_{V_i}$. Applying \thmref{Lem}{emb-mor-aff-welldef} and \thmref{Prop}{relative-glue} again, we find that $\psi=\vphi$.
  
  Hence, if any morphism from an open subspace of $Y$ to $\aff^p$ locally factors through $j$, then $Y$ is regular. In particular, (ii) implies (i). Conversely, assume that $Y$ is regular, and let $\vphi:Y\to\aff^p$ be a morphism. (By the definition of regularity, we may assume that $\vphi$ is defined on all of $Y$.) Let $f_j\defi\vphi^\sharp(t_j)$ and $(U_i)$, $(V_i)$, $g_{ij}$ be as above. Define $\vphi_i:X|_{U_i}\to\aff^p$ by the condition $\smash{\vphi_i^\sharp(t_j)=g_{ij}}$. Then the regularity of $Y$ shows that $\vphi=\vphi_i\circ j$ on $Y|_{V_i}$, so that $\vphi$ locally factors through $j$. 

  Finally, if $j:X\to Y$ is locally retractable, then (ii) clearly holds.
\end{proof}

\begin{Cor}[aff-nbh-reg]
  Let $p,n\in\nats\cup\infty$. The $\knums$-superspaces $(\aff^p)^{(n)}$ are regular.
\end{Cor}

\begin{proof}
  In view of \thmref{Prop}{indlim-regular}, and because inductive limits commute, it is sufficient to prove this for $p$ and $n$ finite. Let $\sh I$ be the ideal of $\Delta$, so that the ideal of the morphism $\smash{\Delta^{(n)}:(\aff^p)^{(n)}}\to\aff^{2p}$ induced by $\Delta$ is $\sh I^{n+1}$.

  In view of \thmref{Prop}{emb-regular} and \thmref{Lem}{finht-canmap-inj}, it will be sufficient to prove that $(\aff^p)^{(n)}$ is tidy. But this is immediate from \thmref{Prop}{aff-nbh} and \thmref{Prop}{weil-tidy}.
\end{proof}

Here is a partial converse of the previous proposition.

\begin{Prop}[reg-thick][regularity of thickenings]
  Let $X$ and $Y$ be $\knums$-superspaces where $X$ is subregular, and $j:Y\to X$ be a locally retractable thickening of finite girth. If $Y$ is regular, then so is $X$.
\end{Prop}

\begin{proof}
  In view of \thmref{Prop}{relative-glue}, the question is local, so that we may assume, by passing to an open subspace, that $j$ admits a retraction $r:X\to Y$ and is of finite girth $N$. Then $X=Y^{(N)}$, where $Y^{(n)}$ denotes the normal infinitesimal neighbourhood of $Y$ in $X$ of order $n$.
  
  Let $f_1,\dotsc,f_p\in\Gamma(\sh O_{X,\Bbbk,\ev})$ and set $g_j\defi j^\sharp(f_j)$. Denote by $\vphi:Y\to\aff^p$ the unique morphism \scth $\vphi^\sharp(t_j)=g_j$, and set $\eps_j\defi r^\sharp(g_j)-f_j$. Since $j^\sharp(\eps_j)=0$, the $\eps_j$ are by assumption nilpotent elements of $\Gamma(\sh O_X)$. By \thmref{Prop}{aff-nbh}, we may define a morphism $\psi^\infty:X\to\smash{(\aff^p)^{(\infty)}}$ by setting $\psi^\infty_0\defi\vphi_0$ and for any $h\in\sh O_{\aff^p}(U)$
  \[
    \psi^{\infty\sharp}(hT^\alpha)\defi r^\sharp(\vphi^\sharp(h))\eps^\alpha,
  \]
  where we use the familiar multi-index notation.
  
  Let $r_\aff^n:(\aff^p)^{(n)}\to\aff^p$ be the morphism $\Delta^{(n)}:(\aff^p)^{(n)}\to\aff^{2p}$, composed with the projection $p_2:\aff^{2p}\to\aff^p$. We define the morphism $\psi\defi r_\aff^\infty\circ\psi^\infty:X\to\aff^p$. Then $\smash{r_\aff^{\infty\sharp}}(t_j)=y_j=x_j-1\cdot(x_j-y_j)$. Considered as a section of $\sh O_{\aff^p}\llbracket T_1,\dotsc,T_p\rrbracket$, the latter quantity equals $t_jT^0-T_j$. Hence,
  \[
    \psi^\sharp(t_j)=r^\sharp(g_j)-\eps_j=f_j,
  \]
  so $\psi$ meets the requirements. Since $X$ is subregular, it is regular.
\end{proof}

\begin{Cor}[weil-ext-reg][regularity of Weil thickenings]
  Let $X$ be a tidy regular $\knums$-superspace and $A$ a Weil $\knums$-superalgebra. Then the Weil thickened $\knums$-superspace $X^A$ is regular. In particular, $\Spec A$ itself and any Weil thickening of $\aff^p$ is regular.
\end{Cor}

\begin{proof}
  By \thmref{Prop}{weil-tidy}, $X^A$ is tidy. By \thmref{Lem}{finht-canmap-inj}, $X^A$ is subregular. Hence, \thmref{Prop}{reg-thick} applies, and the assertion follows, in view of \thmref{Prop}{manifold-mor}.
\end{proof}

\begin{Rem}
  Since by \thmref{Prop}{aff-nbh}, $(\aff^p)^{(n)}$ is a Weil thickening of $\aff^p$, \thmref{Cor}{weil-ext-reg} furnishes an alternative derivation of \thmref{Cor}{aff-nbh-reg}.

  Moreover, \thmref{Cor}{weil-ext-reg} implies that supermanifolds (to be defined below) are regular. We shall presently see an alternative proof (\thmref{Prop}{fg-countht}).
\end{Rem}

\subsection{Locally finitely generated superspaces}\label{sec:locally finitely generated superspaces}

The regular spaces we have considered so far have many desirable permanence properties under the assumption that subregularity is preserved. However, subregularity is not automatically preserved under thickenings and embeddings, and relatedly, (fibre) products do not exist.

We will now single out the full subcategory of locally finitely generated superspaces, which turns out to be better behaved on both accounts: it is stable under Weil thickenings and admits fibre products. 

\begin{Def}[fg-def][finitely generated superspaces]
  Let $X$ be a $\knums$-super\-space. If there exists a tidy embedding $j:X\to\mathbb A^{p|q}$, then $X$ is called \Define[superspace!finitely generated]{finitely generated}. It is called \Define[superspace!locally finitely generated]{locally finitely generated} if it admits a cover by finitely generated open subspaces.
\end{Def}

\begin{Prop}[fg-countht][finite generation implies tidiness]
  Let $X$ be a locally finitely generated $\knums$-superspace. Then $X$ is tidy, and hence, regular.
\end{Prop}

In the \emph{proof}, we need the following lemma for the case of $\Bbbk\neq\knums$. 

\begin{Lem}[cr-regular]
  Let $X$ be a formally Noetherian tidy $\reals$-superspace. If $X$ is regular over $(\reals,\reals)$, then $X_\cplxs$ is regular over $(\cplxs,\reals)$.
\end{Lem}

\begin{proof}
  Let $X_N$ be the tidying of the $N$th order infinitesimal normal neighbourhood of $j_{X_0}:X_0\to X$. \thmref{Prop}{emb-regular} applies, so that the $X_N$ are regular. By \thmref{Prop}{noether-tidy-approx}, $X=\varinjlim_NX_N$. Since complexification is a right adjoint functor by \thmref{Cons}{basefieldch}, it preserves colimits. Thus, in view of \thmref{Prop}{indlim-regular}, we may assume, without loss of generality, that $X$ has finite girth $N$ (say). Moreover, we may assume that $X$ is even. 

  Let $Y\defi X_\cplxs$. By \thmref{Lem}{finht-canmap-inj}, $Y$ is subregular over $(\cplxs,\reals)$. On the other hand, we have $\sh O_{Y,\reals}=\sh O_X\oplus i\,\sh N_X$, so for $f=(f_a)\in\Gamma\Parens0{\sh O_{Y,\reals}^p}$, we may decompose $f=g+h$ where $g_a\in\Gamma(\sh O_X)$ and $h_a\in\Gamma(i\,\sh N_X)$. Let $\vphi:X\to\aff^p$ be the unique morphism of $\reals$-superspaces \scth $\vphi^\sharp(t_a)=g_a$, where we denote by $\aff^p$ the $p$-dimensional $(\reals,\reals)$-affine superspace. 

  Define $\psi:Y\to(\aff^p)_\cplxs$ by $\psi_0\defi\vphi_0$ and $\psi^\sharp(k)\defi\sum_{\Abs0\alpha\sle N}\frac1{\alpha!}\vphi^\sharp\Parens1{\frac{\partial^{\Abs0\alpha}k}{\partial t^\alpha}}h^\alpha$. By the Leibniz rule, and because the girth of $Y$ is $N$, this indeed defines a morphism, and clearly $\psi^\sharp(t_a)=g_a+h_a=f_a$, so $Y$ is regular over $(\cplxs,\reals)$.
\end{proof}

\begin{proof}[\protect{Proof of \thmref{Prop}{fg-countht}}]
  We may assume without loss of generality that there is an embedding $j:X\to\aff^{p|q}$ with ideal $\sh I$. By \thmref{Prop}{noetherianjets-countht}, $X$ is tidy, since this is a local property. In particular, on applying \thmref{Lem}{finht-canmap-inj}, we see that $X$ is subregular. 

  To show that $X$ is in fact regular, let us first assume that $\knums=\Bbbk$, so that any embedding is $\Bbbk$-valued. In this case, \thmref{Prop}{emb-regular} applies, so that by \thmref{Cor}{weil-ext-reg}, we find that $X$ is regular, as claimed. 

  To prove the claim in general, only the case of $(\knums,\Bbbk)=(\cplxs,\reals)$ remains to be considered. Let $Y$ be the affine superspace over $(\reals,\reals)$ of dimension $p|q$, so that $Y_\cplxs=\aff^{p|q}$. Let $\sh I_r\defi\sh I\cap\sh O_Y$ and $j_r:X_r\to Y$ be the embedding of $\reals$-superspaces defined by this ideal. Clearly, the ideal $\sh I_r$ is tidy, so that the $\reals$-superspace $X_r$ is $(\reals,\reals)$-regular by the arguments above. 

  By \thmref{Lem}{weilext-noetherianjets} and \thmref{Lem}{noetherian-emb}, $X_r$ is formally Noetherian, so \thmref{Lem}{cr-regular} applies, and we see that $(X_r)_\cplxs$ is regular. Plainly, there exists a natural surjection $\sh O_{X_r}\otimes\cplxs\to\sh O_X$, and this induces a thickening $X\to(X_r)_\cplxs$. On applying \thmref{Prop}{emb-regular} to it, we arrive by the desired conclusion.
\end{proof}

\begin{Rem}
  \thmref{Prop}{fg-countht} implies in particular that supermanifolds (to be defined below) are regular, a classical result due to Leites \cites{deligne-morgan,leites,vsv}. Note that it furnishes an alternative derivation thereof. Leites's method of proof is closer to that followed in \thmref{Prop}{reg-thick} and \thmref{Cor}{weil-ext-reg}.
\end{Rem}

\begin{Cor}[fg-emb][finite generation and affine embeddings]
  Let $X$ be finitely generated. Then $X\times\Spec A$ admits a tidy embedding into some $\aff^{p|q}$, for any Weil $\knums$-superalgebra $A$.
\end{Cor}

\begin{proof}
  It suffices to show that any $\aff^r\times\Spec A$, where $A$ is a Weil $\knums$-superalgebra, can be embedded into some $\aff^{p+r|q}$.
  
  Since $A$ is finitely generated by nilpotent elements, there exists a surjective homomorphism from $B\otimes C$ onto $A$, where $C\defi \knums[\Theta_1,\dotsc,\Theta_q]$ and
  \[
    B\defi\knums[T_1,\dotsc,T_p]/(T_1^{\alpha_1}\dotsm T^{\alpha_p}_p|\alpha_1+\dotsm+\alpha_p=N),
  \]
  for some non-negative integers $p$, $q$, and $N$, and even resp.~odd indeterminates $T_i$ resp.~$\Theta_j$. We may assume $A=B\otimes C$. Since $\Spec C=\aff^{0|q}$, we may assume $A=B$.
  
  Define $\vphi:\aff^r\times\Spec A\to\aff^{p+r}$ by letting $\vphi_0\defi{\id_{\Bbbk^r}}\times 0$ and $\vphi^\sharp$ be the canonical projection $\sh O_{\aff^{p+r}}\to\sh O_{\aff^{p+r}}/\sh I^N$ where $\sh I(U)\defi\sh O_{\aff^{p+r}}(U)$ if $U\cap(\Bbbk^r\times0)=\vvoid$ and 
  \[
    \sh I(U)\defi\Set1{f\in\sh O_{\aff^{p+r}}(U)}{f(U\cap(\Bbbk^r\times0))=0}
  \]
  otherwise. Then by \thmref{Cor}{hadamard-higherorder}, $\vphi_0^{-1}(\sh O_{\aff^{p+r}}/\sh I^N)\cong\sh O_{\aff^r}\otimes A$. Thus, $\vphi$ is indeed a well-defined embedding $\aff^r\times\Spec A\to\aff^{p+r}$.
\end{proof}

\begin{Cor}[fg-weilext][locally finite generation of Weil thickenings]
  The category of (locally) finitely generated $\knums$-superspaces is stable under Weil thickenings.
\end{Cor}

\begin{proof}
  Let $\vphi:X\to\aff^p\times\Spec A$ be a tidy embedding and $B$ a Weil $\knums$-superalgebra. Applying the functor $(-)^B=(-)\times\Spec B$, we obtain by \thmref{Prop}{weil-extended} a morphism $\vphi^B:X^B\to\aff^p\times\Spec(A\otimes B)$. Evidently, it is an embedding. By \thmref{Prop}{fg-countht}, $X$ is tidy, hence, so is $X^B$, by virtue of \thmref{Prop}{weil-tidy}. Then \thmref{Cor}{tidy-emb} shows that $\vphi^B$ is tidy.
\end{proof}

\begin{Cor}[fg-char][characterisation of finite generation]
  Let $X$ be a $\knums$-superspace. The following are equivalent:
  \begin{enumerate}[wide]
    \item\label{fg-char-i} $X$ is locally finitely generated;
    \item\label{fg-char-ii} $X$ is tidy and admits an open cover by subspaces embeddable into $\aff^p\times\Spec A$, for some integers $p$ and Weil $\knums$-superalgebras $A$;
    \item\label{fg-char-iii} $X$ is tidy and admits an open cover by subspaces embeddable into $\aff^{p|q}$, for some integers $p$ and $q$.
  \end{enumerate}
\end{Cor}

\begin{proof}
  We have \eqref{fg-char-i} $\Rightarrow$ \eqref{fg-char-ii} by \thmref{Prop}{fg-countht}, and \eqref{fg-char-ii} $\Rightarrow$ \eqref{fg-char-iii} by \thmref{Cor}{fg-emb}. Finally, \eqref{fg-char-iii} $\Rightarrow$ \eqref{fg-char-i} follows from \thmref{Cor}{tidy-emb}.
\end{proof}

\begin{Prop}[fg-redbody][finite generation of reduction, body, and even part]
  Let $X$ be locally finitely generated. Then so are $X_0$, $X_\ev$, and $X^\ev$.
\end{Prop}

\begin{proof}
  We may assume that $X$ is finitely generated, so that, taking \thmref{Cor}{fg-emb} into consideration, there is a tidy embedding $\vphi:X\to\aff^{p|q}$. Then $(\aff^{p|q})_0=(\aff^{p|q})_\ev=\aff^p$ is reduced and, in particular, even.

  Since $\vphi\circ j_{X_\ev}=j_{(\aff^{p|q})_\ev}\circ\vphi_\ev$ is an embedding, $\vphi_\ev:X_\ev\to\aff^p$ is an embedding, by \thmref{Prop}{redbodev-emb}. Similarly, so is $\vphi_0:X_0\to\aff^p$. 

  Since $X_0$ is reduced, it is tidy, and so finitely generated, by \thmref{Cor}{fg-char}. Reasoning similarly, it remains to be shown that the identity of $X_\ev$ is tidy, and to that end, it is sufficient that the embedding $\vphi_\ev$ be tidy, by \thmref{Prop}{noetherianjets-countht}.

  So, let $\sh I$ denote the vanishing ideal of $\vphi$ and $\sh J$ that of $\vphi_\ev$. Fix $f\in\sh O_{\aff^p}(U)$ \scth $f_x\in\sh J_x+\ger m_{\aff^p,x}^N$ \fa $x\in U\cap\vphi_0(X_0)$, $N\in\nats$. We may consider local sections $\sh O_{\aff^p}$ as local sections of of $\sh O_{\aff^{p|q}}$, as in \thmref{Not}{std-coord}. The sheaf map $\sh O_{\aff^p}\to\sh O_{\aff^{p|q}}$ considered there is local on the stalks. Thus, if $f_x=g_x+h_x$ where $g_x\in\sh J_x$ and $\smash{h_x\in\ger m^N_{\aff^p,x}}$, then $\smash{h_x\in\ger m_{\aff^{p|q},x}^N}$. Moreover, $\smash{\vphi^\sharp(g_x)=\vphi_\ev^\sharp(g_x)=0}$, so that $f_x\in\smash{\sh I_x+\ger m_{\aff^{p|q},x}^N}$ \fa $x\in U\cap\vphi_0(X_0)$. By the assumption on $\vphi$, this implies that $\smash{\vphi^\sharp(f)=0}$, and in particular, $\smash{\vphi_\ev^\sharp(f)=0}$, proving the claim.

  Finally, we consider the case of the even part. Firstly, let $A$ be the even part of the Grassmann algebra $\sh O_{\aff^{0|M}}$. Then
  \[
    A=\knums\Bracks1{t_{ij}\bigm|1\sle i<j\sle M}/I\ ,\mathtxt{where}I\defi\Parens1{t_{ij}t_{k\ell}\bigm|\{i,j\}\cap\{k,\ell\}\neq\vvoid}.
  \]
  This is an even Weil $\knums$-algebra, and in case $X=\aff^{N|M}$, we have $\smash{X^\ev}=\smash{\aff^N}\times\Spec A$. Applying \thmref{Cor}{fg-weilext}, we find that $X^\ev$ is finitely generated in this case. It follows from the definitions that $\smash{X^\ev}$ is tidy if $X$ is. Thus, according to the above statements about embeddings, \thmref{Cor}{tidy-emb}, and the definitions, $X^\ev$ is locally finitely generated if $X$ is.
\end{proof}

\begin{Prop}[fg-product][products of finitely generated superspaces]
  Let $X$ and $Y$ be (locally) finitely generated. Then $X\times Y$ exists in the category of tidy regular superspaces, and is (locally) finitely generated.
\end{Prop}

\begin{proof}
  The question is local. So, passing to open subspaces, assume that embeddings $i:X\to\aff^p\times\Spec A$ and $j:Y\to\aff^r\times\Spec B$ are given, with the corresponding vanishing ideals $\sh I$ and $\sh J$. Clearly,
  \[
    W\defi\aff^p\times\Spec A\times\aff^r\times\Spec B\cong\aff^{p+r}\times\Spec (A\otimes B)
  \]
  exists in the category of regular superspaces. Let $U\subseteq\aff^p_0\times\aff^r_0$ be open such that $(i_0\times j_0)(X_0\times Y_0)$ is closed in $U$, and $\sh K$ be the ideal of $\sh O_W|_U$ generated by $p_1^\sharp(\sh I)$ and $p_2^\sharp(\sh J)$. Define the $\knums$-superspace $Z\defi(Z_0,\sh O_Z)$ by $Z_0\defi X_0\times Y_0$ and $\sh O_Z\defi(i_0\times j_0)^{-1}\Parens1{\sh O_W|_U/\sh K}$, together with the canonical embedding $k:Z\to W$. Let $t:Z^\circ\to Z$ be the tidying of $Z$. Then $k^\circ\defi k\circ t$ is an embedding and $Z^\circ$ is finitely generated by \thmref{Cor}{fg-char}.
      
  We have $(p_1\circ k^\circ)^\sharp(\sh I)=0$ and $(p_2\circ k^\circ)^\sharp(\sh J)=0$. Hence, the morphisms $p_1\circ k^\circ$ and $p_2\circ k^\circ$ factor uniquely through $i$ and $j$, respectively, giving rise to morphisms $\pi_1:Z^\circ\to X$ and $\pi_2:Z^\circ\to Y$.
  
  Let $R$ be a tidy regular superspace, and $\phi:R\to X$, $\psi:R\to Y$ be morphisms. Since $R$ is regular, $\vphi\defi(i\circ\phi,j\circ\psi):R\to W$ exists. Since $(p_1\circ\vphi)^\sharp(\sh I)=0$ and $(p_2\circ\vphi)^\sharp(\sh J)=0$, we have $\vphi^\sharp(\sh K)=0$. Moreover, $\vphi_0(R_0)\subseteq k_0(Z_0)$, so $\vphi$ factors uniquely through $k$, and since $R$ is tidy, it factors uniquely through $k^\circ$, by \thmref{Prop}{tidying}. Denote the resulting morphism $R\to Z^\circ$ by $\vrho$.
    
  We have $i\circ\pi_1\circ\vrho=p_1\circ k^\circ\circ\vrho=p_1\circ\vphi=i\circ\phi$, so that $\pi_1\circ\vrho=\phi$. Similarly, $\pi_2\circ\vrho=\psi$. Let $\sigma:R\to Z^\circ$ be another morphism \scth $\pi_1\circ\sigma=\phi$ and $\pi_2\circ\sigma=\psi$. Then $k^\circ\circ\sigma=(i\circ\phi,j\circ\psi)=k^\circ\circ\vrho$, and it follows that $\sigma=\vrho$. 
\end{proof}

\begin{Rem}[aff-prod-lfg]
  Given arbitrary finite-dimensional super-vector spaces $V$ and $V'$ over $(\knums,\Bbbk)$, the affine superspace $\aff(V\oplus V')$ is the direct product of $\aff(V)$ and $\aff(V')$ in the category of locally finitely generated $\knums$-superspaces. The superspaces in question are locally finitely generated, so the assertion is immediate from \thmref{Cor}{aff-prod}.
\end{Rem}

\begin{Cor}[tidy-finlimits][limits of finitely generated superspaces]
  The category of locally finitely generated $\knums$-superspaces is finitely complete, \ie all finite limits exist.
\end{Cor}

\begin{proof}
  Since $\aff^p$ is Hausdorff, any locally finitely generated $\knums$-superspace is locally Hausdorff. Hence, the equaliser $Z$ of morphisms $\phi,\psi:X\to Y$ locally finitely generated $\knums$-superspaces exists in the category of tidy $\knums$-superspaces, owing to \thmref{Prop}{tidy-equaliser}.
  
  Such a $Z$ is tidy, and locally admits embeddings into $\aff^p\times\Spec A$, since so does $X$. These embeddings are necessarily tidy, by \thmref{Cor}{tidy-emb}, so $Z$ is locally finitely generated. The claim follows from \thmref{Prop}{fg-product} and standard facts \cite{maclane}.
\end{proof}

The corollary allows to introduce the following notation. 

\begin{Not}[fibprod][Fibre products and fibres]
  Let $\vphi:X\to Z$ and $\psi:Y\to Z$ be morphisms of locally finitely generated superspaces. By virtue of \thmref{Cor}{tidy-finlimits}, the fibred product $X\times_ZY$ exists in the category of locally finitely generated superspaces.
  
  In particular, if $\psi$ is the inclusion of a point $z\in Z_0$ ($Y=*$ and $\psi_0(*)=z$), the fibre product is denoted by $\vphi^{-1}(z)$ or $X_z$. It is called the \Define[fibre!of a morphism]{fibre at $z$} of $\vphi$ (or of $X/Z$).
\end{Not}

A useful fact is that embeddings are preserved by fibre products.

\begin{Cor}[emb-fibprod][fibre products of embeddings]
  Let $j_1:X_1/S\to Y_1/S$ and $j_2:X_2/S\to Y_2/S$ be (open) embeddings in $\ssplfg{S}$. Then the induced morphism $j_1\times_Sj_2:X_1\times_SX_2\to Y_1\times_SY_2$ is an (open) embedding. 

  If $j_1$ and $j_2$ are closed and $S$ is Hausdorff, then it is closed. In this case, if $\sh J_1$ denotes the ideal of $j_1$ and $\sh J_2$ that of $j_2$, then the ideal of $j_1\times_Sj_2$ is the tidying of the ideal generated by $p_1^\sharp(\sh J_1)$ and $p_2^\sharp(\sh J_2)$.
\end{Cor}

\begin{proof}
  Firstly, composites of embeddings are embeddings. Secondly, we observe that if $f\circ g$ is an (open resp.~closed) embedding, then so is $g$, a fact that follows directly from the definitions. Since equalisers in $\ssplfg{\knums}$ are embeddings by \thmref{Prop}{tidy-equaliser}, we may reduce to the case of $S=*$.

  Direct products of topological embeddings are embeddings. By the construction of products, this reduces the question to a local one. Therefore, we may assume that $Y_1=\aff^{p_1|q_1}$, $Y_2=\aff^{p_2|q_2}$. But then the claim follows by construction.
\end{proof}

We note another corollary to \thmref{Prop}{fg-product}.

\begin{Cor}[redbody-lim][reduction, body, and even parts of (co)limits]
  The reduction and body functors preserve limits in $\ssplfg{\knums}$. The even part functor preserves colimits. 
\end{Cor}

\begin{proof}
  By \thmref{Prop}{fg-redbody}, all three functors are endofunctors of $\ssplfg{\knums}$. Reduction is right adjoint to the inclusion of the full subcategory of reduced superspaces in $\ssplfg{\knums}$ by \thmref{Prop}{red-nat}. The body functor is right adjoint to the inclusion of the full subcategory of even superspaces in $\ssplfg{\knums}$ by \thmref{Prop}{even-nat}. Finally, the even part functor is left adjoint to that same inclusion by \thmref{Prop}{ev-part-func}. Since right (resp.~left) adjoints preserve limits (resp.~colimits), the assertion follows.
\end{proof}

We end this subsection by a discussion of morphisms and coordinates in the framework we have introduced. We will discuss the meaning of morphisms from Weil thickened superspaces to a regular superspace with a global chart.

\begin{Scho}[mor-weil][morphisms and products]
  Let $S$ be a locally finitely generated $\knums$-super\-space, and $A$ a Weil $\knums$-superalgebra. Consider the Weil thickened superspace $S^A=S\times\Spec A$ with structure sheaf $\sh O_{S^A}=\sh O_S\otimes A$, defined in \thmref{Cons}{weil-extended-superspace}. By \thmref{Cor}{fg-weilext}, $S^A$ is again locally finitely generated. Recall that $(S^A)_0=S_0$. On one hand, 
  \[
    j_{(S^A)_0}^\sharp(f\otimes1)=j_{S_0}^\sharp(f)\mathfa f\in\sh O_S(U);
  \]
  on the other hand, $j_{\smash{(S^A)_0}}^\sharp=0$ on $\sh O_S\otimes\ger m$.

  Let $X$ be a regular $\knums$-superspace with a global system of coordinates $x=(x_a)$ with $a=1,\ldots,p+q=n$. Then $X$ is isomorphic to an open subspace of $\aff^{p|q}$. By \thmref{Prop}{mor-coord}, the set $X(S^A)$ of all $\vphi:S^A\to X$ is in natural bijection with
  \[
    \Set2{s=(s_a)\in(\Gamma(\sh O_S)\otimes A)^n}{\Abs0{s_a}=\Abs0{x_a}\,,\,s(S_0)\subseteq x(X_0)},
  \]
  where $n\defi p+q$. Under the bijection, $\vphi$ maps to $s$ if and only if $s=\vphi^\sharp(x)$.

  Now, let $m\defi k+\ell$, where $\dim\ger m=k|\ell$, and choose a homogeneous basis $e^1,\dotsc,e^m$ of $\ger m$. For any index $a$, we may write
  \[
    s_a=s_a^0\otimes 1+\sum_{b=1}^ms_a^b\otimes e^b,
  \]
  with unique $s_a^b\in\Gamma(\sh O_S)$, $b=0,\dotsc,m$. Since $j_{(S^A)_0}^\sharp(s_a^b\otimes e^b)=0$ for $b>0$, we have
  \begin{equation}\label{eq:weilthick-points}
    X(S^A)\cong\Set1{(s_a^b)}{\Abs0{s_a^b}=\Abs0{x_a}+\Abs0{e^b},s^0(S_0)\subseteq x(X_0)},   
  \end{equation}
  where we agree to write $e^0\defi1$, and $s^b\defi(s_a^b)_{a=1,\dotsc,n}$ for $b=0,\dotsc,m$. Observe that here, the mapping condition only concerns $s^0$.

  From now, we assume that either $\knums=\Bbbk$ or $pk+q\ell=0$. The latter condition means that either $X$ is purely even and $\ger m$ is purely odd, or \emph{vice versa}, so that all of the $s_a^b$ with $b>0$ are odd. Both conditions guarantee that the $s_a^b$ are $\Bbbk$-valued. 

  Applying \thmref{Prop}{mor-coord}, the tuple $s^0$ corresponds to a morphism $\smash{\vphi^0}:S\to X$, and $s^b$, for $b>0$, corresponds to a morphism $\smash{\vphi^b:S\to\aff^{p|q}}$ if $\Abs0{e^b}=\ev$, and to a morphism $\smash{\vphi^b:S\to\aff^{q|p}}$ otherwise. The correspondences are given \via
  \[
    s^0=(\vphi^0)^\sharp(x)\nd s^b_a=\begin{cases}(\vphi^b)^\sharp(t)&\Abs0{e^b}=\ev,\\(\vphi^b)^\sharp(\Pi t)&\Abs0{e^b}=\odd.\end{cases}
  \]
  Here, $t=(t_a)$ denotes the standard coordinate system of $\smash{\aff^{p|q}}$, and $\Pi t$ denotes the standard coordinate system of $\smash{\aff^{q|p}}$.

  In summary, we have for $\dim X=p|q$ and $\dim A=(k+1)|\ell$ a bijection
  \begin{equation}\label{eq:weil-ext-mor}
    X(S^A)\cong X(S)\times\aff^{p|q}(S)^k\times\aff^{q|p}(S)^\ell=\Parens1{X\times\aff^{kp+\ell q|kq+\ell p}}(S)
  \end{equation}
  which depends on the choice of a (global) coordinate system on $X$.

  In this form, the bijection, however, also depends on the choice of a homogeneous basis of $\ger m$; only the first component is entirely canonical.

  To give a bijection, which only depends on the choice of the coordinate system $(x_a)$ on $X$, consider the affine superspace $\smash{\aff(\ger m^{p|q})}$ associated with the super-vector space $\ger m^{p|q}\defi\ger m\otimes\knums^{p|q}$. The data $(s_a^b)_{b>0}$ correspond to a morphism $\smash{\vphi^+:S\to\aff(\ger m^{p|q})}$ \via $s_a^b=\vphi^{+\sharp}(e_b\otimes t_a)$. We obtain a bijection
  \begin{equation}\label{eq:weil-ext-mor2}
    X(S^A)\cong X(S)\times\aff\Parens1{\ger m^{p|q}}(S)=\Parens1{X\times\aff\Parens1{\ger m^{p|q}}}(S),
  \end{equation}
  where $\vphi$ and $(\vphi^0,\vphi^+)$ are related by
  \begin{equation}\label{eq:weil-adj}
    ({\id}\otimes\eps)(\vphi^\sharp(f))=\vphi^{0\sharp}(f),\quad
    ({\id}\otimes\mu)(\vphi^\sharp(x_a))=\vphi^{+\sharp}(\mu\otimes t_a); 
  \end{equation}
  here, $\eps:A\to\knums$ is the unique algebra morphism, and $f\in\sh O_X(U)$, $\mu\in\ger m^*$ are arbitrary. The obtained bijection is natural in $S$ and the pairs $(X,(x_a))$.
\end{Scho}

As we have seen, the second half of \thmref{Scho}{mor-weil} is valid only in the cases of $\knums=\Bbbk$ or $pk+q\ell=0$. Retaining the above notation, let us discuss why this is so. 

On one hand, morphisms $S^A\to X$ correspond to tuples $(s_a^b)$ of superfunctions where $s_a^0$ are $\Bbbk$-valued, whereas the $s_a^b$ are arbitrary. On the other hand, for any super-vector space $V$ over $(\knums,\Bbbk)$, morphisms $S\to X\times\aff(V)$ correspond to tuples of superfunctions, each of which is $\Bbbk$-valued. The distinction vanishes when $\knums=\Bbbk$, but otherwise (in the \emph{cs} case), it is an indelible fact (if $V$ is sufficiently general). Although not widely noticed, it has been recorded in the literature, see Ref.~\cite{deligne-morgan}*{Example 4.9.5}.

However, using the language of functors, this can be overcome, at the expense of leaving the category $\ssplfg{\knums}$. We now present this in detail.

\begin{Scho}[mor-weil-k][morphisms and products in the \emph{cs} case]
  Retain the notation from \thmref{Scho}{mor-weil}. For any finite-dimensional super-vector space $V$ over $\knums$, we define the set-valued cofunctor $\aff^\knums(V)$ on $\ssplfg{\knums}$ by 
  \begin{equation}\label{eq:aff-k}
    \aff^\knums(V)(T)\defi\Gamma\Parens1{(\sh O_T\otimes V)_\ev},\quad T\in\ssplfg{\knums}.   
  \end{equation}
  We call $\aff^\knums(V)$ the \Define{$\knums$-valued affine superspace functor}. Whenever $\knums=\Bbbk$, then $\aff^\knums(V)$ is represented by $\aff(V)$, by \thmref{Prop}{mor-coord}. Otherwise, it is in general not representable in $\ssplfg{\knums}$. By the discussion in \thmref{Scho}{mor-weil} ending in Equation \eqref{eq:weilthick-points}, we have 
  \begin{equation}\label{eq:weil-ext-mor-k}
    X(S^A)\cong X(S)\times\Gamma\Parens1{(\sh O_S\otimes\ger m^{p|q})_\ev}=\Parens1{X\times\aff^\knums(\ger m^{p|q})}(S)    
  \end{equation}
  The product is in the category $\smash{\Bracks1{\Parens1{\ssplfg{\knums}}^{\mathrm{op}},\Sets}}$ of set-valued cofunctors on $\ssplfg{\knums}$. Moreover, recall that $\ger m$ denotes the maximal ideal of $A$ and $\ger m^{p|q}=\ger m\otimes_\knums\knums^{p|q}$. 

  Explicitly, the bijection is given by mapping 
  \[
    \vphi\mapsto\Parens1{\vphi^0,\vphi^{+1},\dotsc,\vphi^{+(p+q)}}
  \]
  where $\vphi^{+a}\in\Gamma((\sh O_S\otimes\ger m)_\ev)$ are \scth
  \begin{equation}\label{eq:weil-adj-k}
    \vphi^\sharp(x_a)=s_a^0\otimes1+\vphi^{+a},\quad\vphi^{0\sharp}(x_a)=s_a^0.
  \end{equation}
  The bijection in Equation \eqref{eq:weil-ext-mor-k} generalises that in Equation \eqref{eq:weil-ext-mor2}, in that for $\knums=\Bbbk$, it coincides with the latter. In any case, it is natural in $S$ and the pairs $(X,(x_a))$ of superspaces with global coordinate systems. 
\end{Scho}

\subsection{Relative supermanifolds}\label{sec:relative supermanifolds}

In this section, we can finally introduce relative supermanifolds. In what follows, let $S$ be a locally finitely generated $\knums$-superspace. We denote by $\smash{\SSp_S^{\mathrm{lfg}}}$ the full subcategory of $\SSp_S$ consisting of all $X/S$ which are locally finitely generated as $\knums$-superspaces. In the case of $S=*$, we write $\SSp^{\mathrm{lfg}}$.
  
  In this, we are assuming that $\knums$, $\Bbbk$ and $\varpi$ are understood. If this is not the case, we will write $\smash{\SSp^{\varpi,\mathrm{lfg}}_{\knums,\Bbbk,S}}$ and $\smash{\SSp^{\varpi,\mathrm{lfg}}_{\knums,\Bbbk}}$, respectively.

\begin{Def}[rel-aff][relative affine superspace]
  We let $\aff^{p|q}_S\defi S\times\aff^{p|q}$, the product of locally finitely generated superspaces. Let $t=(t_a)$ be the standard coordinate system on $\aff^{p|q}$. By abuse of notation, the collection $\smash{(p_2^\sharp(t_a))}$ of superfunctions on $\smash{\aff^{p|q}_S}$ will again be denoted by $t=(t_a)$.
  
  We call this system of superfunctions the \Define{standard fibre coordinate system} and say it is in \Define[standard fibre coordinate system!standard order]{standard order} if the original system of superfunctions on $\aff^{p|q}$ was. We will also denote this by $(t,\theta)$ if we wish to see the parity explicitly.
\end{Def}

\begin{Prop}[rel-leites][morphisms with relative affine target]
  Let $X/S$ be a relative $\knums$-superspace in $\smash{\SSp_S^{\mathrm{lfg}}}$. There is a natural bijection
  \[
    \Hom[_S]1{X,\aff^{p|q}_S}\to\Set1{(\vphi_a)\in\Gamma(\sh O_X)^{p+q}}{\forall a:\Abs0{\vphi_a}=\Abs0{t_a}}:\vphi\mapsto\vphi^\sharp(t).
  \]
\end{Prop}

\begin{proof}
  Since $X$ is regular, it suffices to remark that there is a natural bijection $\Hom[_S]0{X,S\times Y}=\Hom0{X,Y}$ for locally finitely generated $\knums$-superspaces $Y$.
\end{proof}

\begin{Def}[fib-coord][fibre charts and fibre coordinate systems]
  Let $X/S$ be a relative $\knums$-superspace in $\smash{\SSp_S^{\mathrm{lfg}}}$ and $U\subseteq S$ an open subspace. By \thmref{Prop}{open-fibreprod}, $X\times_SU$ exists and is an open subspace of $X$. Let $V$ be an open subspace of $X\times_SU$.
  
  A \Define[fibre chart!local]{local fibre chart} (defined on $V/U$) is then by definition an open embedding $\smash{\vphi:V/U\to\aff^{p|q}_U/U}$. If $V=X$, then we call $\vphi$ a \Define[fibre chart!global]{global fibre chart}. The tuple $p|q$ is called the \Define[fibre chart!fibre dimension (graded)]{(graded) fibre dimension} of $\vphi$.
  
  Let $\vphi$ be a local fibre chart. The tuple $(x,\xi)\defi\vphi^\sharp(t,\theta)$ (also denoted by $x=(x_a)\defi\vphi^\sharp(t_a)$) is called a system of \Define[fibre coordinate system!local]{local fibre coordinates} (defined on $V/U$). If $\vphi$ is a global chart, we call it a system of \Define[fibre coordinate system!global]{global fibre coordinates}. Due to \thmref{Prop}{rel-leites}, fibre charts and fibre coordinate systems are in bijection. We will also call $p|q$ the \Define[fibre coordinate system!dimension (graded)]{(graded) dimension} of $(x,\xi)$.
  
Given a system of local fibre coordinates $(x,\xi)$ defined on $V/U$, any $f\in\sh O_X(V_0)$ has a unique representation
\begin{equation}\label{eq:fib-coord-repn}
  f=\sum\nolimits_If_I\xi^I\,,\mathtxt{where}\xi^I\defi\xi_{i_1}\dotsm\xi_{i_k}
\end{equation}
and $f_I=x^\sharp(g_I)$ \fs $g_I\in\sh O_{\smash{\aff^p_U}}(\vphi_0(V_0))$. Here, by abuse of notation, $x$ denotes the morphism $V/U\to\aff^p_U/U$ determined by $x^\sharp(t_j)=x_j$.
\end{Def}

\thmref{Prop}{rel-leites} and \thmref{Def}{fib-coord} at once give the following.

\begin{Cor}[mor-coord][relative morphisms \vs fibre coordinates]
  Let $X/S$ and $Y/S$ be in $\smash{\SSp_S^{\mathrm{lfg}}}$ and $y=(y_a)$ a global fibre coordinate system in standard order on $Y$. The following map is a bijection:
  \[
    \Hom[_S]0{X,Y}\to\Set1{x\in\Gamma(\sh O_X)_\ev^p\times\Gamma(\sh O_X)_\odd^q}{x(X_0)\subseteq y(Y_0)}\colon\vphi\mapsto \vphi^\sharp(y).
  \]
\end{Cor}

\begin{Def}[rel-smf][relative supermanifolds]
  Let $X/S$ be in $\SSp_S^{\mathrm{lfg}}$. If, for every point $x\in X_0$, there exist an open subspace $V/U$ of $X/S$ \scth $x\in V_0$ and a local fibre chart of $X/S$ defined on $V/U$, then $X/S$ is called a \Define[super-premanifold!over $S$]{super-premanifold} over $S$. If in addition, $X_0$ is a Hausdorff topological space, then $X/S$ is called a \Define[supermanifold!over $S$]{supermanifold over $S$}. In case $S=*$, $X$ is simply called a \Define{supermanifold}. 
  
   By reason of \thmref{Cor}{man-emb}, the graded dimensions of all local fibre charts at $x\in X_0$ are identical. We denote their common value by $\dim_{S,x}X$ and call this the \Define[supermanifold!graded fibre dimension at $x$]{(graded) fibre dimension at $x$}. For $S=*$, we simply write $\dim_xX$ and speak of the \Define[supermanifold!graded dimension at $x$]{(graded) dimension at $x$}. If $\dim_{S,x}X$ (resp.~$\dim_xX$) is independent of $x$, then we say that $X$ has \Define[supermanifold!pure (fibre) dimension]{pure (fibre) dimension}.
  
  The full subcategory of $\SSp_S$ whose objects are the supermanifolds over $S$ is denoted by $\SMan_S$. For $S=*$, we write $\SMan$. If $\knums$, $\Bbbk$, and $\varpi$ are not understood, we write $\SMan_{\knums,\Bbbk,S}^\varpi$ and $\SMan_{\knums,\Bbbk}^\varpi$, respectively.
\end{Def}

  If $S$ is a supermanifold, then in the literature \cites{leites,deligne-morgan}, supermanifolds over $S$ have been previously considered, under the name of \Define{families of supermanifolds}. 

\begin{Lem}
  Let $X/S$ be a super-premanifold over $S$. When $S$ is a super-premanifold, then so is $X$, and $\dim_xX=\dim_sS+\dim_{S,x}X$, for $p_{X,0}(x)=s\in S_0$.
\end{Lem}

\begin{proof}
  The question is local, and $\aff^{r|s}\times\aff^{p|q}=\aff^{p+r|q+s}$.
\end{proof}

\begin{Lem}
  Let $X$ be a super-premanifold over $S$. Then $X_0$ resp.~$X_\ev$ is a super-premanifold over $S_0$ resp. $S_\ev$. If $S=*$, then $X_0=X_\ev$.
\end{Lem}

\begin{proof}
  The statements are local, and for $\aff^{p|q}_S$, they follow from \thmref{Cor}{redbody-lim}.
\end{proof}

\begin{Prop}[relative-glue-smf][gluing of supermanifolds]
  Let $(U_i,X_i/U_i,\vphi_{ij})$ be gluing data for super-premanifolds over $S$. The glued superspace $X/S$ is a super-premanifold over $S$. Moreover, if $S$ and the $X_i$ are paracompact, then so is $X$.
\end{Prop}

\begin{proof}
  The existence of $X/S$ follows from \thmref{Prop}{relative-glue}. Since we have isomorphisms $X_i/U_i\to (X\times_SU_i)/U_i$, it is clear that $X$ possesses an open cover by local fibre charts, so that it is a super-premanifold over $S$.
  
  To establish the paracompactness assertion, we drop the subscripts $(-)_0$ and work in the category of topological spaces. Let $X$ be the space obtained by gluing the $X_i$, and $p:X\to S$ the projection. Let $\sh V=(V_j)$ be an open cover of $X$. Since $p$ is open and $S$ is paracompact, we may by passing to a refinement assume that $(U_i)$ is locally finite and each $p(V_j)$ is contained in some $U_i$. This gives open covers $\sh V_i$ of $X_i$ \scth $\sh V=\bigcup_i\sh V_i$. Replacing each of the $\sh V_i$ by a locally finite open refinement, we obtain a locally finite open refinement of $\sh V$. Hence, $X$ is paracompact.
\end{proof}

\begin{Prop}[base-change-smf][base change of relative supermanifolds]
  Let $X/S$ be a super-premanifold over $S$ and $T/S$ be in $\ssplfg{S}$. Then $X\times_ST$ is a super-premanifold over $T$, called the \Define{base change} of $X/S$. If $(x_a)$ is a local fibre coordinate system on $X/S$, then $(p_1^\sharp(x_a))$ is a local fibre coordinate system on $(X\times_ST)/T$. If $X/S$ and $T/S$ are supermanifolds over $S$, then $X\times_ST$ is a supermanifold over $T$ and over $S$.
\end{Prop}

\begin{proof}
  For the first statement, it suffices to observe $\aff^{p|q}_S\times_ST=\aff^{p|q}_T$. As to the second, we remark that $X_0\times_{S_0}T_0$ is Hausdorff in case $X_0$, $T_0$ are.
\end{proof}

\begin{Cor}[finprod-sman]
  Finite products in $\SSp_S^{\mathrm{lfg}}$ preserve $\SMan_S$.
\end{Cor}

\subsection{Weil functors}

In this subsection, we construct the Weil functors, and relative versions thereof. These furnish a uniform framework for the study of (relative) even and odd tangent bundles, and their higher order versions. Originally, Weil functors were defined by A.~Weil \cite{weil-pointproches} in the case of ordinary manifolds. For supermanifolds, they were introduced by Balduzzi--Carmeli--Fioresi \cite{balduzzi-carmeli-fioresi}. We show that they (and their relative versions) exist as relative fibre bundles and arise as (relative) \emph{inner homs} on the category $\ssplfg{\knums}$. This point of view is very efficient, as was originally noticed by Kontsevich for the odd tangent bundle in his seminal paper \cite{kontsevich}. 

In the following, we use some standard category theory facts and terminology extensively. We refer the reader to Ref.~\cite{maclane}.

\begin{Def}[inner-hom][Relative inner homs]
  Let $\cat C$ be a finitely complete category, $S$ be an object. We write $\cat C_S$ for the category of objects $X/S$, $Y/S$ over $S$ (\ie morphisms $p_X:X\to S$ and $p_Y:Y\to S$, respectively) and morphisms $\vphi:X/S\to Y/S$ over $S$ (\ie morphism $X\to Y$ in $\cat C$ that relate the given morphisms $p_X$ and $p_Y$). The hom sets in this category will be denoted by $\Hom[_S]0{X,Y}$.

  Given $X/S$ and $Y/S$, we define a functor $\GHom[_S]0{X,Y}$ on $\cat C_S$ with values in the category $\Sets$ of sets and maps by 
  \[
    \GHom[_S]0{X,Y}(T)\defi\Hom[_S]0{T\times_S X,Y}
  \]
  for all $T/S$. An object $Z$ of $\cat C_S$ representing $\GHom[_S]0{X,Y}$ is called a \Define{relative inner hom} from $X$ to $Y$ and denoted by $\GHom[_S]0{X,Y}$. That is, there is a natural bijection
  \[  
    \Hom[_S]0{T,\GHom[_S]0{X,Y}}=\Hom[_S]0{T\times_S X,Y}
  \]
  \fa $T/S$. If it exists, $\GHom[_S]0{X,Y}$ is unique up to canonical isomorphism. In case $S=*$ is a terminal object, we just write $\GHom0{X,Y}$ and drop the `relative'. Of course, for objects $X,Y$ of $\cat C_S$, $\GHom[_S]0{X,Y}$ is just $\GHom0{X,Y}$ in $\cat C_S$.
\end{Def}

\begin{Rem}
  This definition (for $S=*$) mimics the customary identification of maps of sets $f:T\times X\to Y$ and $g:T\to\Maps0{X,Y}$, given by $f(t,x)=g(t)(x)$.

  If $S$ is a set and $p_X:X\to S$, $p_Y:Y\to S$ are sets over $S$, then $\GHom[_S]0{X,Y}$ is represented by the disjoint union
  \[
    \coprod_{s\in S}\Parens1{\{s\}\times\Maps0{X_s,Y_s}},
  \]
  where $X_s=p_X^{-1}(s)$ and $Y_s=p_Y^{-1}(s)$ are the fibres of $X/S$ and $Y/S$, respectively. 
\end{Rem}

Relative inner homs are functorial in several ways, as we now discuss. 

\begin{Cons}[innerhom-functor][Functoriality of inner homs]
  Assume that $\GHom[_S]0{X,Y}$ is representable in the category $\cat C_S$. If $X'/S$ is \scth $\GHom[_S]0{X',Y}$ exists, consider $f:X/S\to X'/S$. We define a morphism 
  \[
    \GHom[_S]0{f,Y}:\GHom[_S]0{X',Y}\to\GHom[_S]0{X,Y}
  \]
  as follows: Let $\vphi\in\GHom[_S]0{X',Y}(T)=\Hom[_S]0{T\times_S X',Y}$. Define
  \[
    \GHom[_S]0{f,Y}(\vphi)\defi\vphi\circ({\id_T}\times_S f)\in\Hom[_S]0{T\times_S X,Y}=\GHom[_S]0{X,Y}(T).
  \]
  Since this construction is natural, it defines the required morphism $\GHom[_S]0{f,Y}$. 
  
  It is clear that $\GHom[_S]0{{\id_X},Y}={\id}$ and 
  \[
    \GHom[_S]0{f\circ f',Y}=\GHom[_S]0{f',Y}\circ\GHom[_S]0{f,Y},
  \]
  for $f:X'/S\to X''/S$ \scth $\GHom[_S]0{X'',Y}$ exists, so that $\GHom[_S]0{-,Y}$ is a contravariant functor on the full subcategory of $\cat C_S$ on which it is defined. 
  
  Similarly, for any $Y'/S$ \scth $\GHom[_S]0{X,Y'}$ is representable, we define, for $g:Y/S\to Y'/S$, a morphism $\GHom[_S]0{X,g}:\GHom[_S]0{X,Y}\to\GHom[_S]0{X,Y'}$ by setting, for any $\psi\in\GHom[_S]0{X,Y}(T)$
  \[
    \GHom[_S]0{X,g}(\psi)\defi g\circ\psi\in\Hom[_S]0{T\times_S X,Y'}=\GHom[_S]0{X,Y'}(T).
  \]
  One sees easily that $\GHom[_S]0{X,-}$ is a functor on its domain of existence in $\cat C_S$.
  
  Finally, assume that $\GHom[_S]0{X,Y}$, $\GHom[_S]0{Y,Z}$ are representable in $\cat C_S$. We define a morphism
  \[
    \circ=\circ_{X,Y,Z}:\GHom[_S]0{Y,Z}\times\GHom[_S]0{X,Y}\to\GHom[_S]0{X,Z},
  \]
  called \emph{inner circle}, as follows. For 
  \begin{align*}
    \vphi\in\GHom[_S]0{X,Y}(T)&=\Hom[_T]0{T\times_S X,T\times_S Y},\\
    \psi\in\GHom[_S]0{Y,Z}(T)&=\Hom[_T]0{T\times_S Y,T\times_S Z},
  \end{align*}
  we obtain a composite $\psi\circ\vphi\in\GHom[_S]0{X,Z}(T)$ as 
  \[
    \psi\circ_\cat C\vphi\in\Hom[_T]0{T\times_S X,T\times_S Z},
  \]
  after identifying with $\GHom[_S]0{X,Z}(T)$. The inner circle $\circ$ is a natural transformation of (tri-)functors. It enjoys the obvious associative law. 

  Moreover, if $\GHom[_S]0{X,X}$ is representable, there is a morphism
  \[
    1_X:S\to\GHom[_S]0{X,X}
  \]
  over $S$, corresponding to the element ${\id_X}\in\Hom[_S]0{X,X}=\GHom[_S]0{X,X}(S)$. (Recall that $S$ is the tensor unit of $\times_S$.) It is a left and right identity for the inner circle. In particular, $\GHom[_S]0{X,X}$ is a monoid object in $\cat C_S$, whenever it exists. 
\end{Cons}

Weil functors are (relative) fibre bundles. In general, the fibres are functors (at least if $\knums\neq\Bbbk$). We briefly digress to give the relevant definitions. 

\begin{Def}[cov-func][coverings of functors]
  Let $X,Y$ be set-valued cofunctors on $\ssplfg{\knums}$. A morphism $j:Y\to X$ is called an \Define[functor!open embedding]{open embedding} if the following two conditions are verified:
  \begin{enumerate}[wide]
    \item $j$ is injective, \ie $j_T:Y(T)\to X(T)$ is injective for all $T\in\ssplfg{\knums}$, and
    \item for any morphism $\vphi:T\to X$, where $T\in\ssplfg{\knums}$, the induced morphism $j\times_X\vphi$, given by the following Cartesian diagram:
    \begin{center}
      \begin{tikzcd}
        Y\times_XT\dar{}\rar{j\times_X\vphi}&T\dar{\vphi}\\
        Y\rar{j}&X
      \end{tikzcd}
    \end{center}
    is representable by an open embedding in $\ssplfg{\knums}$. That is, we have an isomorphism $Y\times_XT\cong T'$ \fs $T'\in\ssplfg{\knums}$ and \via this isomorphism, $j\times_X\vphi$ corresponds to an open embedding $j':T'\to T$ of superspaces.
  \end{enumerate}

  If $Y$ is a \Define{subfunctor}, \ie $Y(T)\subseteq X(T)$ \fa $T\in\ssplfg{\knums}$, then it is called an \Define{open subfunctor} if the inclusion morphism $Y\to X$ is an open embedding. 

  A collection $(\vphi_i:X_i\to X)$ of open embeddings (of functors) is called a \Define[functor!covering]{covering}\index{covering!of a functor} if for any morphism $T\to X$, the open subspaces $T_i$ representing $X_i\times_XT$ form a covering of the superspace $T$. 
\end{Def}

\begin{Ex}
  Let $X$ be a set-valued cofunctor on $\ssplfg{\knums}$. Then $X$ is an open subfunctor of $X$. Indeed, let $\vphi:T\to X$ be a morphism, where $T\in\ssplfg{\knums}$. Then we have an isomorphism $X\times_XT\cong T$, and ${\id}_X\times_X\vphi$ corresponds under this isomorphism to the identity $\id_T$, which is an open embedding. 

  A morphism $Y\to X$ of $X,Y\in\ssplfg{\knums}$ defines an open embedding of functors if and only if it is an open embedding of superspaces. Indeed, open embeddings are stable under fibre products; conversely, we may consider $T=X$ and $\vphi=\id_X$. 
\end{Ex}

\begin{Def}[rel-func][functors and sheaves over superspaces]
  Let $S\in\ssplfg{\knums}$ and $X$ be a set-valued cofunctor on the category $\ssplfg{\knums}$. We say that $X$ is a functor \Define{over $S$} and write $X/S$ if we are given a morphism $X\to S$ in the functor category $\smash{\Bracks1{\Parens1{\ssplfg{\knums}}^{\mathrm{op}},\Sets}}$. Here, we identify $S$ with its point functor. Any superspace $X/S\in\ssplfg{S}$ is a functor over $S$. Similarly, one defines functors \Define{under $S$}.

  A functor $X$ will be called a \Define{sheaf} if it is a sheaf on the site given by the coverings of superspaces, see Refs.~\cite{ks}*{Definitions 17.2.1, 17.3.1}, \cite{moerdijk-reyes}*{Appendix 1, Definitions 1.3, 2.4}. If $X$ is representable, then it is a sheaf. A functor over $S$ that is a sheaf will be called a \Define{sheaf over $S$}. The category of sheaves over $S$ and their morphisms will be denoted by $\Shv_S$.
\end{Def}

We can define gluing data for sheaves over a superspace $S$.

\begin{Def}[fun-glue][gluing data for sheaves over $S$]
  Let $\vphi=(\vphi_i:U_i\to S)$ be a covering of $S\in\ssplfg{\knums}$, and assume given sheaves $X_i/U_i$ over $U_i$ and isomorphisms $\psi_{ij}:X_j\times_{U_j}U_{ij}\to X_i\times_{U_i}U_{ij}$. Then the collection $(\vphi_i:U_i\to S,X_i/U_i,\psi_{ij})$ will be called \Define{gluing data for a sheaf over $S$} if
    \begin{equation}\label{eq:fun-gluing}
      p_{12}^*(\psi_{ij})\circ p_{23}^*(\psi_{jk})\circ p_{13}^*(\psi_{ki})={\id_{X_i\times_{U_i}U_{ijk}}}.
    \end{equation}
    Here, $p_{12}^*{\psi_{ij}}$ is the morphism $X_j\times_{U_j}U_{ijk}\to X_i\times_{U_i}U_{ijk}$ induced by $\psi_{ij}$ by pulling back along the left-hand vertical face of the pullback square
    \begin{center}
      \begin{tikzcd}
        U_{ijk}\rar{p_3}\dar[swap]{p_{12}} & U_k\dar{\vphi_k}\\
        U_{ij}\rar[swap]{\vphi_{ij}} & S
      \end{tikzcd}
    \end{center}
    Similarly for the other quantites in Equation~\eqref{eq:fun-gluing}.
\end{Def}

In this context, \thmref{Prop}{relative-glue} generalises easily. 

\begin{Prop}[fun-glue]
  Let $S$ be a superspace and $(\vphi_i:U_i\to S,X_i/U_i,\psi_{ij})$ gluing data for a sheaf over $S$. Then there exist a sheaf $X/S$ over $S$ and isomorphisms $\psi_i:X_i/U_i\to (X\times_SU_i)/U_i$ with $\psi_j=\psi_i\circ\psi_{ij}$ on $X_j\times_{U_j}U_{ij}$.

  These data are uniquely characterised up to unique isomorphism by the following universal property: For any sheaf $X'/S$ over $S$ and any choice of morphisms $\vrho_i:X_i/U_i\to (X'\times_SU_i)/U_i$ \scth $\vrho_j=\vrho_i\circ\psi_{ij}$ on $X_j\times_{U_j}U_{ij}$, there exists a unique morphism $\vrho:X/S\to X'/S$ \scth $\vrho\circ\psi_i=\vrho_i$.
\end{Prop}

\begin{Def}[fibrebundle-def][relative fibre bundles]
  Let $E/S$ be a sheaf over $S$ (\eg a locally finitely generated superspace over $S$) and let $X/S$ be in $\ssplfg{S}$. Assume given a morphism $p:E/S\to X/S$ over $S$. A (relative) \Define{local trivialisation} $(\vphi:U\to X,F,\tau)$ of $E/X$ (or of $p$) with fibre $F$ consists of an open embedding $\vphi:U\to X$, a sheaf $F/S$ over $S$, and an isomorphism
  \[
    \tau:U\times_SF\to U\times_XE
  \]
  over $U$. A family $(\vphi_i:U_i\to X,F,\tau_i)_{i\in I}$ of local trivialisations is called a \Define{trivialising covering} of $E/X$ if $(\tau_i)_{i\in I}$ is a covering of $E$, see \thmref{Def}{gluing-data} for the representable case and \thmref{Def}{cov-func} for the general case.

  We say that $E/X$ is a \Define{relative fibre bundle} (over $S$) if $E$ is a sheaf and there exists a trivialising covering of $E/X$. The condition that $E$ is a sheaf is automatically verified in case $E$ is representable. 

  In case $S=*$, we simply say that $E/X$ is a \Define{fibre bundle}. If $X/S$ is in $\SMan_S$ and $E/X$ is a relative fibre bundle with a trivialising covering whose fibres are supermanifolds over $S$, then $E$ is a supermanifold over $S$, according to \thmref{Cor}{finprod-sman}.
\end{Def}

In the following, let $A=\knums\oplus\ger m$ be a Weil $\knums$-superalgebra. Moreover, we fix $S\in\ssplfg{\knums}$. We will write
\[
  \Spec_S A\defi S^A=S\times\Spec A.
\]
Given graded dimensions $p|q$ and $r|s$, we let $(p|q)(r|s)\defi (pr+qs)|(ps+qr)$ and say that this is \Define{purely odd} if $pr+qs=0$. 

For any sheaf $X/S$ over $S$, we define the relative inner hom functor
\begin{equation}\label{eq:weilfun-def}  
  T_S^AX\defi\GHom[_S]0{\Spec_S A,X}.
\end{equation}
We call $T_S^A:\Shv_S\to\Shv_S$ the \Define{relative $A$-Weil functor}.

\begin{Prop}[weil-innerhom]
  Let $X/S$ be a relative supermanifold. If $\knums=\Bbbk$ or the graded dimension $(\dim_SX)(\dim\ger m)$ is purely odd, then $T^A_SX$ is representable in the category $\ssplfg{S}$, and lies in the subcategory $\SMan_S$. In general, we have the following facts.
  \begin{enumerate}[wide]
    \item\label{wim-i} For $A=\knums$, we have a natural isomorphism $T_S^\knums X\to X$ of sheaves over $S$.
    \item\label{wim-ii} For any morphism $\vphi:A\to B$ of Weil $\knums$-superalgebras, there is a morphism
    \[
      T_S^\vphi:T_S^AX\to T_S^BX,
    \]
    of sheaves over $S$, and the assignment $\vphi\mapsto T_S^\vphi$ is functorial.
    \item\label{wim-iii} For any $\vphi:A\to B$ as in \eqref{wim-ii}, the $T_S^\vphi$ are natural transformations of functors $T_S^A\to T_S^B$. That is, there is a commutative diagram
    \begin{center}
      \begin{tikzcd}
        T_S^AX\dar[swap]{T_S^\vphi}\rar{T_S^A\psi}&T_S^AY\dar{T_S^\vphi}\\
        T_S^BX\rar[swap]{T_S^B\psi}&T_S^BY
      \end{tikzcd}
    \end{center}
    for every morphism $\psi:X/S\to Y/S$. Here, $T_S^A\psi\defi\GHom[_S]0{\Spec_SA,\psi}$.
    \item\label{wim-iv} The unique morphisms
    \begin{center}
      \begin{tikzcd}
        \knums\rar{\eta}&A\rar{\eps}&\knums
      \end{tikzcd}
    \end{center}
    of $\knums$-superalgebras induce \via \eqref{wim-ii} canonical morphisms
    \begin{center}
      \begin{tikzcd}[column sep=large]
        X\rar{s^A_X\defi T_S^\eta}&T_S^AX\rar{p^A_X\defi T_S^\eps}&X
      \end{tikzcd}
    \end{center}
    of sheaves over $S$. Then $s_X^A$ is a section of $p_X^A$, \ie $p_X^A\circ s_X^A={\id_X}$.
    \item\label{wim-vi} We may consider the sheaf $T_S^AX$ as lying over $X$ by virtue of the projection $p^A_X:T^A_SX\to X$ and under $X$ by the canonical section $s^A_X:X\to T_S^AX$. Then the morphisms $T_S^A\psi:T_S^AX\to T_S^AY$ are over and under $\psi$. The projection and canonical section define natural transformations $p_X:\smash{T_S^{(-)}}\to\id$ and $s_X:\id\to\smash{T_S^{(-)}}$, so the morphisms $\smash{T_S^\vphi}$ from \eqref{wim-ii} are over and under $X$.
    \item\label{wim-v} With the morphism $p^A_X:T_S^AX\to X$ as the projection, $T_S^AX/X$ is a relative fibre bundle, whose fibre at $x\in X_0$ with $p|q=\dim_{S,x}X$ is isomorphic to the sheaf $\smash{\aff^\knums(\ger m^{p|q})}=\smash{\aff^\knums(\ger m\otimes_\knums\knums^{p|q})}$. Thus, if $\knums=\Bbbk$ and $\dim A=(r+1)|s$, then
    \[
      \dim_{X,y}T_S^AX=pr+qs|ps+qr\mathtxt{whenever}\dim_{S,x}X=p|q\,,\,p_X^A(y)=x.
    \]
  \end{enumerate}
\end{Prop}

\begin{Rem}
  As stated above, the Weil functors $T_S^AX$ are not representable in $\ssplfg{\knums}$ if $\knums\neq\Bbbk$ unless $(\dim_SX)(\dim\ger m)$ is purely odd. In fact, this is true in general of (relative) fibre bundles with a $\knums$-affine fibre that is not purely odd. This phenomenon is reflected by the following fact: Let $X$ be the real supermanifold associated with a not purely odd \emph{cs} manifold. Then the complex structure on the odd part of $TX$ defines a non-involutive distribution, as observed in Ref.~\cite{bgls}.
\end{Rem}

In a first step, we construct the various natural morphisms presented above by use of the adjunction defining the inner hom functor.

\begin{Cons}[weil-extend]
  Consider a morphism of Weil $\knums$-superalgebras, say $\vphi:A\to B$. This induces a morphism
  \[
    \Spec\vphi\defi(*,\vphi):\Spec B\to\Spec A,
  \]
  which we promote to a morphism
  \[
    \Spec_S\vphi:\Spec_SB\to\Spec_SA
  \]
  over $S$. We may define
  \[
    T_S^\vphi\defi\GHom[_S]0{\Spec_S\vphi,X}:T_S^AX\to T_S^BX
  \]
  by simply applying \thmref{Cons}{innerhom-functor}. If $\GHom[_S]0{\Spec_SA,X}$ and $\GHom[_S]0{\Spec_SB,X}$ are representable, then it is a morphism of superspaces over $S$, by the Yoneda lemma.
\end{Cons}

\begin{Lem}[innerhom-sheaf]
  Let $X/S$ and $Y/S$ be in $\ssplfg{S}$. Then the inner hom functor $\GHom[_S]0{X,Y}$ is a sheaf over $S$.
\end{Lem}

\begin{proof}
  The statement is immediate from the definition, combined with \thmref{Prop}{open-fibreprod}, \thmref{Prop}{relative-glue}, and \thmref{Cor}{emb-fibprod}.
\end{proof}

\begin{proof}[\protect{Proof of \thmref{Prop}{weil-innerhom}}]
  We first assume that $X$ is an open subspace of the relative affine superspace $\smash{\aff_S^{p|q}=S\times\aff^{p|q}}$. Consider the tensor product $\ger m^{p|q}=\ger m\otimes_\knums\knums^{p|q}$.

  Recall the facts discussed in \thmref{Scho}{mor-weil} and \thmref{Scho}{mor-weil-k}. For $T\in\ssplfg{S}$, we have, by the obvious extension of Equation \eqref{eq:weil-ext-mor2} to the relative case, a bijection
  \[
    \Hom[_S]0{T\times_S\Spec_S A,X}=\Hom[_S]0{T\times\Spec A,X}\cong\Hom[_S]0{T,X\times\aff^\knums(\ger m^{p|q})},
  \]
  which is natural in $T/S$. Thus, $\GHom[_S]0{\Spec_S A,X}$ is isomorphic to the direct product $X\times\aff^\knums(\ger m^{p|q})$. 
  
  The general case now follows by passing to an open cover, using the functoriality of $\GHom[_S]0{\Spec_S A,-}$ from \thmref{Cons}{innerhom-functor}, and applying \thmref{Prop}{relative-glue} and \thmref{Prop}{fun-glue}, respectively. By the same token, the assignment $X\mapsto T_S^AX$ defines a functor.

  For $A=\knums$, we note that $\Spec_SA=S$, so we have natural bijections
  \[
    \GHom[_S]0{\Spec_SA,X}(T)=\Hom[_S]0{T\times_SS,X}=\Hom[_S]0{T,X}
  \]
  for $T/S$, so that $X\cong T_S^\knums X$ by the Yoneda Lemma. This proves item \eqref{wim-i}.

  Item \eqref{wim-ii} follows immediately by \thmref{Cons}{weil-extend}. Similarly, \eqref{wim-iii} is implied by the associativity of the inner circle, \vq \thmref{Cons}{innerhom-functor}. Item \eqref{wim-iv} is a special case of \eqref{wim-ii}. For \eqref{wim-vi}, the first part is a special case of \eqref{wim-iii}, and the second follows from 
  \[
    \eps_B\circ\vphi=\eps_A\nd\eta_B=\vphi\circ\eta_A
  \]
  for any morphism $\vphi:A\to B$ of Weil $\knums$-superalgebras. The second equation is a consequence of unitality, and the first is one of locality, which is automatic. (For the latter, one may argue as in the proof of \thmref{Prop}{mor-local}.)

  As for \eqref{wim-v}, it is by the functoriality and \eqref{wim-vi} sufficient to consider an open subspace $\smash{X\subseteq\aff^{p|q}_S}$, which was treated in detail at the beginning of the proof. 
\end{proof}

\begin{Def}[weil-functor][relative Weil functor]
  The functor $T_S^A$ which maps $X/S\mapsto T^A_SX$ and morphisms $\psi:X/S\to Y/S$ to morphisms $\psi^A\defi T_S^A\psi:T_S^AX\to T_S^AY$ over $\psi$ is called the (relative) \Define[Weil functor]{$A$-Weil functor}. Endowed with the projection $p^A_X$, $T_S^AX$ is called the (relative) \Define{$A$-Weil bundle} of $X$.
\end{Def}

\begin{Ex}
  Let $X\in\SMan_S$ and consider the algebra $A\defi\sdual_\ev$ of dual numbers. There is a natural bijection from $\Hom[_S]0{T,T_S^AX}$ to the set of pairs
  \[
    \Set0{(\vphi,\delta)\in\Hom[_S]0{T,X}\times\Hom[_{\smash{p_{T,0}^{-1}\sh O_S}}]0{\vphi_0^{-1}\sh O_X,\sh O_T}\,}{\,\delta(fg)=\vphi^\sharp(f)\delta(g)+\delta(f)\vphi^\sharp(g)}
  \]
  of morphisms $\vphi:T/S\to X/S$ and even derivations along $\vphi$, which are linear over $\sh O_S$. One easily verifies that the transition functions of $T_S^AX$ are given by first derivatives. Hence, $T_S^AX$ is a $\knums$-vector bundle whose sheaf of sections is seen to be 
  \[
    \sh T_{X/S}\defi\ShGDer[_{p_{X,0}^{-1}\sh O_S}]0{\sh O_X,\sh O_X},
  \]
  the \emph{relative tangent sheaf} of $X$, where $\underline{\sh Der}$ is the sheaf of graded derivations. Moreover, if $E$ is a $\knums$-vector bundle over $X$, then by definition, its sheaf of sections is the unique $\sh O_X$-module sheaf $\sh E$ such that for every $Y/X$, we have 
  \[
    \Hom0{Y,Y\times_XE}=\Hom0{\sh O_Y,p_Y^*(\sh E)}.
  \]
  Therefore, $T_S^AX$ is the relative tangent bundle of $X$.
\end{Ex}

\begin{Prop}[weilfunc-basic1]
  The relative Weil functor $T_S^A$ enjoys the following properties:
  \begin{enumerate}[wide]
    \item\label{wf-i} $T_S^A$ preserves finite fibre products over $S$, and there is a natural base change 
    \[
      R\times_ST_S^AX\cong T_R^A(R\times_SX).
    \]
    \item\label{wf-ii} If $U$ is an open subspace of $X$, then $T_S^AU$ is an open subfunctor of $T_S^AX$.
    \item\label{wf-iii} If $B$ is a Weil $\knums$-superalgebra, then there is a natural isomorphism 
    \[
      T_S^{A\otimes B}=T_S^A\circ T_S^B.
    \]
    In particular, we have a natural isomorphism $T_S^A\circ T_S^B\cong T_S^B\circ T_S^A$. 
  \end{enumerate}
\end{Prop}

\begin{proof}
  \eqref{wf-i} Preservation of binary products follows from the natural bijection
  \begin{align*}
    \Hom[_S]1{T,T_S^A(X\times_SY)}&=\Hom[_S]1{T\times_S\Spec_SA,X\times_SY}\\
    &=\Hom[_S]1{T\times_S\Spec_SA,X}\times\Hom[_S]1{T\times_S\Spec_SA,Y}\\
    &=\Hom[_S]1{T,T_S^AX}\times\Hom[_S]1{T,T_S^AY}   
   \end{align*} 
  after an application of Yoneda's Lemma. The case of the empty product is similar:
  \[
    \Hom[_S]1{T,T_S^A(S)}=\Hom[_S]1{T\times_S\Spec_S A,S}=*,
  \]
  so that $T_S^A(S)\cong S$. For the base change, we compute
  \begin{align*}
    \Hom[_R]1{T,R\times_ST_S^AX}&=\Hom[_S]1{T,T_S^AX}=\Hom[_S]1{T\times_S\Spec_SA,X}\\
    &=\Hom[_R]1{T\times_R\Spec_RA,R\times_SX}\\
    &=\Hom[_R]1{T,T_R^A(R\times_SX)},
  \end{align*}
  where we recall $\Spec_SA=S\times\Spec A$, and the superspace $T/R$ is considered as lying over $S$ by the composition with $R/S$.
  
  Item \eqref{wf-ii} is clear by construction, and for \eqref{wf-iii}, we observe that 
  \[
    \Spec A\otimes B=\Spec A\times\Spec B
  \]
  by \thmref{Prop}{weil-extended}, so 
  \begin{align*}
    \Hom[_S]1{T,T_S^{A\otimes B}X}&=\Hom1{T\times_S\Spec_S A\times_S\Spec_S B,X}\\
    &=\Hom1{T\times_S\Spec_S A,T_S^BX}=\Hom[_S]0{T,T_S^A(T_S^BX)},
  \end{align*}
  for any locally finitely generated $\knums$-superspace $T/S$.
\end{proof}

The following lemma is a generalisation of \thmref{Prop}{weil-innerhom} \eqref{wim-v}. It shows the utility of the general point of view on Weil functors, but we will not be applying it.

\begin{Lem}
  Let $\smash{\vphi:A\to B}$ be a morphism of Weil $\knums$-superalgebras. If $\vphi$ is surjective, then $\smash{T_S^\vphi X:T_S^AX\to T_S^BX}$ is a relative fibre bundle with fibre at $x$ isomorphic to $\aff^\knums((\ker\vphi)^{p|q})$, where $\dim_{S,x}X=p|q$. If $\vphi$ is injective, then locally, $T_S^\vphi$ is isomorphic to $\smash{({\id},0):T_S^AX\to T_S^AX\times\aff^\knums(\ger m_B/\vphi(\ger m_A))^{p|q}}$.
\end{Lem}

\begin{proof}
  Both statements are local. Let $X$ be an open subspace of $\aff_S^{p|q}$. Then $\smash{T_S^\vphi:X\times\aff^\knums(\ger m_A^{p|q})\to X\times\aff^\knums(\ger m_B^{p|q}})$ is given by $\smash{T_S^A={\id_X}\times\aff^\knums(\vphi^{p|q}})$, where we define $\smash{\vphi^{p|q}\defi\vphi|_{\ger m_A}\otimes{\id_{\knums^{p|q}}}}$. If $\vphi$ is injective, then there is a linear isomorphism $\ger m_B\cong\ger m_A\times\big(\ger m_B/\vphi(\ger m_A)\big)$. Similarly, if it is surjective, then there is a linear isomorphism $\ger m_A\to\ker\vphi\times\ger m_B$.
\end{proof}

Remarkably, the natural transformations $T^\vphi_S$ already exhaust the full set of natural transformations between (relative) Weil functors. In the setting of manifolds, this was already noticed by A.~Weil \cite{weil-pointproches}.

\begin{Prop}[weil-nat]
  Let $A$ and $B$ be Weil $\knums$-superalgebras. Consider $T_S^A$ and $T_S^B$ as functors $\SMan_S\to\Shv_S$. The map $\vphi\mapsto T_S^\vphi$ is a bijection between algebra morphisms $A\to B$ to transformations $T_S^A\to T_S^B$ natural under base change in $T/S$.
\end{Prop}

In the \emph{proof}, we note the following lemma.

\begin{Lem}[weil-nat-unique]
  Any natural transformation $T_S^A\to T_S^B$ is uniquely determined by its value at $X=\aff_S^1$.
\end{Lem}

\begin{proof}
  Indeed, let $\sigma,\tau:T_S^A\to T_S^B$ be natural. We first show that they are determined by their values on $\smash{\aff_S^{1|1}}$. So, assume that $\sigma_{\smash{\aff_S^{1|1}}}=\tau_{\smash{\aff_S^{1|1}}}$.

  The naturality of $\sigma$ and $\tau$ gives $\sigma_X=\tau_X$ for $X=\smash{\aff^1_S}$ and $X=\smash{\aff_S^{0|1}}$. Since $\smash{T_S^A}$ and $\smash{T_S^B}$ preserve fibre products over $S$, we obtain equality for $X=\smash{\aff_S^{p|q}}$, where $p$ and $q$ are arbitrary. By item \eqref{wf-ii} of \thmref{Prop}{weilfunc-basic1}, this equality extends to the case of an open subspace $X$ of $\smash{\aff_S^{p|q}}$, and thus, to an arbitrary supermanifold over $S$, in view of \thmref{Prop}{relative-glue}.

  This proves $\sigma=\tau$ under the assumption of $\sigma_X=\tau_X$ for $\smash{X=\aff_S^{1|1}}$. Naturality and item \eqref{wf-iii} of \thmref{Prop}{weilfunc-basic1} show that $\smash{T_S^C\tau_X=\tau_{T_S^CX}}$ for any Weil superalgebra $C$. In particular, for $\smash{C=\knums[\tau]}$ ($\tau$ being an odd indeterminate), we have $\smash{T_S^C\aff^{1|0}=\aff_S^{1|1}}$, by Equation \eqref{eq:weil-ext-mor2}. Then we have $\tau_{\smash{\aff_S^{1|1}}}=T_S^C\tau_{\smash{\aff_S^{1|0}}}$, finally proving the assertion.
\end{proof}

\begin{proof}[\protect{Proof of \thmref{Prop}{weil-nat}}]
  Observe that $T^A(\aff^1)=\aff^\knums(A)$. Thus, the structure morphisms $a,m:\aff^1\times\aff^1\to\aff^1$ of addition and multiplication, respectively, of the algebra $\knums$, give rise to operations $T^A(a)$ and $T^A(m)$ on $A$. By construction of $T^A$, these are the same as $\aff^\knums(-)$, applied to the given addition and multiplication of $A$.

  Using the base change formula in \thmref{Prop}{weilfunc-basic1}, we have the identities $T_S^A(\aff_S^1)=S\times T^A(\aff^1)$, $T_S^A(a)={\id_S}\times T^A(a)$, and $T_S^A(m)={\id_S}\times T^A(m)$. Any natural transformation $T_S^A\to T_S^B$ defines a morphism $\psi:S\times\aff^\knums(A)\to S\times\aff^\knums(B)$ over $S$, which commutes with the promotion of $a$ and $m$ to operations on $S\times\aff^\knums(A)$ and $S\times\aff^\knums(B)$. The unitality of $\psi$ is obtained by applying naturality to 
  \[
    \eta_S\defi{\id}_S\times\aff^\knums(\eta):S\times\aff^\knums(\knums)\to S\times\aff^\knums(A),
  \]
  where $\eta(1)=1_A$ is the unit of $A$. Again by naturality, $\psi$ is of the form ${\id_S}\times\aff^\knums(\vphi)$ for some $\knums$-algebra morphism $\vphi:A\to B$. By \thmref{Lem}{weil-nat-unique}, this proves our claim.
\end{proof}

\end{document}